\documentclass[11pt]{amsart}
\usepackage{amssymb}
\usepackage{amsmath, amsfonts, amsthm, amssymb, cite}
\usepackage{color}
\usepackage[all]{xy}
\usepackage{fancyhdr}
\begin{document}



\newtheorem{theorem}{Theorem}[section]
\newcommand{\mar}[1]{{\marginpar{\textsf{#1}}}}
\newtheorem*{theorem*}{Theorem}
\newtheorem{prop}[theorem]{Proposition}
\newtheorem*{prop*}{Proposition}
\newtheorem{lemma}[theorem]{Lemma}
\newtheorem{corollary}[theorem]{Corollary}
\newtheorem*{conj*}{Conjecture}
\newtheorem*{qst}{Question}
\newtheorem*{corollary*}{Corollary}
\newtheorem{definition}[theorem]{Definition}
\newtheorem*{definition*}{Definition}
\newtheorem{remarks}[theorem]{Remarks}
\newtheorem{example}[theorem]{Example}
\newtheorem{rems}[theorem]{Remarks}
\newtheorem{rem}[theorem]{Remark}
\newtheorem*{rem*}{Remark}
\newtheorem*{rems*}{Remarks}
\newtheorem{cor}[theorem]{Corollary}
\newtheorem{defn}[theorem]{Definition}
\newtheorem{thm}[theorem]{Theorem}

\newtheorem*{not*}{Notation}
\newcommand\pa{\partial}
\newcommand\cohom{\operatorname{H}}
\newcommand\Td{\operatorname{Td}}
\newcommand\Trig{\operatorname{Trig}}
\newcommand\Hom{\operatorname{Hom}}
\newcommand\End{\operatorname{End}}
\newcommand\Ker{\operatorname{Ker}}
\newcommand\Ind{\operatorname{Ind}}
\newcommand\cker{\operatorname{coker}}
\newcommand\oH{\operatorname{H}}
\newcommand\oK{\operatorname{K}}
\newcommand\codim{\operatorname{codim}}
\newcommand\Exp{\operatorname{Exp}}
\newcommand\CAP{{\mathcal AP}}
\newcommand\T{\mathbb T}
\newcommand{\M}{\mathcal{M}}
\newcommand\ep{\epsilon}
\newcommand\te{\tilde e}
\newcommand\Dd{{\mathcal D}}

\newcommand\what{\widehat}
\newcommand\wtit{\widetilde}
\newcommand\mfS{{\mathfrak S}}
\newcommand\cA{{\mathcal A}}
\newcommand\maA{{\mathcal A}}
\newcommand\maF{{\mathcal F}}
\newcommand\maN{{\mathcal N}}
\newcommand\cM{{\mathcal M}}
\newcommand\maE{{\mathcal E}}
\newcommand\cF{{\mathcal F}}
\newcommand\maG{{\mathcal G}}
\newcommand\cG{{\mathcal G}}
\newcommand\cH{{\mathcal H}}
\newcommand\maH{{\mathcal H}}
\renewcommand\H{{\mathcal H}}
\newcommand\cO{{\mathcal O}}
\newcommand\cR{{\mathcal R}}
\newcommand\cS{{\mathcal S}}
\newcommand\cU{{\mathcal U}}
\newcommand\cV{{\mathcal V}}
\newcommand\cX{{\mathcal X}}
\newcommand\cD{{\mathcal D}}
\newcommand\cnn{{\mathcal N}}
\newcommand\wD{\widetilde{D}}
\newcommand\wL{\widetilde{L}}
\newcommand\wM{\widetilde{M}}
\newcommand\wV{\widetilde{V}}
\newcommand\Ee{{\mathcal E}}
\newcommand{\npartial}{\slash\!\!\!\partial}
\newcommand{\Heis}{\operatorname{Heis}}
\newcommand{\Solv}{\operatorname{Solv}}
\newcommand{\Spin}{\operatorname{Spin}}
\newcommand{\SO}{\operatorname{SO}}
\newcommand{\ind}{\operatorname{ind}}
\newcommand{\Index}{\operatorname{Index}}
\newcommand{\ch}{\operatorname{ch}}
\newcommand{\rank}{\operatorname{rank}}
\newcommand{\G}{\Gamma}
\newcommand{\HK}{\operatorname{HK}}
\newcommand{\Dix}{\operatorname{Dix}}

\newcommand{\tM}{\tilde{M}}  
\newcommand{\tS}{\tilde{S}}
\newcommand{\tH}{\tilde{\mathcal H}}
\newcommand{\tg}{\tilde{g}}
\newcommand{\tx}{\tilde{x}}
\newcommand{\ty}{\tilde{y}}
\newcommand{\ox}{\otimes}


\newcommand{\abs}[1]{\lvert#1\rvert}
 \newcommand{\A}{{\mathcal A}}
        \newcommand{\D}{{\mathcal D}}
\newcommand{\HH}{{\mathcal H}}
        \newcommand{\LL}{{\mathcal L}}
        \newcommand{\B}{{\mathcal B}}
        \newcommand{\K}{{\mathcal K}}
\newcommand{\oo}{{\mathcal O}}
         \newcommand{\PP}{{\mathcal P}}
         \newcommand{\Q}{{\mathcal Q}}
        \newcommand{\s}{\sigma}
        \newcommand{\coker}{{\mbox coker}}
        \newcommand{\dd}{|\D|}
        \newcommand{\n}{\parallel}
\newcommand{\bma}{\left(\begin{array}{cc}}
\newcommand{\ema}{\end{array}\right)}

\newcommand{\bca}{\left(\begin{array}{c}}
\newcommand{\eca}{\end{array}\right)}
\newcommand{\sr}{\stackrel}
\newcommand{\da}{\downarrow}
\newcommand{\tD}{\tilde{\D}}
        \newcommand{\R}{\mathbb R}
        \newcommand{\C}{\mathbb C}
        \newcommand{\h}{\mathbb H}
\newcommand{\Z}{\mathcal Z}
\newcommand{\N}{\mathbb N}
\newcommand{\tto}{\longrightarrow}
\newcommand{\ZZ}{{\mathcal Z}}
\newcommand{\ben}{\begin{displaymath}}
        \newcommand{\een}{\end{displaymath}}
\newcommand{\be}{\begin{equation}}
\newcommand{\ee}{\end{equation}}
        \newcommand{\bean}{\begin{eqnarray*}}
        \newcommand{\eean}{\end{eqnarray*}}
\newcommand{\nno}{\nonumber\\}
\newcommand{\bea}{\begin{eqnarray}}
        \newcommand{\eea}{\end{eqnarray}}
\newcommand{\x}{\times}

\newcommand{\Ga}{\Gamma}
\newcommand{\e}{\epsilon}
\renewcommand{\L}{\mathcal{L}}
\newcommand{\supp}[1]{\operatorname{#1}}
\newcommand{\norm}[1]{\parallel\, #1\, \parallel}
\newcommand{\ip}[2]{\langle #1,#2\rangle}
\newcommand{\nc}{\newcommand}
\nc{\gf}[2]{\genfrac{}{}{0pt}{}{#1}{#2}}
\nc{\mb}[1]{{\mbox{$ #1 $}}}
\nc{\real}{{\mathbb R}}
\nc{\comp}{{\mathbb C}}
\nc{\ints}{{\mathbb Z}}
\nc{\Ltoo}{\mb{L^2({\mathbf H})}}
\nc{\rtoo}{\mb{{\mathbf R}^2}}
\nc{\slr}{{\mathbf {SL}}(2,\real)}
\nc{\slz}{{\mathbf {SL}}(2,\ints)}
\nc{\su}{{\mathbf {SU}}(1,1)}
\nc{\so}{{\mathbf {SO}}}
\nc{\hyp}{{\mathbb H}}
\nc{\disc}{{\mathbf D}}
\nc{\torus}{{\mathbb T}}
\newcommand{\tk}{\widetilde{K}}
\newcommand{\boe}{{\bf e}}\newcommand{\bt}{{\bf t}}
\newcommand{\vth}{\vartheta}
\newcommand{\CGh}{\widetilde{\CG}}
\newcommand{\db}{\overline{\partial}}
\newcommand{\tE}{\widetilde{E}}
\newcommand{\tr}{{\rm tr}}
\newcommand{\ta}{\widetilde{\alpha}}
\newcommand{\tb}{\widetilde{\beta}}
\newcommand{\txi}{\widetilde{\xi}}
\newcommand{\hV}{\hat{V}}
\newcommand{\IC}{\mathbf{C}}
\newcommand{\IZ}{\mathbf{Z}}
\newcommand{\IP}{\mathbf{P}}
\newcommand{\IR}{\mathbf{R}}
\newcommand{\IH}{\mathbf{H}}
\newcommand{\IG}{\mathbf{G}}
\newcommand{\IS}{\mathbf{S}}
\newcommand{\CC}{{\mathcal C}}
\newcommand{\CD}{{\mathcal D}}
\newcommand{\CS}{{\mathcal S}}
\newcommand{\CG}{{\mathcal G}}
\newcommand{\CL}{{\mathcal L}}
\newcommand{\CO}{{\mathcal O}}
\nc{\ca}{{\mathcal A}}
\nc{\cag}{{{\mathcal A}^\Gamma}}
\nc{\cg}{{\mathcal G}}
\nc{\chh}{{\mathcal H}}
\nc{\ck}{{\mathcal B}}
\nc{\cd}{{\mathcal D}}
\nc{\cl}{{\mathcal L}}
\nc{\cm}{{\mathcal M}}
\nc{\cn}{{\mathcal N}}
\nc{\cs}{{\mathcal S}}
\nc{\cz}{{\mathcal Z}}
\nc{\sind}{\sigma{\rm -ind}}

\newcommand\clFN{{\mathcal F_\tau(\mathcal N)}}       
\newcommand\clKN{{\mathcal K_\tau(\mathcal N)}}       
\newcommand\clQN{{\mathcal Q_\tau(\mathcal N)}}       %
\newcommand\tF{\tilde F}
\newcommand\clC{\mathcal C}
\newcommand\clE{\mathcal E}
\newcommand\clF{\mathcal F}
\newcommand\clK{\mathcal K}
\newcommand\clL{\mathcal L}
\newcommand\clN{\mathcal N}
\newcommand\clM{\mathcal M}
\newcommand\clS{\mathcal S}
\newcommand\clB{\mathcal B}
\newcommand\Del{\Delta}
\newcommand\g{\gamma}
\newcommand\eps{\varepsilon}
\newcommand\vf{\varphi}
\newcommand\E{\mathcal E}

\newcommand{\CDA}{\mathcal{C_D(A)}} 
\newcommand{\dslash}{{\pa\mkern-10mu/\,}}

\newcommand{\sepword}[1]{\quad\mbox{#1}\quad} 
\newcommand{\comment}[1]{\textsf{#1}}   

\nc{\nt}{\newtheorem}
\nc{\bra}{\langle}
\nc{\ket}{\rangle}
\nc{\cal}{\mathcal}
\nc{\frk}{\mathfrak}

\parindent=0.0in

 \title{Dixmier traces and extrapolation description of noncommutative  Lorentz spaces}

\author{Victor GAYRAL}
\address{Laboratoire de Math\'ematiques\\
Universit\'e de Reims Champagne-Ardenne\\
Moulin de la Housse-BP 1039, 51687 Reims, FRANCE\\
e-mail: victor.gayral@univ-reims.fr}
\author{Fedor  SUKOCHEV}
\address{School of Mathematics and Statistics,
University of New South Wales\\
Kensington NSW 2052, AUSTRALIA\\
e-mail: f.sukochev@unsw.edu.au}

\begin{abstract}
We study the connections  between Dixmier traces, $\zeta$-functions and traces of heat semigroups
beyond the dual of the  Macaev ideal and 
in the general context of  semifinite von Neumann
algebras. We show that the correct framework for this investigation is that of operator Lorentz spaces 
possessing  an extrapolation description. We demonstrate the applicability of our results to
H\"ormander-Weyl pseudo-differential calculus on $\R^n$. In this context, we prove that the Dixmier 
trace of a pseudo-differential operator coincides with the `Dixmier integral' of its symbol.
\end{abstract}

\maketitle

{\Small{\bf Keywords:}
Dixmier traces, singular traces, Lorentz (Marcinkiewicz) spaces, $\zeta$-functions, heat 
kernels, extrapolation, H\"ormander-Weyl pseudo-differential operators }

\qquad 

\qquad

\section*{Introduction}

This article is devoted  to the study of singular traces
beyond $\mathcal L^{1,\infty}(\cH)$, the  dual of the Macaev ideal. 
 The theory of singular traces has been intensively developed recently,   especially 
after Alain Connes' application of that theory to noncommutative geometry
and physics \cite{Co0,Co4,Co}. 
In particular, singular traces now have a broad application to several fields of 
mathematics and theoretical physics, including  noncommutative symmetric spaces, geometry of 
Banach spaces,  pseudo-differential operators, index theory and geometric analysis, fractal 
geometry, noncommutative integration theory, renormalization and quantum field theory...
To illustrate that, and without any pretension to exhaustivity,   see for instance
\cite{BF,CGRS1,CPS2,CEL,CiS,CGS,ER,G1,GIV,GI,GI2,GIL,I,L,LMP,NS,N,NR,P,Po,VGB}.

So far, known results and applications concentrate around singular traces arising from  the 
logarithmic  divergence of a normal  trace. 
However,   we strongly believe that 
 singular traces associated with more general
  divergences  should have natural applications in all the above mentioned areas.
 In the context of quantum groups, we give such an example after Corollary \ref{OUOU}.
We also demonstrate in section 5 that general singular traces have natural applications to the theory
  of  H\"ormander-Weyl pseudo-differential calculus on $\R^n$.

  More specifically, the main motivation of the present article is  to characterize the class of noncommutative Lorentz  
spaces of $\tau$-compact operators within a semifinite 
 von Neumann algebra, for which one can generalize the `classical' relationships between:
  Dixmier traces,  residues of $\zeta$-functions   and  
  asymptotics of  traces of  heat semigroups. 
 Our main discovery is the existence of a strong connection between 
  singular traces and extrapolation theory. At the present stage, extrapolation
  only appears as a technical device   for the study of singular traces. We leave a more profound 
  investigation of their interconnections for  a future work.  
 
 \quad
 
Let us first review   the `classical' results we wish to extend to more general Lorentz spaces.
 Let
  $\clM^{1,\infty}(\mathcal H)$  be the  dual of the  Macaev 
ideal   of compact operators on a separable Hilbert space $\mathcal H$:
 $$
\clM^{1,\infty}(\mathcal H):=\Big\{T\in \mathcal K(\mathcal H)\;:\;\|T\|_{1,\infty}:=\sup_{n\in\N}\frac1{\log(1+n)}\sum_{k=1}^n\mu(T,k)<\infty\Big\}.
$$
Here  $\mathcal K(\mathcal H)$ denotes the ideal of compact operators 
and $\{\mu(T,n) \}_{n\in\N}$ is the sequence of singular values of a compact operator $T$.
We will always embed $\ell^\infty(\N)$ onto $L^\infty(\R_+^*)$ via the map
$$
\{x_n\}_{n\in\N}\mapsto \sum_{n\in\N} x_n\chi_{[n,n+1)},
$$
and, according to this embedding, a state on $L^\infty(\R_+^*)$ defines a state on 
$\ell^\infty(\N)$.
In his seminal paper \cite{Dix},  Jacques Dixmier constructed singular traces
(i.e$.$ traces  vanishing  on finite rank operators), which for
$\clM^{1,\infty}(\mathcal H)$
  are given by the linear extension of the unitarily invariant  weights
\begin{align}
\label{TRW}
{\rm Tr}_\omega:\mathcal B(\mathcal H)^+\to [0,\infty],\quad
T\mapsto \omega\Big(\Big\{\frac1{\log(1+n)}\sum_{k=1}^n\mu(T,k)\Big\}_{n\in\N}\Big),
\end{align}
where $\omega$ is an arbitrary dilation invariant (see \eqref{invariance} for the definition)
 state on $L^\infty(\R_+^*)$.
 We stress that  the 
computation of a Dixmier trace can be an
highly non-trivial task in  concrete examples. This is not really stipulated by the presence of an
extended limit $\omega$
in the formula \eqref{TRW}. To  our knowledge and at least within the framework of spectral triples in  noncommutative 
geometry,
 all the examples of
computations of Dixmier traces involve   `measurable operators', that is operators in 
$\clM^{1,\infty}(\mathcal H)$ for which all the Dixmier traces take the same value. For
instance, this is the case for classical pseudo-differential operators on a closed manifold\footnote{
A very short and simple proof of that fact is displayed in \cite[Theorem 3]{KLPS}.
Note also that a concrete example a	 of non-measurable operator in the context of
 pseudo-differential operators on $\mathbb R^n$ with non-homogeneous symbol, is also given there.}.
The real problem for the computation of a Dixmier trace 
is that, typically,  one   has no access 
to the sequence of singular values of an operator.
 This explains why in concrete applications of Dixmier traces,  the defining formula \eqref{TRW}
 is of little utility. 
A remedy to this problem  exists and consists in expressing a Dixmier traces ${\rm Tr}_\omega$
 in term of 
the ordinary trace $\rm Tr$, which is a computable object.
The following theorem (for which we refer to \cite{Co4}, see also \cite{BF,CPS2,CRSS,SUZ,SZ})
is an example to such a remedy and a main sample of `classical' result which we set to generalize.
 \begin{thm}
 \label{zero}
 Let $T\in\clM^{1,\infty}(\mathcal H)^+$, $B\in\clB(\mathcal H)$
  and $\omega$ be an exponentiation  invariant (see \eqref{invariance})  state on $L^\infty(\R_+^*)$.
  Then we have  
  \begin{equation}
  \label{triple}
  {\rm Tr}_\omega(BT)=\omega\Big(\Big[r\mapsto \frac{\zeta_B(T,1+\log(r)^{-1})}{\log(1+r)}
  \Big]\Big)=\omega\Big(\Big[\lambda\mapsto  \frac1{\log(1+\lambda)}\int_0^\lambda 
  \xi_B(T,t^{-1})\,\frac{dt}{t^2}\Big]\Big).
  \end{equation}
\end{thm}
 Here, $\zeta_B(T,.)$ and $ \xi_B(T,.)$ are the generalized $\zeta$-function and heat-trace function,
respectively defined  by
 $$
 \zeta_B(T,z):={\rm Tr}(BT^z),\quad z\in \C, \;\Re(z)>1\qquad
 \mbox{and}\qquad  \xi_B(T,t):={\rm Tr}\big(Be^{-tT^{-1}}\big),\quad t>0.
 $$

The  first version of this result (for $T$ measurable, $B=1$ and with different assumptions on 
$\omega$)
 is due to Alain Connes and appeared already in his book \cite{Co4}. 
It has then been followed by an important list of   improvements and generalizations 
of many kinds  \cite{BF,CPS2,CRSS,SZ}. In the form we phrase it here, this result
 first appeared in \cite{SUZ}. We do not want to do a complete review of the different steps leading
 to this `final' answer. The reader interested in the history of this result may consult the 
  book \cite{LSZ}. 
  
  \quad
  
We should also mention two  other approaches to obtain alternative
 expressions for a Dixmier trace. The first one is  a spectral formula, 
 which is somewhat analogous
 to the Lidskii's formula for the ordinary trace \cite{SSZ}:
 $$
  {\rm Tr}_\omega(T)=\omega\Big(\Big\{\frac1{\log(1+n)}\sum_{\lambda\in\sigma(T),\;|\lambda|\geq\log(n)/n}\lambda\Big\}_{n\in\N}\Big),
 $$
 where $T$ is arbitrary in $\clM^{1,\infty}(\mathcal H)$.
 One of the main features of this formula  is that it extends to more general operator ideals than 
 $\clM^{1,\infty}(\mathcal H)$, in particular to the Lorentz ideal $\clM_\psi(\mathcal H)$ 
 (with $\psi\in\Omega$ satisfying condition \eqref{SSZ-cond} -- see below).    The second alternative approach  relies on the 
 results of \cite{KLPS}, and allows one to compute Dixmier traces from expectation values 
 \cite[Corollary 8.2.4]{LSZ}:
 $$
  {\rm Tr}_\omega(T)=\omega\Big(\Big\{\frac1{\log(1+n)}\sum_{k=1}^n\langle Te_k,e_k\rangle
  \Big\}_{n\in\N}\Big). 
 $$
Here $T$ is required to belong to $\mathcal L^{1,\infty}(\mathcal H)$ (the non-closed subspace of 
$\clM^{1,\infty}(\mathcal H)$ whose elements have their  sequence of singular values
dominated by the harmonic sequence), and is such that there exists 
$0\leq V\in\mathcal L^{1,\infty}(\mathcal H)$ with $\sup_{t>0}t^{1/2}\|T(1+tV)^{-1}\|_2<\infty$
and $\{e_n\}_{n\in\N}$ is an eigenbasis for $V$, ordered in such a way  that the corresponding sequence
of eigenvalues is decreasing. This formula has a great impact on the theory of pseudo-differential
operators \cite{KLPS}, and especially on Connes' Trace Theorem \cite{Co0},
 which connects  Dixmier traces 
and
Wodzicki residue \cite{W1,W2}.

 \quad

We now explain to what extent we  generalize Theorem \ref{zero}. Fix $(\clN,\tau)$  a
 semifinite von Neumann algebra and a normal faithful semifinite trace. In that context,
 there is an analogous notion of singular numbers, the generalized $s$-numbers
  (see for instance \cite{FK}). The latter object is a
 (right-continuous and decreasing) function on the positive half-line, $\mu(T)$, attached
 to every  operator $T$ in $\clN$
(see  \eqref{SN} for the precise definition). The generalized $s$-numbers
extend both singular numbers of compact operators
  and decreasing rearrangement of measurable functions
 on $\sigma$-finite measure spaces. Let then $\psi:[0,\infty)\to [0,\infty)$ be a concave and increasing 
 function, vanishing at zero and diverging at infinity. We denote by
  $\Omega$ the set of such functions.
  Then, one can generalize the  ideal $\clM^{1,\infty}(\mathcal H)$
 by setting:
 $$
\clM_\psi(\clN,\tau):=\Big\{T\in \clN\;:\;\|T\|_{\clM_\psi}:=\sup_{t>0}\frac1{\psi(t)}\int_0^t\mu(T,s)
\,ds<\infty\Big\},\qquad\psi\in\Omega.
$$
 Of course, when $\clN=\mathcal B(\mathcal H)$ and
 $\tau={\rm Tr}$, $\mu(T)$ is a step function  and if moreover
  $\psi(t)=\log(1+t)$, we have $\clM_\psi(\clN,\tau)=\clM^{1,\infty}(\mathcal H)$ isometrically.

A  comment on terminology is in order. In many sources, including most of our own papers
on this subject, 
these ideals are called (operator) Marcinkiewicz spaces. However,
 in a commutative setting, these spaces were 
first introduced by G. Lorentz in 1950 \cite{Lorentz}. In this paper, 
 we refer to these spaces as  (operator) Lorentz spaces, 
which is consistent with the terminology employed in the book of Bennett
and Sharpley \cite{BS}. Again, we refer to the 
book \cite{LSZ}, where historical details of the development of the theory of Lorentz
spaces are given.\\
With these definitions at hand, 
the extension of Dixmier's construction of singular traces\footnote{For type $I_\infty$
factors and  $\psi\in\Omega$
satisfying $\lim_{t\to\infty}\psi(2t)/\psi(t)=1$, the following construction already appearded  in \cite{Dix}.}, 
to more general Lorentz spaces is rather direct. Namely, given $\psi\in\Omega$ and 
$\omega$  an arbitrary dilation invariant state on $L^\infty(\R_+^*)$, one  considers the
 unitarily invariant  functional:
\begin{align}
\label{DT}
\tau_{\psi,\omega}:\clM_\psi(\clN,\tau)^+\to [0,\infty),\quad
T\mapsto \omega\Big(\Big[t\mapsto\frac1{\psi(t)}\int_0^t\mu(T,s)\,ds\Big]\Big).
\end{align}
It is then shown in \cite{KSS} that the latter is positively additive (i.e$.$ it defines a weight on $\clN$)
if and only if $\omega$ is compatible with $\psi$,  in the sense that
$$
\omega\Big(\Big[t\mapsto \frac{\psi(at)}{\psi(t)}\Big]\Big)=1,\qquad \forall a>1.
$$
In this case, the linear extension of the
functional \eqref{DT} defines a trace on $\clM_\psi(\clN,\tau)$,
by which we mean a unitarily invariant linear functional. 
Moreover, these traces are  singular, that is    they 
vanish on projections of 
 finite trace $\tau$. 
 The main difference with the classical case of $\mathcal B(\mathcal H)$,
 is that in the general semifinite setting,
these singular traces may be supported either at zero or at infinity (or at both). By this 
we mean that  they may either
come from the lack of integrability of the generalized $s$-numbers function $\mu(T)$ at infinity, or
of the unboundness of $\mu(T)$ at zero (or both). (Equivalently,  a singular trace 
supported at infinity vanishes on $\mathcal L^1(\clN,\tau)$ and a singular trace supported at zero
vanishes on $\clN$.)
Here however, we are mainly concerned with
 singular traces supported at infinity.
 
Again, having in general no practical  access to the generalized $s$-numbers  of an
operator in a semifinite von Neumann algebra, 
the computation of these singular traces
remains an elusive problem.
This is mostly why we aim to generalize Theorem \ref{zero}. 
The hope is that 
Theorem \ref{zero} holds  for more general Lorentz spaces, 
with $\tau$ instead of ${\rm Tr}$ (including in the definition of the $\zeta$-function $\zeta_B$
 and of the heat-kernel
function $\xi_B$)
and with $\psi$ instead of $\log$ in the denominators of the double equality \eqref{triple}.
(This is known to be true when $\psi$ behaves like a logarithm at infinity and for any semifinite 
von Neumann algebra and related ideas can be found in \cite{CRSS,S1,S2}.) 
More   precisely, the  question we are addressing  in this paper is the following:

\quad

{\it  What is the class of functions $\psi\in\Omega$ for which
we have
  \begin{align}
  \label{un}
  {\tau}_{\psi,\omega}(BT)=C_\zeta(\psi)\,\omega\Big(\Big[r\mapsto \frac{\zeta_B(T,1+\log(r)^{-1})}{\psi(r)}
  \Big]\Big)=C_\xi(\psi)\,\omega\Big(\Big[\lambda\mapsto  \frac1{\psi(\lambda)}\int_0^\lambda 
  \xi_B(T,t^{-1})\,\frac{dt}{t^2}\Big]\Big),
  \end{align}
  for all  semifinite von Neumann algebras $\clN$ with a fixed
   normal semifinite faithful trace $\tau$, 
  all positive elements $T$ in the Lorentz space $\clM_\psi(\clN,\tau)$,
 all  $B$ in $\clN$, all   exponentiation  invariant  states $\omega$ on $L^\infty(\R_+^*)$
 and for certain positive constants $C_\zeta(\psi)$ and $C_\xi(\psi)$ which depend  on $\psi$
 only?}
 
 \quad
 
 Here the functions $\zeta_B$ and $\xi_B$ are defined exactly in the same fashion as their `classical'
counterparts in \eqref{triple} via replacing $\rm Tr$ with $\tau$.

We do not fully answer this question as  we have only found  sufficient 
conditions on $\psi\in\Omega$ for \eqref{un} to hold true. However, these conditions are
wide  enough to cover new applications. For instance, they apply to $\psi(t)=\log(1+t^{1/\beta})
^\beta$, $\beta>0$, and to $\psi(t)=\log(1+\dots\log(1+t)\dots)$.
The first condition that we have found, and this is the main novelty of our approach, 
relies on  an extrapolation description  of the Lorentz space 
$\clM_\psi(\clN,\tau)$. The extrapolated space we need to consider  
 first appeared in  \cite{AL1,AL2} and is of the following nature.
  Let $\|.\|_p$, $1\leq p<\infty$, be 
the  $L^p$-norm associated to $\tau$. Then, we let
$\mathfrak L_\psi(\clN,\tau)$  be the Banach  space of  measurable operators (relative to 
 $\clN$ and $\tau$)  such that
$$
\|T\|_{\mathfrak L_\psi}:=
\sup_{p>1}\frac{\|T\|_p}{\psi(e^{(p-1)^{-1}})}<\infty.
$$
It can be proven that  $ \mathfrak L_\psi(\clN,\tau)$ always embeds in $\clM_\psi(\clN,\tau)$ 
and our first requirement is
\begin{align}
\label{C1}
 \mathfrak L_\psi(\clN,\tau)=\clM_\psi(\clN,\tau).
 \end{align}
The second condition we need to impose, can be interpreted as  the existence of 
an exponentiation index
for $\psi$:
\begin{align}
\label{C2}
\quad \mbox{for all}\quad\alpha>1\,,\quad
\mbox{the limit}\quad \lim_{t\to\infty}\frac{\psi(t^\alpha)}{\psi(t)}\quad
\mbox{exists}.
\end{align}
We can now state our answer to the question above, which is the combination of 
Theorem \ref{zeta} and Theorem \ref{psi2}, the main results of this paper:
\begin{thm}
\label{deux}
Let $(\clN,\tau)$ be a semifinite von Neumann algebra with a normal semifinite faithful trace
and let  $\psi$ be an element of $\Omega$ satisfying conditions \eqref{C1} and \eqref{C2}.
Then the double equality \eqref{un} holds for all $T\in  \clM_\psi(\clN,\tau)^+$,  all  $B\in\clN$  
and all  exponentiation  invariant (see \eqref{invariance}) states $\omega$  on $L^\infty(\R_+^*)$ with
 $$
 C_\zeta(\psi)=\Gamma\Big(1+\log\Big(\lim_{t\to\infty}\frac{\psi(t^e)}{\psi(t)}\Big)\Big)^{-1}
 \qquad \mbox{\rm and}\qquad C_\xi(\psi)=1.
 $$
\end{thm}

Some comments are in order.
Observe first that  the condition \eqref{C1}  is  necessary 
for the middle term in \eqref{un} to be well defined for all $T\in\clM_\psi(\clN,\tau)^+$,
 and that condition  \eqref{C2} warrants the 
existence of the constant $C_\zeta(\psi)$. Hence, we believe that conditions 
\eqref{C1} and \eqref{C2} are   necessary as well, for the first equality in \eqref{un}
to hold true. However, we  do not believe that conditions \eqref{C1} and \eqref{C2} are   
necessary for the equality between the term on the left and the term on the right in \eqref{un}.
The reason for this  is  that the term on the right  in \eqref{un}  is well defined 
under a strictly  weaker condition than \eqref{C2} (Corollary \ref{equal}).   

\quad

We  stress that Theorem 
\ref{deux} is new even for type $I_\infty$ factors. To illustrate it, we give 
in Section \ref{psido} an application to Weyl pseudo-differential operators on $\R^n$:
$$
{\rm OP}_W: \clS'(\mathbb R^{2n})\to\clL\big( \clS(\mathbb R^{n}), \clS'(\mathbb R^{n})\big),
$$
with symbols in the general H\"ormander classes $S(m,g)$. 
To formulate the results obtained in this section, we set $\clM_\psi(\R^{2n})$ for the commutative 
Lorentz
space  associated
with $L^\infty(\R^{2n})$ and Lebesgue integral  and
 $\clM_\psi\big(L^2(\R^{n})\big)$ for the type $I_\infty$ Lorentz space associated 
with $\mathcal B\big(L^2(\R^{n})\big)$ and operator trace. We also set $\int_{\psi,\omega}$
for a Dixmier trace on $\clM_\psi(\R^{2n})$, that we may call a `Dixmier integral'.
In this context, we are able to deduce the following consequence
of the relation between the Dixmier trace and the trace of the heat semigroup, as stated in Theorem
\ref{deux}.

\begin{thm}
Let $\psi\in\Omega$ satisfying   conditions \eqref{C1} and \eqref{C2} and let also
$(g,m)$  be an H\"ormander pair (see Definition \ref{HP}) such that
$m\in  \clM_\psi(\mathbb R^{2n})$.
Then, 
$$
{\rm OP}_W:S(m,g)\to \clM_\psi\big(L^2(\mathbb R^n)\big),
$$
 continuously. Moreover, for any symbol $f\in S(m,g)$ 
and any exponentiation invariant state $\omega$
on $L^\infty(\R_+^*)$, we have
$$
{\rm Tr}_{\psi,\omega}\big({\rm OP}_W(f)\big)=\int_{\psi,\omega} f.
$$ 
\end{thm}
Here, $\int_{\psi,\omega}$ is the Dixmier trace on the commutative 
von Neumann algebra $L^\infty(\R^{2n})$ constructed out of the semifinite
trace given by the Lebesgue integral.
This result is proven in the text in Theorems \ref{Nique} and  \ref{ils}.
Observe that it uses the full strength of the theory of singular traces for
semifinite von Neumann algebras. It also
 nicely complements the classical relation between the operator trace of a
 Weyl pseudo-differential operator and the integral of its symbol:
$$
{\rm Tr}\big({\rm OP}_W(f)\big)=\int f,
$$ 
whenever ${\rm OP}_W(f)$ is trace class and $f$  is integrable\footnote{This occurs for instance when $f$
is a Schwartz function.}.
 In particular,
 the  relation between Dixmier traces and
$\zeta$-functions of Theorem \ref{deux}, applied to $L^\infty(\R^{2n})$,
gives  a very elementary way to compute
the Dixmier trace of an H\"ormander-Weyl pseudo-differential operator:
\begin{theorem}
\label{quatre}
Let $\psi\in\Omega$ satisfying conditions \eqref{C1} and \eqref{C2} and let also
$(g,m)$  be an H\"ormander pair such that
$m\in  \clM_\psi(\mathbb R^{2n})$.
Then, for any  symbol $ f\in S(m,g)$ and any exponentiation invariant state $\omega$
on $L^\infty(\R_+^*)$, we have:
$$
{\rm Tr}_{\psi,\omega}\big({\rm OP}_W(f)\big)= C_\zeta(\psi)\,
\omega\Big(\Big[r\mapsto\frac1{\psi(r)}\int_{\R^{2n}}f(x,\xi) \;|f(x,\xi)|^{1/\log(r)}\;d^nxd^n\xi\Big]\Big),
$$
where $ C_\zeta(\psi)$ is the constant (associated to $\psi$) given in Theorem  \ref{deux}.
\end{theorem}
 
We should also mention that  Theorem \ref{quatre} gives  very simple alternative proofs
of \cite[Proposition 4.17]{G1} and \cite[Theorem 4.7]{N}.

\qquad

{\bf Acknowledgments.} It is a pleasure to thank Sergei Astashkin, Bruno Iochum, Steven Lord,
Evgen\^ii Semenov, 
Alexandr Usachev
 and Dmitriy Zanin  for various
discussions and suggestions during the development  of this project.

\section{Preliminary results}
\label{1}

In this section, we first recall the basis of the theory of  operator Lorentz spaces in the semifinite 
setting. We introduce then a subclass of Lorentz spaces on which we  focus.
This subclass is characterized by a single condition, which is
displayed in  \eqref{our-cond}.
In Proposition \ref{hier}, we compare this condition with those which  are most
 commonly used  in the literature (conditions \eqref{SSZ-cond} and \eqref {ex-sing-trace-easy}
 in the text). 
 Next, we recall the fundamentals of singular traces on Lorentz spaces.
The main role is played by the notion of $\psi$-compatible state and
 we construct  singular traces in the context of  exponentiation 
 invariant $\psi$-compatible  states  (Proposition \ref{lin}). 
We conclude
  this section with a spectral formula for singular traces, again in the context of exponentiation invariant (Proposition \ref{calc}).
\subsection{Notations }

\noindent
Set  $\N:=\{1,2,3,\dots\}$, $\N_0:=\{0,1,2,\dots\}$, 
$\R_+=[0,\infty)$ and $\R_+^*=(0,\infty)$.
Let $\clN$ be a semifinite von Neumann algebra endowed with $\tau$, a faithful semifinite normal 
trace. We denote by $\clL_0(\clN,\tau)\equiv \clL_0$ the $*$-algebra of $\tau$-measurable operators
affiliated with $\clN$. Recall \cite{FK} that a densely defined closed operator $T$ affiliated 
with $\clN$ is   $\tau$-measurable if 
$$
\lim_{\lambda\to\infty}\tau\big(1-E_{|T|}(\lambda)\big)<\infty,
$$
 where 
$E_{|T|}$ denotes the spectral family of $|T|$. For
$p\in[1,\infty)$, we set $\clL^p(\clN,\tau)\equiv \clL^p$ for the  noncommutative 
$L^p$-space:
$$
\clL^p:=\big\{T\in\clL_0\; :\;\|T\|_p:=\tau(|T|^p)^{1/p}<\infty\big\}.
$$ 
We also consider $\clF(\clN,\tau)\equiv\clF$,  the subset  of $\clN$ consisting of finite linear 
combinations of projections of finite trace  $\tau$ and we denote by $\clK(\clN,\tau)\equiv\clK$ the
norm closure of $\clF$ in $\clN$. The latter is termed the algebra of $\tau$-compact operators. 
For $T\in\clL_0$, a $\tau$-mesurable operator, we let $\mu(T):=[t\mapsto\mu(T,t)]$ 
be the generalized $ s$-numbers  function. Recall that the latter is  defined by
\begin{equation}
\label{SN}
\mu(T,t):=\inf\big\{\|PT\|\;:\;\tau(1-P)\leq t\big\}\,,\quad t>0,
\end{equation}
where the infimum runs over the  lattice of projections in $\clN$. 
Observe that  $\mu(T)$  depends on the choice of the
trace $\tau$. To lighten the notations, we do not make this dependance explicit.
We also define the distribution function, $n(T)=[t\mapsto n(T,t)]$, of a $\tau$-measurable
operator $T$ by:
$$
n(T,t):=\tau\big(\chi_{(t,\infty)}(T)\big)\,,\quad t>0,
$$
where $\chi_{(t,\infty)}$ is the indicator function of the interval $(t,\infty)$.
It can be shown that
$$
\mu(T,s)=\inf\big\{t>0\;:\;n(T,t)\leq s\big\}\qquad \mbox{and}\qquad
n(T,t)\leq s\;\Leftrightarrow \;\mu(T,s)\leq t.
$$
If $T$ and $ S$ belong to $\clL_0$, we say
that $T$ majorizes   $S$ in the sense of Hardy-Littlewood-P\'olya, 
denoted $S\prec\prec T$, if and only if
$$
\int_0^t\mu(S,s)\,ds\leq \int_0^t\mu(T,s)\,ds\,,\quad \forall t>0.
$$
Specializing this to the commutative von Neumann algebra of essentially bounded
functions on a $\sigma$-finite measure space $(X,\mu)$ with normal semifinite faithful  trace
given by the Lebesgue integral, for measurable functions $f,g:X\to\mathbb C$ which are
finite outside a measurable set of finite measure, we accordingly use the notation 
$f\prec\prec g$ and recall that in this case, $\mu(f)$ coincides we the nondecreasing 
rearrangement of $f$.

\noindent 
Let $\clE(\clN,\tau)\equiv\clE$ be a Banach subspace of $\clL^1+\clN$, with norm $\|.\|_\clE$.
Then $\clE$ is called a {symmetric operator space} if for all $A\in\clE^+$ and 
$B\in(\clL^1+\clN)^+$ with $B\leq A$, we have $B\in\clE$ with $\|B\|_\clE\leq\|A\|_\clE$
and if for all $ A\in\clE$ and 
$ B\in\clL^1+\clN$ with $\mu(B)=\mu( A)$, we have $B\in\clE$ with $\|B\|_\clE=\|A\|_\clE$.
Observe that $\clE$  is a bimodule over $\clN$ (or an ideal in $\clN$ if $\clE\subset\clN$)
since the condition above implies that for all $A\in\clE$ and all $B,C\in\clN$, we have $BAC\in\clE$
with $\|BAC\|_{\clE}\leq \|B\|\,\|A\|_\clE\,\|C\|$.
If moreover for all $ A\in\clE$ and 
$ B\in\clL^1+\clN$ with $B\prec\prec  A$, we have $B\in\clE$ with $\|B\|_\clE\leq\|A\|_\clE$,
then $\clE$ is called a fully symmetric operator space. Symmetric (resp. fully symmetric)
operator spaces are the noncommutative counterparts  of symmetric (resp. fully
symmetric) Banach spaces of  (classes of) measurable functions on $\sigma$-finite
measure spaces. Note that given $(E,\|.\|_E)$, a Banach space of 
Lebesgue measurable functions on the interval 
$[0,\infty)$, we can construct an operator space, by setting
$$
E(\clN,\tau):=\big\{T\in\clL_0\,:\,\mu(T)\in E\big\}.
$$
Normed with $T\mapsto\| \mu(T)\|_E$,
the latter is  symmetric (resp. fully symmetric) when $E$ is 
symmetric  (resp. fully symmetric).
 If moreover $\clN$ is either non-atomic or a type $I_\infty$ factor, then every  symmetric 
 (resp. fully symmetric) operator space
arises this way \cite{CS, DDP, KS, SC}. 
When a  symmetric operator space $\clE$ is of the form $E(\clN,\tau)$, for some Banach space
$E$ of functions on $[0,\infty)$, we define its 
fundamental function $\varphi_\clE$  as the corresponding fundamental function of $E$:
 $$
 \varphi_{\clE}(t):=\|\chi_{(0,t)}\|_E.
 $$
In this situation,  we also define the lower and upper Boyd indices of $\mathcal E$, 
  as the ordinary Boyd indices of 
$E$. 
Hence, with $D_a$, $a>0$,  the dilation operator 
acting on functions on $[0,\infty)$ via $D_a f(t):=f(at)$,
we set\footnote{We warn the reader that in many sources, the  dilation operator $D_a$ is 
defined  with $a^{-1}$ instead of $a$. }
$$
\mathfrak p_\clE:=\lim_{a\downarrow 0}\frac{\log a^{-1}}{\log\|D_a\|_{E\to E}}\,,\qquad \mathfrak q_
\clE:=\lim_{a\to\infty}\frac{\log a^{-1}}{\log\|D_a\|_{E\to E}}.
$$
Recall that we always have $1\leq \mathfrak p_\clE\leq \mathfrak q_\clE\leq\infty$.

\noindent 
Whenever $f$ and $g$ are real valued function on a set $X$, we write $f\asymp g$ when
there exist $0<c\leq C$ such that for all $x\in X$, we have $cf(x)\leq g(x)\leq Cf(x)$.

\subsection{Lorentz spaces}
An important class of fully symmetric operator spaces is given by the 
 Lorentz spaces. To introduce this class of operator spaces, we consider the following class of
 concave functions:
 \begin{definition}
Let $\Omega$ be the set of   all increasing concave  functions
 $\psi:[0,\infty)\to [0,\infty)$ 
such that
\begin{align*}
\psi(0)=0\qquad\mbox{and}\qquad \lim_{t\to\infty}\psi(t)=\infty,
\end{align*}
and let $\Omega_b$ be the subset of $\Omega$ whose elements satisfy further: 
$$
\psi(t)=O(t),\quad t\to 0.
$$
\end{definition}
  Associated with an element $\psi\in\Omega$, one defines the operator Lorentz  space 
 $ \clM_\psi:=\clM_\psi(\R_+)(\clN,\tau)$ to be  the collection of all $\tau$-measurable 
 operators, such that
 $$
\|T\|_{\clM_\psi}:=\sup_{t>0} \frac1{\psi(t)}\int_0^t\mu(T,s)\,ds<\infty.
 $$  
 We  denote by $\clM_\psi^0$, the closure of $\clF$  in $\clM_\psi$ and by $\clM_\psi^+$, the
positive cone of $\clM_\psi$. We also use the notation $\clM_\psi(\mathbb R_+)$ for the commutative
Lorentz space of (classes of) functions on the positive half-line, i.e$.$ for $\clN=L^\infty(\R_+)$
and $\tau=\int_0^\infty\cdot \,dx$.
 Note  that $t/\psi(t)$ is the fundamental function of $\clM_\psi$ and that 
  $\clM_\psi$   is continuously embedded in $\clN$ if and only if  $\psi\in\Omega_b$
 or equivalently if $\psi'\in L^\infty(\R_+)$.
  Indeed, as $\psi'\in\clM_\psi(\mathbb R_+)$ and $\|\psi'\|_\infty=\lim_{t\downarrow0}\psi'(t)$, we deduce that 
  $$
  \clM_\psi\subset\clN\quad\Rightarrow\quad\psi'\in L^\infty(\R_+)
  \quad\Rightarrow\quad\psi(t)\leq \|\psi'\|_\infty\,t.
  $$   
  Conversely, if $\psi(t)\leq C t$ for $ t\in( 0,1)$, we obtain 
  $$
  \|T\|_{\clM_\psi}\geq\frac{\int_0^t\mu(T,s)\,ds}{\psi(t)}\geq\frac {\mu(T,t)}{C},\quad \forall t\in(0,1),
  \quad \forall T\in\clM_\psi,
  $$
  which entails that $T\in\clN$ as 
  $$
  \|T\|=\lim_{t\downarrow 0}\mu(T,t).
  $$
 This justifies our notation for $\Omega_b$. Observe that when $\psi\in\Omega_b$, then
 the Lorentz space $\clM_\psi$ consists of bounded operators and therefor the latter is an ideal in
 $\clN$.
  The distinction between $\Omega$ and $\Omega_b$,
 plays a minor role in this paper. In fact, this  distinction will only appear in Section \ref{2} (which
 is also of independent interest) while in Sections \ref{3} and \ref{4}, even for $\psi\in\Omega
 \setminus \Omega_b$, we are forced to consider the elements of $\clM_\psi\cap\clN$ only.
 
\quad

\noindent
In this article, we are mostly interested in   Lorentz  spaces
associated with the subset of  $\Omega$ whose elements are 
characterized by the  additional condition:
\begin{align}
\label{our-cond}
\forall\alpha>1\,,\quad
\mbox{the limit}\quad A_\psi(\alpha):=\lim_{t\to\infty}\frac{\psi(t^\alpha)}{\psi(t)}\quad
\mbox{exists}.
\end{align}
We will see in a moment that  condition \eqref{our-cond} is strictly   stronger than  the condition
 used in \cite{SSZ}:
\begin{equation}
\label{SSZ-cond}
\lim_{t\to\infty}\frac{\psi\big(t\psi(t)\big)}{\psi(t)}=1,
\end{equation}
and that the latter is strictly stronger than the  standard  condition (used in particular in 
\cite{Dix} and \cite{CRSS}):
\begin{equation}
\label{ex-sing-trace-easy}
\forall \,a>1\,,\qquad\lim_{t\to\infty}\frac{\psi(at)}{\psi(t)}=1.
\end{equation}

\noindent We first  observe that  Theorem 6 of  \cite{KSS}, proved there for extended limits,
 can be adapted to the
 ordinary limit:
\begin{lemma}
\label{eq-conds-1}
Let $\psi\in\Omega$. Then
  condition \eqref{ex-sing-trace-easy} is equivalent to any of the
following two conditions:
$$
\exists \,a>1\;:\;\lim_{t\to\infty}\frac{\psi(at)}{\psi(t)}=1\,,\qquad \qquad
\lim_{t\to\infty}\frac{t\psi'(t)}{\psi(t)}=0.
$$
\end{lemma}
\noindent
\begin{proof}
Call (\ref{ex-sing-trace-easy}$'$) the condition 
$\exists \,a>1$ : $\lim_{t\to\infty}\frac{\psi(at)}{\psi(t)}=1$
and (\ref{ex-sing-trace-easy}$''$) the condition $\lim_{t\to\infty}\frac{t\psi'(t)}{\psi(t)}=0$.
Observe first that the concavity and the monotony of $\psi$,
imply  that  for any $a>1$ there exists $C>0$  such that:
 $$
 0\leq\frac{at\psi'(at)}{\psi(at)}\leq C\frac{\psi(at)-\psi(t)}{\psi(t)},
 $$
 and it also imply that for any $a>1$:
  $$
 1\leq\frac{\psi(at)}{\psi(t)}\leq 1+(a-1)\frac{t\psi'(t)}{\psi(t)},
 $$ 
 Hence we deduce that  $(\ref{ex-sing-trace-easy})\Rightarrow (\ref{ex-sing-trace-easy}')
 \Rightarrow (\ref{ex-sing-trace-easy}'')\Rightarrow (\ref{ex-sing-trace-easy})$.
   \hfill
  \end{proof}

\noindent
We next observe a simple but important consequence of  condition \eqref{our-cond}:
\begin{lemma}
\label{F}
Let $\psi\in\Omega$ satisfying condition \eqref{our-cond}. Then we have 
$A_\psi(\alpha)=\alpha^{\log A_\psi(e)}$ and in particular
$$\lim_{\alpha\downarrow 1}A_\psi(\alpha)=1.$$
\end{lemma}
\noindent
\begin{proof}
Observe first that if $A_\psi(\alpha)$ exists for all $\alpha>1$, then it exists for all $\alpha>0$.
Indeed, we clearly have $A_\psi(1/\alpha)=1/A_\psi(\alpha)$. Analogously, we also observe that 
$A_\psi(\alpha\beta)=
A_\psi(\alpha)A_\psi(\beta)$ for all $\alpha,\beta>0$. 
These two observations entail that $A_\psi(\alpha)=\alpha^{\log A_\psi(e)}$
for all $\alpha\in\exp\{\mathbb Q\}$. As the  function 
$[\alpha\mapsto A_\psi(\alpha)]$ is  increasing, we conclude from what precedes that
it is also continuous  and thus $A_\psi(\alpha)=\alpha^{\log A_\psi(e)}$ for all 
$\alpha>0$.
\hfill\end{proof}
\begin{defn}
\label{kpsi}
In the sequel, for $\psi\in\Omega$ satisfying condition \eqref{our-cond}, we set
$$
k_\psi:=\log\big(A_\psi(e)\big),
$$
and we observe that we always have  $A_\psi(e)\geq 1$ and thus $k_\psi\geq 0$. 
\end{defn}

\begin{rem}
 Note also the equivalence of the following two conditions:
$$
\forall\alpha>1\,,\quad \lim_{t\to\infty}\frac{\psi(t^\alpha)}{\psi(t)}=1\qquad\Longleftrightarrow\qquad
\exists\,\alpha>1\;:\; \lim_{t\to\infty}\frac{\psi(t^\alpha)}{\psi(t)}=1.
$$
Indeed, assume that there exists $\alpha_0>1$ such that 
$\lim_{t\to\infty}\frac{\psi(t^{\alpha_0})}{\psi(t)}=1$. Then for  $\alpha\in(1,\alpha_0]$ and $t\geq1$,
we deduce
$$
1\leq\frac{\psi(t^\alpha)}{\psi(t)}\leq \frac{\psi(t^{\alpha_0})}{\psi(t)},
$$
which implies that $\lim_{t\to\infty}\frac{\psi(t^{\alpha})}{\psi(t)}=1$. For $\alpha\in(\alpha_0,\infty)$,
fix $\beta\in(1,\alpha_0]$ and $n\in\mathbb N$ such that $\beta^n\geq\alpha$. Then we get
$$
1\leq\frac{\psi(t^\alpha)}{\psi(t)}\leq\frac{\psi(t^{\beta^n})}{\psi(t)}=
\frac{\psi(t^{\beta^n})}{\psi(t^{\beta^{n-1}})}\,\frac{\psi(t^{\beta^{n-1}})}{\psi(t^{\beta^{n-2}})}\,
\dots\,\frac{\psi(t^{\beta})}{\psi(t)},
$$
and the assertion follows from the preceding argument.
\end{rem}
\begin{example}
\label{ex}
Consider  the functions
\begin{equation}
\label{ex1}
\psi_{n,\beta}(t):= \big(\log(1+\log(1+\cdots\log(1+t^{1/\beta})\dots))\big)^{\beta}\,, \qquad \beta>0\,,\quad
\mbox{$n$ iterations},
\end{equation}
 and\footnote{ This example generalizes \cite[Example 6]{SSZ}.} for a suitable constant $C>0$
\begin{equation}
\label{ex2}
\tilde\psi_{n,\beta}(t):=\begin{cases} C^{-1}\, t\,,\quad &t\in[0,C]\\
\exp\big\{\big(\log(\dots\log( t))\dots\big)^\beta\big\}\,,\quad &t>C
\end{cases},\qquad\beta\in(0,1)\,,\quad
\mbox{$n$ iterations}.
\end{equation}
It is clear that  all these functions belong to  $\Omega_b$. Moreover, $\psi_{n,\beta}$ satisfies
condition \eqref{our-cond}, with $k_\psi=\beta$ when $n=1$ and $k_\psi=0$ when
$n>1$. Also, $\tilde\psi_{n,\beta}$ satisfies
condition \eqref{our-cond} for $n>1$ only, with  $k_\psi=0$.
\end{example}
The Lorentz spaces associated to the   elements of $\Omega_b$ given above appear naturally in many  contexts.
Indeed,  with $\psi_{1,2}(t)=\log(1+t^{1/2})^2$, the space  $\clM_{\psi_{1,2}}$
is a
central tool in \cite{N} for the construction of trace functionals for holomorphic families of 
pseudo-differential operators on $\R^n$. (There $\clM_{\psi_{1,2}}$ is denoted by  
$\clL_{\log}^{(1,\infty)}$.)
Also, $\clM_{\psi_{1,n}}$, $\psi_{1,n}(t)=\log(1+t^{1/n})^n$ $n\in\N$, is the natural receptacle of the
embedding of $\clM_{\psi_{1,1}}(\clN_1,\tau_1)\otimes\dots\otimes\clM_{\psi_{1,1}}(\clN_n,\tau_n)$
 into $\clN_1\otimes\dots\otimes\clN_n$. Indeed, let $(\clN_j,\tau_j)$, $j=1,2$, be two semifinite
von Neumann algebras with normal semifinite faithful traces. From the relation:
$$
\int_0^t\mu(A_1\otimes A_2,s)\,ds=\sup_{t_1t_2=t}\, \int_0^{t_1}\mu(A_1,s_1)\,ds_1\;
 \int_0^{t_2}\mu(A_2,s_2)\,ds_2,
$$
we deduce that for any $\psi_1,\psi_2\in\Omega$
\begin{align}
\label{ellele}
\clM_{\psi_1}(\clN_1,\tau_1)\otimes \clM_{\psi_2}(\clN_2,\tau_2)\subset 
\clM_{\psi}(\clN_1\otimes\clN_2,\tau_1\otimes\tau_2)\quad\mbox{where}
\quad\psi(t)=\sup_{t_1t_2=t}\psi_1(t_1)\psi_2(t_2).
\end{align}
Now, the concavity of the function $[s\mapsto\log(\log(1+e^s))]$, implies that for every $n\in\N$
and $t_1,t_2>0$
$$
\log\big(1+t_1^{1/n}\big)^n\log(1+t_2)\leq \log\big(1+(t_1t_2)^{1/(n+1)}\big)^{n+1}.
$$
This, together with an inductive argument,   leads us to
$$
\sup_{t_1t_2=t}\log\big(1+t_1^{1/n}\big)^n\,\log(1+t_2)=\log\big(1+t^{1/(n+1)}\big)^{n+1},
$$
and thus
$$
\clM_{\psi_{1,1}}(\clN_1,\tau_1)\otimes\dots\otimes\clM_{\psi_{1,1}}(\clN_n,\tau_n)\subset \clM_{\psi_{1,n}}
(\clN_1\otimes\dots\otimes\clN_n,\tau_1\otimes\dots\otimes\tau_n).
$$
 Finally, we prove in Corollary \ref{OUOU} that if $0\leq T\in\clN$ satisfies
$$
{\tau\big(e^{-tT^{-1}}\big)} \sim C (T)\,t^{-1}\,|\log t|^{-1},\quad t\downarrow 0,
$$
then  $T$ belongs to $\mathcal \clM_{\psi_{2,1}}$ with $\psi_{2,1}(t)=\log(1+\log(1+t))$.

\quad

\noindent
Our next result establishes a hierarchy in the conditions \eqref{our-cond}, \eqref{SSZ-cond} and
 \eqref{ex-sing-trace-easy}.
\begin{prop}
\label{hier}
Condition \eqref{our-cond} is strictly stronger than condition \eqref{SSZ-cond} and the latter 
is strictly stronger than condition \eqref{ex-sing-trace-easy}.
\end{prop}
\noindent
\begin{proof}
First, we prove that condition \eqref{our-cond}  implies condition 
\eqref{ex-sing-trace-easy}. So, fix $a>1$.
Since $\psi$ is increasing, we have
$$
1\leq\liminf_{t\to\infty}\frac{\psi(at)}{\psi(t)}.
$$
Then, choose  $\alpha>1$ arbitrary. Since $at\leq t^\alpha$ for $t\geq a^{1/(\alpha-1)}$,
we deduce by  Lemma \ref{F}:
$$
\limsup_{t\to\infty}\frac{\psi(at)}{\psi(t)}\leq \limsup_{t\to\infty} \frac{\psi(t^\alpha)}{\psi(t)}=
\lim_{t\to\infty} \frac{\psi(t^\alpha)}{\psi(t)}=
\alpha^{k_\psi}.
$$
As $\alpha>1$ is arbitrary, we get
$$
\limsup_{t\to\infty}\frac{\psi(at)}{\psi(t)}\leq1,
$$
and the claim follows. Next, observe that
$$
1\leq\liminf_{t\to\infty} \frac{\psi\big(t\psi(t)\big)}{\psi(t)}.
$$
Then, choose  $\eps>0$ arbitrary. By \cite[Lemma 4.1]{CRSS}, 
condition \eqref{ex-sing-trace-easy} implies that $\psi(t)\leq C_\eps\, t^\eps$. 
Thus from
Lemma \ref{F} and what we have already proven, we get:
$$
\limsup_{t\to\infty} \frac{\psi\big(t\psi(t)\big)}{\psi(t)}\leq \limsup_{t\to\infty}
\frac{\psi(C_\eps t^{1+\eps})}{\psi(t)}
=\lim_{t\to\infty}\frac{\psi(C_\eps t^{1+\eps})}{\psi(t^{1+\eps})}\frac{\psi( t^{1+\eps})}{\psi(t)}= 
(1+\eps)^{k_\psi}.
$$
As $\eps>0$ is arbitrary, we deduce that
$$
\limsup_{t\to\infty} \frac{\psi\big(t\psi(t)\big)}{\psi(t)}\leq1,
$$
and thus condition  \eqref{our-cond}  implies condition  \eqref{SSZ-cond}. 
That condition  \eqref{SSZ-cond}  implies condition 
\eqref{ex-sing-trace-easy} is obvious. To see that all these implications are strict, we consider the 
 family of elements $\psi_\beta:=\tilde\psi_{1,\beta}\in\Omega_b$, $\beta\in(0,1)$,
where $\tilde\psi_{n,\beta}$  is defined in \eqref{ex2}.
Condition \eqref{ex-sing-trace-easy} is satisfied for all 
$\beta\in(0,1)$. Indeed, for $t>1$ and $a>1$, we have
$$
\frac{\psi_\beta(at)}{\psi(t)}=\exp\big\{\log( t)^\beta\big(\big(1+\log(a)/\log(t)\big)^\beta-1\big)\big\}
\sim\exp\big\{\beta\log(a)\log( t)^{\beta-1}\big\}\to 1,\quad t\to\infty.
$$
Condition \eqref{SSZ-cond} is satisfied for $\beta\in(0,1/2)$ but not for $\beta\in[1/2,1)$:
$$
 \frac{\psi_\beta\big(t\psi_\beta(t)\big)}{\psi_\beta(t)}= 
 \exp\big\{\log( t)^\beta\big(\big(1+\log(t)^{\beta-1}\big)^\beta-1\big)\big\}\sim
 \exp\big\{\beta\log(t)^{2\beta-1}\big\}\to
 \begin{cases}
 1,&\beta\in(0,\tfrac 12)\\
 \sqrt e,&\beta=\tfrac12\\
 \infty,&\beta\in(\tfrac12,1)
 \end{cases}\!\!\!,\quad t\to\infty.
$$
Condition \eqref{our-cond} is never satisfied. Indeed, for $\alpha>1$, we have:
$$
\frac{\psi_\beta\big(t^\alpha\big)}{\psi_\beta(t)}= \exp\big\{(\alpha^\beta-1)\log( t)^\beta\big\}\to\infty
,\quad t\to\infty.
$$
This concludes the proof.
\hfill\end{proof}

\subsection{Dixmier traces and $\psi$-compatible states}

\noindent
We say that a positive continuous linear functional $\varphi:\clM_\psi\to\mathbb C$ is  singular 
if its restriction to
$\clM_\psi^0$ is zero, symmetric if $\varphi(A)=\varphi(B)$ for all $A,B\in\clM_\psi^+$ such that
$\mu(A)=\mu(B)$ and  fully symmetric if moreover $\varphi(A)\leq\varphi(B)$ for all 
$A,B\in\clM_\psi^+$ such that $A\prec\prec B$. Note that every fully symmetric functional
is symmetric and bounded and if $\clM_\psi\not\subset\clL^1$, then every symmetric
functional is singular. 

\noindent
For $a>0$, we define 
the translation $T_a$, dilation $D_a$ and exponentiation $E_a$ 
operators on $L^\infty(\mathbb R_+^*)$ by
$$
T_a f(t):=f(t+a)\,,\quad D_a f(t):=f(at)\,,\quad E_a f(t):=f(t^a)\,,\qquad f\in L^\infty(\mathbb R_+^*).
$$
A continuous linear functional $\omega\in\big( L^\infty(\mathbb R_+^*)\big)^*$ is said to be  translation invariant, 
dilation invariant or exponentiation invariant if for all $f\in L^\infty(\mathbb R_+^*)$ and all $a>0$, 
we have
\begin{align}
\label{invariance}
\omega(T_a f)=\omega(f)\,,\quad\omega(D_a f)=\omega(f)\quad\mbox{or}\quad
\omega(E_a f)=\omega(f).
\end{align}

\begin{defn}
A state $0\leq\omega\in\big( L^\infty(\mathbb R_+^*)\big)^*$ is called  singular 
 if  it vanishes
on functions converging to zero at infinity.
\end{defn}
\noindent
We observe that translation, dilation and exponentiation invariant states are singular. 
Singular states may also be called extended limits, as they are extensions (via the Hahn-Banach
Theorem) of the ordinary limit functional on the linear subset of $L^\infty(\mathbb R_+^*)$
consisting of
functions converging at infinity.

\noindent
It is proven in
 \cite[Theorem 3.4]{DPSS},  that a given Lorentz  space $\clM_\psi$
possesses  a nontrivial fully symmetric  singular trace supported at infinity (i.e$.$ which vanishes
on $\clM_\psi\cap\clL^1$) if and only if
\begin{equation}
 \label{ex-sing-trace}
\exists \,a>1\quad :\quad\liminf_{t\to\infty}\frac{\psi(at)}{\psi(t)}=1.
\end{equation}
Thus, Proposition \ref{hier} shows that the Lorentz  spaces associated
with the elements of $\Omega$ satisfying condition \eqref{our-cond} do possess  nontrivial  singular
 traces supported at infinity.
\begin{rem}
\label{rem-eq-cond}
As in Lemma \ref{eq-conds-1}, the arguments of \cite[Theorem 6]{KSS} can be generalized  
to show that  condition \eqref{ex-sing-trace} is equivalent to any of the
following two conditions:
$$
\forall \,a>1\,,\;\;\liminf_{t\to\infty}\frac{\psi(at)}{\psi(t)}=1\,,\qquad \qquad
\liminf_{t\to\infty}\frac{t\psi'(t)}{\psi(t)}=0.
$$
\end{rem}

\noindent
The most important  class of fully symmetric singular traces is certainly the class of Dixmier traces.
The most general way to construct  Dixmier traces comes from \cite{KSS},  where we borrow
the following notion:
\begin{defn}
\label{def-psi-comp}
A  state $0\leq \omega\in\big( L^\infty(\mathbb R_+^*)\big)^*$ is said to be 
$\psi$-compatible  if there exists $a>1$ such that
$$
\omega\Big(\Big[t\mapsto \frac{\psi(at)}{\psi(t)}\Big]\Big)=1.
$$
\end{defn}
\noindent
It is proven in \cite[Proposition 10]{KSS} that given  $\omega$, a $\psi$-compatible and dilation 
invariant  state on
 $ L^\infty(\mathbb R_+^*)$, the (unitarily invariant and singular) functional
 \begin{equation}
\label{Dix-Trace}
\tau_{\psi,\omega}:\clM_\psi ^+\to \mathbb R_+\;,\qquad T\mapsto \omega\Big(\Big[t\mapsto
\frac1{\psi(t)}\int_0^t\mu(T,s)\,ds\Big]\Big),
\end{equation}
is positively additive (and positively homogeneous too), 
and its extension by linearity  on ${\clM_\psi} $ is called a Dixmier trace.
As proven in \cite[Proposition 10]{KSS}, it is a fully symmetric normalized tracial functional 
on ${\clM_\psi} $.
Surprisingly, the converse also holds true (see \cite[Theorem 11]{KSS}): 
every  fully symmetric 
normalized\footnote{A functional $\varphi$ on $\clM_\psi$ is normalized if $\varphi(T)=1$
for all positive $T\in\clM_\psi$ such that $\mu(T)=\psi'$.}
  tracial linear functional on ${\clM_\psi} $ is a  Dixmier trace. 

\begin{rem}
Observe that when $\psi\in\Omega$ satisfies condition \eqref{ex-sing-trace-easy}, then
every  singular
 state is $\psi$-compatible. Hence, Proposition \ref{hier} shows that it is still the case when
 condition \eqref{SSZ-cond} or  condition \eqref{our-cond} is satisfied.
\end{rem}

\noindent
The next result is strongly inspired from \cite[Proposition 5]{SUZ}. It  relates an analogous  notion of
$\psi$-compatibility in the context of exponentiation  invariant states, 
to the existence of Dixmier traces:
\begin{prop}
\label{lin}
Let $\omega$ be an exponentiation invariant state of $L^\infty(\mathbb R_+^*)$ 
and $\psi\in\Omega$,  such that
\begin{align}
\label{cond-mod}
\lim_{\alpha\,\downarrow \,1}\omega\Big(\Big[t\mapsto\frac{\psi(t^\alpha)}{\psi(t)}\Big]\Big)=1.
\end{align}
Then the functional
\begin{align}
\label{exp-exists}
\clM_\psi^+\to[0,\infty),\qquad T\mapsto \omega\Big(\Big[t\mapsto
\frac1{\psi(t)}\int_0^t\mu(T,s)\,ds\Big]\Big),
\end{align}
is positively additive. Since its extension by linearity is fully symmetric and normalized,
there exists  a $\psi$-compatible and
dilation invariant state $\omega'$, such that  this functional coincides with
the Dixmier trace $\tau_{\psi,\omega'}$.
\end{prop}
\noindent\begin{proof}
We just have to prove additivity  on $\clM_\psi^+$. So let $T,S\in \clM_\psi^+$.
By \cite[Theorem 4.4 (ii)]{FK}, we have
$$
\int_0^{t/2}\big( \mu(T,s)+\mu(S,s)\big)\, ds\leq \int_0^{t} \mu(T+S,s)\,ds\leq
\int_0^{t}\big( \mu(T,s)+\mu(S,s)\big)\, ds,\qquad\forall\,t>0,
$$
which yields by positivity of $\omega$:
$$
 \omega\Big(\Big[t\mapsto
\frac1{\psi(t)} \int_0^{t} \mu(T+S,s)\,ds\Big]\Big)\leq
\omega\Big(\Big[t\mapsto
\frac1{\psi(t)}\int_0^{t}\big( \mu(T,s)+\mu(S,s)\big)\, ds\Big)\Big],
$$
and
\begin{align*}
\omega\Big(\Big[t\mapsto
\frac1{\psi(t)}\int_0^{t}\big( \mu(T,s)+\mu(S,s)\big)\, ds\Big]\Big)&\leq
 \omega\Big(\Big[t\mapsto 
\frac1{\psi(t)} \int_0^{2t} \mu(T+S,s)\,ds\Big]\Big).
\end{align*}
Since for every $\eps>0$ and $t$ sufficiently large we have $2t\leq t^{1+\eps}$,
we deduce
$$
 \omega\Big(\Big[t\mapsto 
 \frac1{\psi(t)}\int_0^{t}\big( \mu(T,s)+\mu(S,s)\big)\, ds\Big]\Big)\leq
  \omega\Big(\Big[t\mapsto
\frac1{\psi(t)} \int_0^{t^{1+\eps}} \mu(T+S,s)\,ds\Big]\Big).
$$
On the other hand, for an arbitrary $0\leq g\in L^\infty(\mathbb R_+^*)$, as 
$\psi(t^{1+\eps})/\psi(t)\geq 1$ for $t\geq 1$ and $\eps>0$, we deduce from our hypotheses on $\omega$:
$$
0\leq  \omega\Big(\Big[t\mapsto \Big(\frac{\psi(t^{1+\eps})}{\psi(t)}-1\Big)g(t)\Big]\Big)\leq\|g\|_\infty
\;\omega\Big(\Big[t\mapsto \frac{\psi(t^{1+\eps})}{\psi(t)}-1\Big]\Big)\to0,\quad \eps\to 0.
 $$
 Since,
 $$
 \sup_{t>0}\frac1{\psi(t^{1+\eps})} \int_0^{t^{1+\eps}} \mu(T+S,s)\,ds=\|T+S\|_{\clM_\psi}\,,\quad\forall\eps>0,
 $$
 we deduce
 \begin{align*}
  \omega\Big(\Big[t\mapsto
\frac1{\psi(t)}\int_0^{t}\big( \mu(T,s)+\mu(S,s)\big)\, ds\Big]\Big)&\leq \|T+S\|_{\clM_\psi} \omega\Big(\Big[t\mapsto \frac{\psi(t^{1+\eps})}{\psi(t)}-1\Big]\Big)\\
&\quad+ \omega\Big(\Big[t\mapsto
\frac1{\psi(t^{1+\eps})} \int_0^{t^{1+\eps}} \mu(T+S,s)\,ds\Big]\Big).
\end{align*}
By exponentiation invariance of $\omega$, we obtain
$$
\omega\Big(\Big[t\mapsto
\frac1{\psi(t^{1+\eps})} \int_0^{t^{1+\eps}} \mu(T+S,s)\,ds\Big]\Big)
=\omega\Big(\Big[t\mapsto
\frac1{\psi(t)} \int_0^{t} \mu(T+S,s)\,ds\Big]\Big)
$$
Thus,  taking the limit $\eps\to 0$,
we get
$$
\omega\Big(\Big[t\mapsto\frac1{\psi(t)}\int_0^{t}\big( \mu(T,s)+\mu(S,s)\big)\, ds\Big]\Big)\leq
  \omega\Big(\Big[t\mapsto
\frac1{\psi(t)} \int_0^{t} \mu(T+S,s)\,ds\Big]\Big).
$$
The last statement follows from  \cite[Theorem 11]{KSS}.
\hfill\end{proof}
\begin{corollary}
Let $\psi\in\Omega$ satisfying condition \eqref{our-cond}. Then for any exponentiation
invariant state $\omega$, the linear extension of
the map \eqref{exp-exists} defines a Dixmier trace on $\clM_\psi$.
\end{corollary}
\noindent\begin{proof}
This follows from Lemma \ref{F}. 
\hfill\end{proof}

\begin{rem}
\label{comp-mod}
Observe that condition \eqref{cond-mod} is not satisfied for the elements $\tilde\psi_{1,\beta}$,
$\beta\in(0,1)$ given in \eqref{ex2}. Indeed, we have already seen that
 in this case  $\lim_{t\to\infty}\psi(t^\alpha)/
\psi(t)=\infty$ for all $\alpha>1$.
\end{rem}

\noindent
The goal of the next result is to clarify the relationship between the notion of  $\psi$-compatibility  
and  conditions
 \eqref{ex-sing-trace} and \eqref{ex-sing-trace-easy}. The first part of this result is implicitly contained
 in \cite{KSS} and largely discussed in  \cite{LSZ}. For convenience of the reader, we supply the detailed argument.
 \begin{prop}
 \label{psi-adm}
 Let $\psi\in\Omega_b$. \\
 $\rm (i)$ Condition  \eqref{ex-sing-trace} is satisfied if and only if there exists a $\psi$-compatible 
 dilation invariant state.\\
$\rm (ii)$  Condition  \eqref{ex-sing-trace-easy} is satisfied if and only if  every translation 
invariant state  is $\psi$-compatible.
\end{prop}
\noindent
\begin{proof}
By Remark \ref{rem-eq-cond}, we may assume that $a=2$ in condition  \eqref{ex-sing-trace}.
If a $\psi$-compatible dilation invariant state $\omega$ exists, then by \cite[Proposition 10]{KSS}
the map $\tau_{\psi,\omega}$ is a fully symmetric  singular trace on $\clM_\psi$,
 which by \cite[Theorem 3.4]{DPSS} implies (since $\psi\in\Omega_b$ is equivalent to  
 $\clM_\psi\subset\clN$
 and thus nontrivial singular traces are necessarily supported at infinity) that $\liminf_{t\to\infty}\psi(2t)/\psi(t)
 =1$. (See also the discussion prior to Definition 9 in \cite{KSS} for a direct proof of this fact.)
Now assume that condition  \eqref{ex-sing-trace} holds. Then, by \cite[Lemma 6.3.4]{LSZ}, there
exists a dilation invariant state $\omega\in \big(L^\infty(\mathbb R_+^*)\big)^*$, such that
$\omega\big(\big[t\mapsto t\psi'(t)/\psi(t)\big]\big)=0$. Arguing as in Lemma \ref{eq-conds-1}, we see that 
$\omega$ is $\psi$-compatible as well.
 This proves the first equivalence.
We next prove the second equivalence. 
Clearly, if the condition \eqref{ex-sing-trace-easy}  holds,
 then every singular state is  $\psi$-compatible. In particular 
every translation invariant state is $\psi$-compatible. 
Next assume that every translation invariant state is $\psi$-compatible.
By \cite[Theorem 6]{KSS}, the condition of $\psi$-compatibility of Definition \ref{def-psi-comp}
is equivalent to $\omega\big(\big[t\mapsto {\psi(at)}/{\psi(t)}\big]\big)=1$, for all $a>1$. By Lemma
\ref{eq-conds-1}, we can therefore prove condition  \eqref{ex-sing-trace-easy} with $a=2$.
 This means
that the map $[t\mapsto \psi(2t)/\psi(t)]$ is almost convergent (in the sense of Lorentz
\cite{Lor}). By \cite[Theorem 3.3]{LSS}, this implies that
\begin{align}
\label{LSSt}
\lim_{u\to\infty}\frac1u\int_0^u \frac{\psi(2t)}{\psi(t)}\,dt=1.
\end{align}
Last, by \cite[Section 6.8]{H} if there exists $C>0$ such that
\begin{align}
\label{Ht}
t\frac d{dt}\Big(\frac{\psi(2t)}{\psi(t)}\Big)\geq -C,
\end{align}
then \eqref{LSSt} is equivalent to \eqref{ex-sing-trace-easy}. Hence, we just need to
prove the inequality \eqref{Ht} to conclude the proof. For this,  we remark that
$$
t\frac d{dt}\Big(\frac{\psi(2t)}{\psi(t)}\Big)=t\frac{2\psi'(2t)\psi(t)-\psi(2t)\psi'(t)}{\psi(t)^2}
\geq -\frac{t\psi(2t)\psi'(t)}{\psi(t)^2}=-\frac{t\psi'(t)}{\psi(t)}\frac{\psi(2t)}{\psi(t)},
$$
and  observing that by the concavity of $\psi$, we have:
$$
0\leq\frac{t\psi'(t)}{\psi(t)}\leq 1,\qquad 1\leq \frac{\psi(2t)}{\psi(t)}\leq 1+\frac{t\psi'(t)}{\psi(t)}
\leq 2,
$$
the result follows with $C=2$.
\hfill\end{proof}

\noindent
As the property \eqref{ex-sing-trace-easy} also implies that 
 every dilation
invariant state of $L^\infty(\mathbb R_+^*)$ is $\psi$-compatible, we deduce the following:
\begin{cor}
 If  every translation 
invariant state of $L^\infty(\mathbb R_+^*)$ is $\psi$-compatible, then  every dilation
invariant state of $L^\infty(\mathbb R_+^*)$ is $\psi$-compatible.
\end{cor}

\noindent 
We next show that the existence of dilation invariant and $\psi$-compatible states is
also related to the value of the Boyd indices of the associated Lorentz spaces.

\begin{prop}
\label{la-tres-belle}
Let $\psi\in\Omega_b$. \\
$\rm (i)$  There exists a dilation invariant and
$\psi$-compatible state if and only if $\mathfrak p_{\clM_\psi}=1$.\\
$\rm (ii)$ Every dilation invariant state is
$\psi$-compatible if and only if $\mathfrak p_{\clM_\psi}=\mathfrak q_{\clM_\psi}=1$.
\end{prop}
\noindent
\begin{proof}
For the first part, we observe   by \cite[Proposition 2.3]{AS}, that condition \eqref{ex-sing-trace} is
equivalent to $\mathfrak p_{\clM_\psi}=1$. One then concludes with Proposition \ref{psi-adm} (i).\\
For the second assertion, we start by  \cite[Theorem 8]{KSS},  which shows that
every dilation-invariant state is $\psi$-compatible if and only if 
\begin{align}
\label{every}
\forall \eps>0\,,\quad\exists C>0\,\,:\,\,\forall t>0,\,\forall s>1\,,\quad \psi(st)\leq C s^\eps\psi(t).
\end{align}
On the other hand, we have by definition of the upper Boyd index of $\clM_\psi$:
$$
\mathfrak q_{\clM_\psi}=\lim_{a\to \infty}\frac{\log a^{-1}}{\log\|D_a\|_{\clM_\psi\to\clM_\psi}}.
$$
Thus $\mathfrak q_{\clM_\psi}=1$ if and only if $\|D_a\|_{\clM_\psi\to\clM_\psi}=a^{-1+o(1)}$, 
 $a\to \infty$. Next, we observe that 
\begin{equation}
\label{recall}
\|D_a\|_{\clM_\psi\to\clM_\psi}=\sup_{t>0}\frac{\psi(at)}{a\psi(t)}.
\end{equation}
Indeed, as
 $\psi'$ belongs to the unit sphere of $\clM_\psi(\mathbb R_+)$ and as $\psi'$ is decreasing, 
we first deduce:
$$
\|D_a\|_{\clM_\psi\to\clM_\psi}\geq\|D_a\psi'\|_{\clM_\psi}.
$$
Also, since for an arbitrary $f\in\clM_\psi(\mathbb R_+)$ we have $\mu(f)\prec\prec \|f\|_{\clM_\psi}\psi'$
and $\mu(D_af)=D_a\mu(f)$, 
we obtain
$$
\frac{\|D_af\|_{\clM_\psi}}{\|f\|_{\clM_\psi}}\leq  \|D_a\psi'\|_{\clM_\psi},
$$
and \eqref{recall} follows from the computation:
$$
\|D_a\psi'\|_{\clM_\psi}=\sup_{t>0} \frac1{\psi(t)}\int_0^t\psi'(as)\,ds=
\sup_{t>0}\frac{\psi(at)}{a\psi(t)}.
$$

 Hence $\mathfrak q_{\clM_\psi}=1$ if and only if
   $$
  \sup_{t>0}\frac{\psi(at)}{\psi(t)}=a^{o(1)}, \quad a\to \infty,
  $$
which is equivalent to
  \begin{align}
\label{everybis}
\forall \eps>0\,,\quad\exists a>1\,\,:\,\,\forall t>0,\quad \psi(at)\leq  a^\eps\psi(t).
\end{align}
 All that remains to do is to show that the conditions \eqref{every} and \eqref{everybis}
 are equivalent. We first assume that condition  \eqref{every} is satisfied.
 Then fix $\eps>0$ and let $C(\eps/2)$ be the constant given in condition \eqref{every}
 for $\eps/2$ (we may assume this constant to be strictly greater than one, since in the opposite case
 the conclusion is obvious). Define then
 $$
 a(\eps):=C(\eps/2)^{2/\eps}>1.
 $$
 We then get that for all $t>0$
 $$
 \psi\big(a(\eps)t\big)\leq C(\eps/2) a(\eps)^{\eps/2}\psi(t)=  C(\eps/2)^2\psi(t)=a(\eps)^\eps\psi(t),
 $$
 which is condition \eqref{everybis}. Next, assume that condition \eqref{everybis} holds.
 Fix $\eps>0$ and let $a(\eps)>1$ be as given by condition \eqref{everybis}. For an arbitrary $s>1$,
 let $n\in\mathbb N$ be such that $s\in[a(\eps)^n,a(\eps)^{n+1})$. By monotony of $\psi$ and
 using repeatedly condition \eqref{everybis}, we get  for arbitrary $t>0$:
 $$
 \psi(st)\leq\psi(a(\eps)^{n+1}t)\leq a(\eps)^\eps\psi(a(\eps)^{n}t) \leq \dots
 \leq a(\eps)^{(n+1)\eps}\psi(t)\leq a(\eps)^\eps s^\eps\psi(t),
 $$
 which is condition \eqref{every} with constant $C(\eps)=a(\eps)^\eps$.
\hfill
\end{proof}

\noindent
In the last result of this section, we generalize \cite[Proposition 4.3]{CRSS}:
\begin{prop}
\label{calc}
Let $\omega$ be an exponentiation invariant state of $L^\infty(\mathbb R_+^*)$ 
and $\psi\in\Omega$.  If condition  \eqref{cond-mod}  holds and if \begin{equation}
\label{weak}
\forall \eps\in(0,1)\,,\;\exists C>0\;:\; \forall t>1\,,\; \psi(t)< C\, t^\eps,
\end{equation}
then for every $T\in\clM_\psi^+$, we have
$$
\tau_{\psi,\omega}(T)=\omega\Big(\Big[t\mapsto  \frac{\tau\big(T\chi_{(\frac1t,\infty)}(T)\big)}{\psi(t)}
\Big]\Big)=
\omega\Big(\Big[t\mapsto \frac{\tau\big(T\chi_{(\frac1t,1)}(T)\big)}{\psi(t)}\Big]\Big).
$$
\end{prop}
\noindent
\begin{proof}
Under the condition \eqref{weak} and by exactly the
same arguments than those of \cite[Proposition 4.3]{CRSS}, 
we get that for any $\eps>0$ and $t>0$
sufficiently large:
\begin{align*}
\int_0^t\mu(T,s)\,dt\leq \tau\big(T\chi_{(\frac1t,\infty)}(T)\big)+C\leq \int_0^{t^{1+\eps}}\mu(T,s)\,dt+C,
\end{align*}
  for an absolute constant $C>0$.
Dividing these inequalities by $\psi(t)$, applying $\omega$ and using the same argument
as in Proposition \ref{lin}, we get the result. The argument for the second equality is identical 
to those of \cite[Corollary 4.4]{CRSS}.
\hfill\end{proof}
\begin{rem}
\label{ouf}
Gathering  Remark \ref{comp-mod}, Proposition \ref{hier} and \cite[Lemma 4.1]{CRSS}, we deduce
that the condition \eqref{weak} is satisfied for $\psi\in\Omega$ verifying 
condition \eqref{our-cond}.
\end{rem}

\section{An extrapolation description of Lorentz  spaces close to $\clL^1$}
\label{2}

This section contains our chief innovation to study singular traces on generic Lorentz spaces. 
It relies on a characterization of Lorentz spaces
which can be described as an extrapolated space.
Our results  (partly) extend  \cite{AL1,AL2} in the noncommutative 
setting and  for Lorentz spaces which are close to $\clL^1$.
(The authors of  \cite{AL1,AL2} study the commutative and dual  situation of Lorentz spaces 
which are close to $L^\infty$.)
For that, we 
construct an extrapolation method, $\mathfrak L$, which assigns to each  Banach lattice ideals of 
measurable functions on the interval
 $(1,\infty)$, a fully symmetric operator spaces on any semifinite von Neumann 
algebra.  
We apply our extrapolation method to a pair $(F_\psi,F^\psi)$ of Banach spaces of functions on
$(1,\infty)$, canonically 
associated to any element $\psi\in\Omega$ and call $(\mathfrak L_\psi,\mathfrak L^\psi)$
the resulting pair of fully symmetric operator spaces. The first important result of this section is 
Proposition \ref{approx}, where we prove that one always has
$ \mathfrak L_\psi\subset\clM_\psi\subset \mathfrak L^\psi$ with continuous inclusions. We then give
different characterizations of when  $\clM_\psi= \mathfrak L^\psi$ and $ \mathfrak L_\psi=\clM_\psi$  
(Proposition \ref{eq-conds} and Proposition \ref{conds-2}).
We conclude by exhibiting a convenient  sufficient  condition to have the equality 
$ \mathfrak L_\psi=\clM_\psi$, which in turn, is one of the the main step in the  analysis we perform in 
section \ref{3}  (Proposition \ref{prop-eq}).

\subsection{Lorentz  spaces close to $\clL^1$}
\begin{defn}
\label{close-p0}
A  symmetric operator space
 is said to be close to $\clL^1$, if  it is continuously embedded in $\clL^p$ for all $p> 1$,
and if it is  not contained   in $\clL^1$.
\end{defn}

\noindent
In the context of  Lorentz  spaces, we have an easy characterization of those which are 
close to $\clL^1$:
\begin{lemma}
\label{eq-close}
Let $\psi\in\Omega$. Then $\clM_\psi$ is close to $\clL^{1}$ if and only if
$\psi'\in L^p(\mathbb R_+)$ for all  $p>1$ and $\psi'\not\in L^1(\mathbb R_+)$.
\end{lemma}
\noindent 
\begin{proof}
Note first  that for $T\in\clM_\psi$ and since $\psi(0)=0$:
$$
\int_0^t\mu(T,s)\, ds\leq\|T\|_{\clM_\psi}\,\psi(t)=\|T\|_{\clM_\psi}\,\int_0^t\psi'(s)\,ds,
\quad \forall t>0.
$$
Thus 
$
\mu(T)\prec\prec\|T\|_{\clM_\psi}\,\psi'\;,
$
 which, by \cite[Chapter II, Lemma 3.4]{GK},  implies that for all $p\geq 1$,
$
\mu(T)^{p}\prec\prec \|T\|_{\clM_\psi}^{p}\,\psi'^{p},
$
 that is
$$
\int_0^t\mu(|T|^p,s)\, ds=\int_0^t\mu(T,s)^p \,ds\leq\|T\|_{\clM_\psi}^p\,\int_0^t\psi'(s)^p\,ds\,,\quad\forall t>0.
$$
So taking the limit $t\to \infty$ for $p>1$, we deduce:
\begin{equation}
\label{inequality1}
\|T\|_{p}\leq\|T\|_{\clM_\psi} \,\|\psi'\|_{p}\;,\qquad\forall \,T\in\clM_\psi.
\end{equation}
Assume first that $\clM_\psi$ is close to $\clL^1$.  We then observe 
that $\psi'\in L^p(\mathbb R_+)$ for all $p>1$ as $\psi'\in\clM_\psi(\mathbb R_+)$. Now, as $\clM_\psi\not\subset\clL^1$,
there exists $T_0\in\clM_\psi$ such that $\|T_0\|_1=\infty$. Taking the limit $p\to 1$ in
\eqref{inequality1}, gives
$$           
\lim_{p\to 1}\|\psi'\|_p\geq\lim_{p\to 1}\frac{\|T_0\|_{p}}{\|T_0\|_{\clM_\psi}}=\infty,
$$
and the first implication is proven. Next, assume the second condition. Then the inequality
\eqref{inequality1} shows that $\clM_\psi\subset \clL^p$ for all $p>1$. But 
 $\clM_\psi\not\subset \clL^1$ as $\psi'\in\clM_\psi(\mathbb R_+)$ 
 and $\psi'\notin L^1(\mathbb R_+)$.
 This proves the second implication.
\hfill\end{proof}
\begin{rem}
In the light of the previous result, it is  clear that the Lorentz spaces associated with 
the elements of $\Omega_b$ given in  \eqref{ex1} and \eqref{ex2}, are close to $\clL^1$.
\end{rem}

 \begin{lemma}
 \label{close-L1}
 Let $\psi\in\Omega$  be such that the property \eqref{weak} holds.
Then $\clM_\psi$ is close to $\clL^1$.
\end{lemma}
\noindent
\begin{proof}
Observe that since $\psi(0)=0$, we have
$\|\psi'\|_1=\lim_{t\to\infty}\psi(t)=\infty$ and thus
$\psi'\notin L^1(\mathbb R_+)$.
By the concavity of $\psi$, we deduce that
$t\psi'(t)\leq \psi(t)$ and thus 
 for $t>1$, we get $\psi'(t)\leq C(\eps)t^{-1+\eps}$ for all $\eps\in(0,1)$.
Hence $\psi'\in L^p(\mathbb R_+)$ for all $p>(1-\eps)^{-1}$. As $\eps\in(0,1)$ is
arbitrary, $\psi'\in L^p(\mathbb R_+)$ for all $p>1$. The proof then follows by Lemma 
\ref{eq-close}.
\hfill\end{proof}
\noindent
Since the  property  \eqref{weak} is implied by \eqref{every}, and that \eqref{every}
is equivalent to  the $\psi$-compatibility of every dilation invariant state of $L^\infty(\mathbb R_+^*)$, 
we  deduce by Proposition \ref{la-tres-belle} (ii):
\begin{cor}
\label{ew}
Let $\psi\in\Omega_b$ be 
such that $\mathfrak p_{\clM_\psi}=\mathfrak q_{\clM_\psi}=1$, then $\clM_\psi$ 
is close to $\clL^1$. In particular,  this happens
when condition \eqref{our-cond} holds.
\end{cor}

\subsection{An extrapolation method}
\noindent
Consider $(F,\|.\|_F)$,  a Banach ideal lattice of (class of) measurable
functions  on the open interval $(1,\infty)$. 
Given $T\in\clL_0$, we define $\eta_T$ to be the $[0,\infty]$-valued measurable function 
on $(1,\infty)$ given by
\begin{equation}
\label{eta}
\eta_T:=\big[p\in(1,\infty)\mapsto\|T\|_p\big].
\end{equation}
We then define a Banach space of measurable operators by
\begin{equation}
\label{LF}
 \mathfrak L_F:=\big\{T\in\clL_0\,:\,\eta_T\in F\big\}.
\end{equation}
It is obvious that  $\mathfrak L_F$  is a bi-module over $\clN$ and that it shares most of the
 properties of $F$. 
More specifically, we have:
\begin{lemma}
\label{basic-prop}
 Let $F$ be a Banach ideal lattice   on the interval $(1,\infty)$.
Normed with
$$
\|T\|_{\mathfrak L_F}:=\big\|\eta_T\big\|_F,
$$
 the  operator space $\mathfrak L_F$, becomes a fully symmetric
 operator space and it has  the Fatou property if and only if  $F$ does.
 Moreover, if the evaluation maps
 $$
F\to\mathbb C\,,\quad f\mapsto f(p)\,,\qquad p\in(1,\infty),
 $$
 are continuous, then $\mathfrak L_F$ embeds continuously in $\clL^p$, for all $p\in(1,\infty)$.
\end{lemma}
\noindent
\begin{proof}
First note that the $\clL^p$-spaces are symmetric and
have the Fatou property. Hence, $\mathfrak L_F$ has
 the Fatou property if and only if  $F$ does. Next observe that if $T\in\clL_0$ and $S\in \mathfrak L_F$
are such that $T\prec\prec S$, then for $p\geq 1$,  $|T|^p\prec\prec |S|^p$ and thus for $p
>1$, $\eta_T(p)=\|T\|_p\leq\|S\|_p=\eta_S(p)$, which entails that $T\in\mathfrak L_F$,
with  $\|T\|_{\mathfrak L_F}
\leq\|S\|_{\mathfrak L_F}$, since $F$ is   Banach lattice. This proves the first claim.
Last, if for all $p>1$, there exists
a constant $C_p>0$ such that for all $f\in F$ we have $|f(p)|\leq C_p\|f\|_F$, we deduce that
for all $T\in \mathfrak L_F$, we have $\|T\|_p=\eta_T(p)\leq C_p\,\|\eta_T\|_F=
C_p\,\|T\|_{\mathfrak L_F}$, proving the second claim.
\hfill\end{proof}
\noindent
Next, we turn to the notion of extrapolated 
 operator  spaces close to $\clL^1$.
\begin{defn}
\label{def-extra}
 Let $\clE$  be a fully symmetric operator  space. 
Then $\clE$ is called an {\em extrapolated operator   space}  if $\clE=\mathfrak L_F$ with 
equivalent norms, 
 where $F$ is a Banach ideal lattice of measurable functions
 on  the interval $(1,\infty)$, which contains the constant unit function and such that  the evaluation maps $ F\to\mathbb C$, $f\mapsto f(p)$, $p>1$, are continuous. 
We let $\mathcal X_1$ be
the  class of all extrapolated operator spaces close to $\clL^1$.
\end{defn}

\noindent
We will see in Corollary \ref{obs} that  
the class $\mathcal X_1$ contains the standard Lorentz space $\clM^{1,\infty}$ and, more generally,
the Lorentz spaces associated with the elements  of $\Omega_b$ given in \eqref{ex1}.
Our next observation is that an extrapolated operator space 
close to $\clL^1$ has always trivial Boyd indices:
\begin{prop}
\label{BI}
Let $\mathcal E$ be a fully symmetric operator  space such that $\clE\in \mathcal X_1$ and such that
 $\clE\subset\clN$ continuously. Then 
$\mathfrak p_\clE=\mathfrak q_\clE=1$.
\end{prop}
\noindent
\begin{proof}
By assumption, there exists a Banach lattice $F$ on the interval $(1,\infty)$ such that
$\clE=\mathfrak L_F$ with equivalent norms. Thus, the Boyd indices of $\clE$ and of 
$\mathfrak L_F$ coincide.  Moreover, by assumption too, $\mathfrak L_F\subset\clN$
continuously. This implies that $\mathfrak L_F=\tilde{\mathfrak L}_F$ with equivalent norms, where
$$
\tilde{ \mathfrak L}_F:=\big\{T\in\clN\,:\,\tilde\eta_T\in F\big\}\quad\mbox{and}
 \quad\tilde\eta_T(p):=\eta_T(p)+\|T\|.
$$
Indeed since $\eta_T\leq\tilde\eta_T$ and since $F$ is a Banach lattice, for $T\in \tilde{ \mathfrak L}_F$ we have 
$$
\|T\|_{{ \mathfrak L}_F}=\|\eta_T\|_F\leq\|\tilde\eta_T\|_F=\|T\|_{\tilde{ \mathfrak L}_F}.
$$
On the other hand, we know by assumption that there exists $C>0$ such that for
all $T\in { \mathfrak L}_F$, we have $\|T\|\leq C\,\|T\|_{{ \mathfrak L}_F}$. This 
entails that
$$
\|T\|_{\tilde{ \mathfrak L}_F}=\|\tilde\eta_T\|_F\leq \|\eta_T\|_F+\|T\|\,\|1\|_F\leq\big(
1+C\,\|1\|_F\big)\|T\|_{{ \mathfrak L}_F},
$$
which is enough to conclude as the constant unit function belongs to $F$ by
assumption (see Definition
\ref{def-extra}). Thus, the Boyd indices of $\clE$ and of $ \tilde{ \mathfrak L}_F$
coincide as well.
Now, fix $q>1$. We consider the closed  subspace 
$F_q$ of $F$ defined by
$$
F_q:=\{f \chi_{(1,q]}\,:\, f\in F\},
$$
where $ \chi_{(1,q]}$ is the indicator function of the interval $(1,q]$.
Equipped with the norm
$$
\|f\|_{F_q}:=\|f \chi_{(1,q]}\|_{F},
$$
$F_q$ becomes a Banach lattice as well. Moreover, $\tilde{\mathfrak L}_F
=\tilde{\mathfrak L}_{F_q}$
with equivalent norms. Indeed as $F$ is a Banach lattice, for $T\in \tilde{\mathfrak L}_{F}$, 
we have
$$
\|T\|_{\tilde{\mathfrak L}_{F_q}}=\|\tilde\eta_T\|_{F_q}=\|\tilde\eta_T\chi_{(1,q]}\|_{F}\leq 
\|\tilde\eta_T\|_F=\|T\|_{\tilde{\mathfrak L}_{F}},
$$
 while for $T\in \tilde{\mathfrak L}_{F_q}$, we deduce 
$$
\|T\|_{\tilde{\mathfrak L}_{F}}=\|\tilde\eta_T\|_{F}\leq\|\tilde\eta_T\chi_{(1,q]}\|_{F}+\|\tilde\eta_T\chi_{(q,\infty)}\|_{F}
\leq \|T\|_{\tilde{\mathfrak L}_{F_q}}+\|\chi_{(q,\infty)}\|_{F}\sup_{p>q}\tilde\eta_T(p).
$$
Next, we observe that for $p>q$, the Young inequality yields:
$$
\tilde\eta_T(p)=\|T\|+\|T\|_p\leq\|T\|+\|T\|^{1-q/p}\|T\|_{q}^{q/p}\leq 2\|T\|+\|T\|_{q}
\leq 2\tilde\eta_T(q).
$$
Observe also that as the evaluation map is continuous on $F$ (by assumption)
and as $F_q\subset F$,
the evaluation map is also continuous on $F_q$ for its own topology. Indeed, for $f\in F_q$
and $p\in (1,q]$ we have
$$
f(p)=f(p)\chi_{(1,q]}(p)\leq C_p\|f\chi_{(1,q]}\|_F= C_p\|f\|_{F_q}.
$$
Consequently, we obtain
 $$
\|T\|_{\tilde{\mathfrak L}_{F}}\leq \Big(1+2\|\chi_{(q,\infty)}\|_{F}\,
C_{q}\Big) \|T\|_{\tilde{\mathfrak L}_{F_q}},
$$
and the equivalence of the norms of $\tilde{\mathfrak L}_F$ and 
$\tilde{\mathfrak L}_{F_q}$
is proven, since by assumption  $F$ contains the constant unit function and thus
contains $\chi_{(q,\infty)}$. \\
\indent
All these preliminary considerations show that the Boyd indices of $\clE$ coincide
with those of $\tilde{\mathfrak L}_{F_q}$, with $q>1$ arbitrary. So, let now $f$ be
 a measurable function on 
$(1,\infty)$ with finite $\tilde{\mathfrak L}_{F_q}$-norm. Then, for $a>1$, we have
$$
\|D_af\|_{\tilde{\mathfrak L}_{F_q}}=\big\|\big[p\mapsto\|D_af\|_p\big]\big\|_{F_q}
=\big\|\big[p\mapsto a^{-1/p}\|f\|_p\big]\chi_{(1,q]}\big\|_{F}\geq a^{-1/q}
\|f\|_{\tilde{\mathfrak L}_{F_q}}.
$$
This entails  $\mathfrak q_\clE\leq q$, which gives the result since $q$ can be chosen
arbitrarily  close to $1$.\hfill
\end{proof}

\noindent
{In the context of
extrapolated Lorentz spaces} contained in $\clN$,  we can relate the property of being close
to $\clL^1$ to the values of the Boyd indices.
\begin{cor}
Let  $\psi\in\Omega_b$ be such that $\clM_\psi$
is an extrapolated operator space. Then $\clM_\psi$ is close to $\clL^1$ (and thus
$\clM_\psi\in\mathcal X_1$) if and only if $\mathfrak p_{\clM_\psi}=\mathfrak q_{\clM_\psi}=1$. 
In this case, $\clM_\psi$ admits nontrivial singular traces supported at infinity, since then
every dilation invariant state is $\psi$-compatible.
\end{cor}
\noindent
\begin{proof}
This a combination of Corollary \ref{ew} and of Proposition \ref{BI}.
\hfill\end{proof}

\subsection{The lower and upper approximations}
\noindent
Our next task is to construct  lower and  upper approximations of a Lorentz space by
extrapolated spaces.
\begin{defn}
\label{UD}
Let $\psi\in\Omega$ be such that $\clM_\psi$ is close to $\clL^1$.
We then let
  \begin{align*}
F^\psi  &:=\Big\{f:(1,\infty)\to \mathbb C\,,\;\mbox{measurable} \;:\;\|f\|_{F^\psi}:=
{\rm ess}\sup_{p>1}
\frac{|f(p)|}{ \|\psi'\|_p}<\infty\Big\},\\
  F_\psi&:= \Big\{f:(1,\infty)\to \mathbb C\,,\;\mbox{measurable} \;:\;\|f\|_{F_\psi}:=
  {\rm ess}\sup_{p>1}
\frac{|f(p)|}{ \psi(e^{(p-1)^{-1}})}<\infty\Big\}.
  \end{align*}
Accordingly,
  we let  $\mathfrak L^\psi\equiv\mathfrak L_{F^\psi}$ and
    $\mathfrak L_\psi\equiv\mathfrak L_{F_\psi}$  be the associated 
  extrapolated operator spaces. Recall that their norms are given by
\begin{equation}
\label{upper}
\|T\|_{\mathfrak L^\psi}:=\|\eta_T\|_{F^\psi}=
{\rm ess}\sup_{p>1}\frac{\|T\|_p}{\|\psi'\|_p},
\qquad
\|T\|_{\mathfrak L_\psi}:=\|\eta_T\|_{F_\psi}=
{\rm ess}\sup_{p>1}\frac{\|T\|_p}{\psi(e^{(p-1)^{-1}})}.
\end{equation}
\end{defn}

\noindent
We also use the notations $\mathfrak L^\psi(\mathbb R_+)$ and
    $\mathfrak L_\psi(\mathbb R_+)$ to denote the extrapolated spaces
    described above for the commutative von Neumann algebra  $L^\infty(\mathbb R_+)$ 
 with trace given by the Lebesgue integral.

\begin{rem}
\label{rem2}
Observe that the fundamental functions of $\mathfrak L^\psi$   and $\mathfrak L_\psi$ 
are respectively given by:
\begin{align}
\label{CF}
\varphi^\psi\equiv\varphi_{ \mathfrak L^\psi}(t)=\sup_{p>1}\frac{t^{1/p}}{\|\psi'\|_p}\,,\qquad
\varphi_\psi\equiv\varphi_{ \mathfrak L_\psi}(t)=\sup_{p>1}\frac{t^{1/p}}{\psi(e^{(p-1)^{-1}})}.
\end{align}
Since the evaluation maps are continuous on $F_\psi$ and $F^\psi$:
$$
|f(p)|\leq \psi(e^{(p-1)^{-1}})\,\|f\|_{F_\psi}\,,\qquad |f(p)|\leq \|\psi'\|_p\,\|f\|_{F^\psi}.
$$
we  deduce from Lemma \ref{basic-prop},  that $\mathfrak L^\psi$
and  $\mathfrak L_\psi$
embed continuously in $\clL^p$ for all $p>1$. 
We also observe that $\mathfrak L_\psi\subset\clN$ for all $\psi\in\Omega$. Indeed, suppose
that there exists an element $T\in\mathfrak L_\psi$ such that $T\notin\clN$. Then $\|T\|_p\to\infty$,
$p\to\infty$, while $\|T\|_p\leq\|T\|_{\mathfrak L_\psi}\psi(e^{(p-1)^{-1}})\to\|T\|_{\mathfrak L_\psi}
 \psi(1)$, $p\to\infty$, a contradiction. In particular, we get for $T\in\mathfrak L_\psi$:
 $$
 \|T\|\leq\psi(1)\,\|T\|_{\mathfrak L_\psi}.
 $$
 When $\psi\in\Omega_b$, by the same reasoning we have 
further that
 $\mathfrak L^\psi\subset\clN$, as then $\psi'\in L^\infty(\mathbb R_+)$.
Last, as the constant unit function belongs to both $F_\psi$ and $F^\psi$, we see that 
$\mathfrak L^\psi$ and  $\mathfrak L_\psi$ belong to the class $\mathcal X_1$.
Hence, Proposition \ref{BI} applies, and shows
that  $\mathfrak L_\psi$  has trivial Boyd indices. This is also true for $\mathfrak L^\psi$ 
 when $\psi\in\Omega_b$.
\end{rem}

\noindent
We make an important observation. Since $\psi(0)=0$, we get from the
 H\"older inequality: 
\begin{align*}
\psi(e^r)=\int_0^{e^r}\psi'(t)\,dt\leq \|\chi_{[0,e^r]}\|_{1+r}\,\|\psi'\|_{1+1/r}=e^{r/(1+r)}  \|\psi'\|_{1+1/r}
\leq e  \|\psi'\|_{1+1/r}\,,\quad r>0.
\end{align*}
Specifying this to the case $r=(p-1)^{-1}$, $p>1$,  leads us to the following  embedding
result:

 \begin{lemma}
 Let $\psi\in \Omega$ be such that $\clM_\psi$ is close to $\clL^{1}$. 
 Then $F_\psi\subset F^\psi$ with $\|.\|_{F^\psi}
 \leq e\|.\|_{F_\psi}$. Consequently, $\mathfrak L_\psi\subset \mathfrak L^\psi$ with 
 $\|.\|_{\mathfrak L^\psi}\leq e\|.\|_{\mathfrak L_\psi}$.
\end{lemma}
\noindent
Our  approximation result for a Lorentz space close to $\clL^1$ is that one always has
\begin{equation}
\label{double}
\mathfrak L_\psi\subset \clM_\psi\subset\mathfrak L^\psi,
\end{equation}
with continuous inclusions:
\begin{prop}
\label{approx}
Let $\psi\in\Omega$ be such that $\clM_\psi$ is  close to 
$\clL^1$. Then $ \mathfrak L_\psi\subset\clM_\psi\subset \mathfrak L^\psi$
with 
$$
\|.\|_{ \mathfrak L^\psi}\leq \|.\|_{\clM_\psi}\leq \max\Big\{e,\frac{\psi(1)}{\psi'(1)}\Big\}\,\|.\|_{\mathfrak L_\psi}.
$$ 
\end{prop}
\noindent
\begin{proof}
In the proof of Lemma \ref{eq-close}, we have obtained for $T\in\clM_\psi$ and
for $p>1$ the inequality:
$$
\|T\|_{p}\leq\|T\|_{\clM_\psi} \,\|\psi'\|_{p}.
$$
This implies the first inequality.
Take next $T\in\mathfrak L_\psi$. 
As noticed before, we have $\mathfrak L_\psi\subset \clN$ with 
$ \|T\|\leq\psi(1)\,\|T\|_{\mathfrak L_\psi}$.
Thus for $t\in(0,1]$, we can use the estimates
$$
\int_0^t\mu(T,s) \,ds\leq \|T\|\,t\leq\psi(1)\,\|T\|_{\mathfrak L_\psi}\,t
\quad \mbox{and}\quad \psi(t)=\int_0^t\psi'(s)\,ds\geq \psi'(1)\,t,
$$
to get
$$
\frac{1}{\psi(t)}\int_0^t\mu(T,s) \,ds\leq \frac{\psi(1)}{\psi'(1)}\,\|T\|_{\mathfrak L_\psi},\quad t\in(0,1].
$$
For $t>1$, the H\"older inequality gives
for any $\eps\in(0,1)$:
\begin{align*}
\int_0^t\mu(T,s) \,ds\leq\Big(\int_0^t\,\mu\big(|T|^{1+\eps},s\big)\,ds
\Big)^{1/(1+\eps)}\Big(\int_0^t\,ds\Big)^{\eps/(1+\eps)}
&\leq\|T\|_{1+\eps}\,t^{\eps/(1+\eps)}\\ &\leq\,\|T\|_{\mathfrak L_\psi}\,t^{\eps/(1+\eps)}\,
\psi(e^{\eps^{-1}}).
\end{align*}
Choosing $\eps=1/\log(t)$, one gets
$$
\frac{1}{\psi(t)}\int_0^t\mu(T,s) \,ds\leq e \|T\|_{\mathfrak L_\psi},\quad t>1,
$$
proving the second inequality. \hfill
\end{proof}

\noindent
Before going any further, let us discuss two (elementary)  applications of Proposition \ref{approx}. 
The first one is in the context of tensor products of Lorentz spaces.  Let
 $\psi_j\in\Omega_b$, $j=1,2$, be such  that $\clM_{\psi_j}$
 is close to $\clL^1$. Then from the equality $\|A_1\otimes A_2\|_p=\|A_1\|_p\|A_2\|_p$, for all
 $A_j\in\clL^p(\clN_j,\tau_j)$ and $p\geq 1$, we easily deduce that
 $$
 \mathfrak L_{\psi_1}(\clN_1,\tau_1)\otimes\mathfrak L_{\psi_j}(\clN_2,\tau_2)
 \subset\mathfrak L_{\psi_1\psi_2}(\clN_1\otimes\clN_1,\tau_1\otimes\tau_2),
 $$
 continuously.  Using this result one can improve
 the embedding \eqref{ellele}. Indeed, by Proposition \ref{approx}, we get
 \begin{prop} 
Let $\psi_j\in\Omega_b$, $j=1,2$, be such $\mathfrak L_{\psi_j}=\clM_{\psi_j}$.
Then, we have a continuous embedding:
 $$
 \mathcal M_{\psi_1}(\clN_1,\tau_1)\otimes\mathcal M_{\psi_j}(\clN_2,\tau_2)
 \subset\mathcal M_{\psi_1\psi_2}(\clN_1\otimes\clN_1,\tau_1\otimes\tau_2).
 $$
 \end{prop}

Our second application is related
with  operators of the form $f(X)g(-i\nabla)$ (see \cite[Chapter 4]{Simon} for a detailed account
of this theory). For $\psi\in\Omega$, define 
$\clM_\psi(\R^n)$, $\mathfrak L_\psi(\R^n)$ and $\mathfrak L^\psi(\R^n)$ to be the Lorentz and
the associated extrapolated spaces constructed out of the
commutative von Neumann algebra $L^\infty(\mathbb R^{2n})$ with Lebesgue integral
 and define $\clM_\psi\big(L^2(\R^n)\big)$, $\mathfrak L_\psi\big(L^2(\R^n)\big)$ and 
 $\mathfrak L^\psi\big(L^2(\R^n)\big)$ to be their analogues  constructed out of the
von Neumann algebra $\mathcal B\big(L^2(\R^n)\big)$ with  operator trace.
For Borel functions $f,g$ on $\R^n$, let $f(X)$ be the operator of point-wise
multiplication by $f$,  that is $(f(X)g)(t)=f(t)g(t), t\in \mathbb R^{2n}$ and let $g(-i\nabla)$ be the operator of convolution  by the Fourier transform of $g$, 
that is $g(-i\nabla)=\mathcal{F}g(X)\mathcal{F}^{-1}$, where $\mathcal{F}$ is the Fourier transform.
 Let
 $\psi_1,\psi_2\in\Omega_b$.  Using the standard equalities $A^{1/2}BA^{1/2}=|B^{1/2}A^{1/2}|^2$ for positive bounded operators, we immediately infer that 
$\|A^{1/2}BA^{1/2}\|_p=\|B^{1/2}A^{1/2}\|_{2p}^2$, and obtain
 for $0\leq f,g\in L^p(\mathbb R^{2n})$, $p\geq 1$:
 \begin{align*}
& \|f^{1/2}(X)g(-i\nabla) f^{1/2}(X)\|_{\mathfrak L_{\psi_1\psi_2}}=
 {\rm ess}\sup_{p>1}\frac{\|f^{1/2}(X)g(-i\nabla) f^{1/2}(X)\|_p}{\psi_1(e^{(p-1)^{-1}})\psi_2(e^{(p-1)^{-1}})}\\
& \qquad\qquad\qquad \qquad\qquad\qquad \quad\;\;\,=
 {\rm ess}\sup_{p>1}\frac{\|g^{1/2}(-i\nabla) f^{1/2}(X)\|_{2p}^2}{\psi_1(e^{(p-1)^{-1}})\psi_2(e^{(p-1)^{-1}})}\\
 &
 \leq (2\pi)^{-\frac{n}p}{\rm ess}\sup_{p>1}\frac{\|g^{1/2}\|_{2p}^2\| f^{1/2}\|_{2p}^2}{\psi_1(e^{(p-1)^{-1}})\psi_2(e^{(p-1)^{-1}})}
 =(2\pi)^{-\frac{n}p}{\rm ess}\sup_{p>1}\frac{\|g\|_{p}\| f\|_{p}}{\psi_1(e^{(p-1)^{-1}})\psi_2(e^{(p-1)^{-1}})}\\
 &\leq (2\pi)^{-\frac{n}p} \|f\|_{\mathfrak L_{\psi_1}} \|g\|_{\mathfrak L_{\psi_2}},
 \end{align*}
where the first inequality follows from \cite[Theorem 4.1]{Simon}. Hence, Proposition \ref{approx}
yields:
\begin{prop} For $j=1,2$, let $\psi_j\in\Omega_b$ be such that
$\mathfrak L_{\psi_j}=\clM_{\psi_j}$.
Then, for  $0\leq f\in
 \clM_{\psi_1}(\R^n)$ and $0\leq g\in \clM_{\psi_2}(\R^n)$ we
 have $f^{1/2}(X)g(-i\nabla) f^{1/2}(X)\in \clM_\psi\big(L^2(\R^n)\big)$,
with
$$
 \|f^{1/2}(X)g(-i\nabla) f^{1/2}(X)\|_{\clM_\psi}\leq C \|f\|_{\clM_{\psi_1}} \|g\|_{\clM_{\psi_2}}.
 $$
 \end{prop}

\noindent 
Our next task  is to   understand when $\clM_\psi$ coincides with the 
extrapolated space $\mathfrak L_\psi$.  As seen before this is a very natural question
and, moreover, 
 as the main results of this article (Theorem \ref{zeta} and 
 Theorem \ref{psi2}) are valid on the
 positive cone  of $\mathfrak L_\psi$, only. It is also important to understand when
$\clM_\psi=\mathfrak L^\psi$. As shown next,  this equality of spaces
 is a necessary condition to have  $\clM_\psi=\mathfrak L_\psi$.
\begin{prop}
\label{eq-conds}
 Let $\clM_\psi$  be a Lorentz space close to $\clL^1$. 
 Then, the  three following  conditions are equivalent:\\
{\rm (i)} $\clM_\psi\in\mathcal X_1$,\\
{\rm (ii)} $\clM_\psi= \mathfrak L^\psi$ with equivalent norms,\\
{\rm (iii)} We have
 $$
 \psi(t)\asymp \Big(\sup_{p>1}\frac{t^{1/p-1}}{\|\psi'\|_p}\Big)^{-1}.
 $$
 \end{prop}
 \noindent
\begin{proof}
Assume that (ii) holds true. Then, the fundamental functions of  $\clM_\psi$ and 
$ \mathfrak L^\psi$ are equivalent, which implies (iii) using the first equality in
\eqref{CF} and the fact that the 
fundamental function of $\clM_\psi$ is $t/\psi(t)$. Next, assume that (iii) holds true. 
Note that by Proposition \ref{approx}
$\|.\|_{ \mathfrak L^\psi}\leq \|.\|_{\clM_\psi}$, so that the fundamental function  of 
$\mathfrak L^\psi$ is always smaller than that of $\clM_\psi$. Since (iii) also contains the 
converse inequality, (iii) implies that the fundamental functions of 
$\mathfrak L^\psi$ and $\clM_\psi$
are equivalent. As $\clM_\psi$ is the largest symmetric  operator space with a given 
fundamental 
function (see \cite[Chapter II, Theorem 5.7]{KPS}),
 we deduce that $\clM_\psi\supset \mathfrak L^\psi$, which implies (ii) by Proposition 
\ref{approx} again. Thus (ii) $\Leftrightarrow$ (iii). 

That (i) $\Leftarrow $ (ii) is trivial. Thus, we just need  prove that (i) $\Rightarrow $ (ii). 
So assume that the exists a Banach lattice $F$ on the interval $(1,\infty)$ such that 
$\clM_\psi=\mathfrak L_F$. Take $f\in F^\psi$. Note that $|f(p)|\leq \|f\|_{F^\psi}\|\psi'\|_p$, 
$p>1$. Since moreover
 $\psi'\in\clM_\psi(\mathbb R_+)= \mathfrak L_F(\mathbb R_+)$, we deduce that 
 $$
 \big[p\in(1,\infty)\mapsto \|\psi'\|_p\big]\in F.
 $$
 As $F$ is a Banach lattice, we deduce that $f\in F$. Hence $F^\psi\subset F$ and consequently 
 $\mathfrak L^\psi\subset \mathfrak L_F=\clM_\psi$. As the reversed inclusion
  $\mathfrak L^\psi\supset\clM_\psi$ is always true by Proposition \ref{approx}, 
  we get 
 (i) $\Leftrightarrow$ (ii). 
\hfill
\end{proof}

\noindent
Proposition \ref{eq-conds} tells us that we cannot have 
$\mathfrak L_\psi= \clM_\psi\subsetneqq\mathfrak L^\psi$ (because $\mathfrak L_\psi= \clM_\psi$
implies that $\clM_\psi\in\mathcal X_1$ and $ \clM_\psi\subsetneqq\mathfrak L^\psi$ implies that
 $\clM_\psi\notin\mathcal X_1$),
but  the three other possibilities $\mathfrak L_\psi= \clM_\psi=\mathfrak L^\psi$, 
$\mathfrak L_\psi\subsetneqq\clM_\psi=\mathfrak L^\psi$ and $\mathfrak L_\psi\subsetneqq\clM_\psi
\subsetneqq\mathfrak L^\psi$, are a priori possible. 
In particular, $\mathfrak L_\psi= \clM_\psi$ is equivalent to 
$\mathfrak L_\psi= \clM_\psi=\mathfrak L^\psi$.
Observe also that $\mathfrak L^\psi\subset \clN$ if and only if $\psi\in\Omega_b$. 
Indeed, if $\mathfrak L^\psi\subset \clN$, then $\clM_\psi\subset \clN$ and thus 
$\psi\in\Omega_b$. Conversely, if $\psi\in\Omega_b$, then $\psi'\in L^\infty(\mathbb R_+)$
which as observed in Remark \ref{rem2}, entails that $\mathfrak L^\psi\subset \clN$.
Hence, if $\psi\in\Omega\setminus\Omega_b$, we get 
$\mathfrak L_\psi\subsetneqq \clM_\psi$. 
 Last we observe that $\psi'$ is always
 in the unit ball of $\mathfrak L^\psi$ but it need not to belong to $\mathfrak L_\psi$. In fact,
 $\psi'\in\mathfrak L_\psi$ is equivalent to $\mathfrak L_\psi
= \clM_\psi=\mathfrak L^\psi$:
\begin{prop} 
\label{conds-2}
 Let $\clM_\psi$  be a  Lorentz space close to $\clL^1$.
Then, the three conditions below are equivalent:\\
{\rm (i)} $\mathfrak L_\psi= \clM_\psi$, with equivalent norms,\\
{\rm (ii)} $\psi'\in\mathfrak L_\psi(\mathbb R_+)$,\\
{\rm (iii)} There exists $C>0$, such that for all $p>1$ we have $\|\psi'\|_p\leq C \,\psi(e^{(p-1)^{-1}})$.\\
Moreover, any of these conditions  implies that $ \clM_\psi=\mathfrak L^\psi$ and also
that $\psi\in\Omega_b$.
\end{prop}
\noindent\begin{proof}
That (i) $\Rightarrow$ (ii), follows from the fact that $\psi'\in\clM_\psi(\mathbb R_+)$. Next, by definition of the norm
of $\mathfrak L_\psi$ (see the second equality in \eqref{upper}), we have
$$
\|\psi'\|_p\leq\|\psi'\|_{\mathfrak L_\psi}\,\psi\big(e^{(p-1)^{-1}}\big), \quad\forall p>1,
$$
hence (ii) $\Rightarrow$ (iii).
For the last part, note that in  the proof of Proposition \ref{approx}, we have obtained 
$\|T\|_p\leq\|T\|_{\clM_\psi}\|\psi'\|_p$, $\forall T\in\clM_\psi$, $\forall p>1$. Thus, if (iii) holds
we get for all $T\in\clM_\psi$:
$$
\frac{\|T\|_p}{\psi\big(e^{(p-1)^{-1}}\big)}\leq C\|T\|_{\clM_\psi},
$$
and thus $\|T\|_{\mathfrak L_\psi}\leq C\|T\|_{\clM_\psi}$, which entails that 
$\clM_\psi\subset \mathfrak L_\psi$. 
As the converse inclusion has been proven in  Proposition \ref{approx},
we deduce that (iii) $\Rightarrow$ (i), completing the proof.
\hfill\end{proof}

\noindent
From Proposition \ref{eq-conds}, we also deduce a  criterion for the equality
$\mathfrak L_\psi= \clM_\psi$ when $ \clM_\psi=\mathfrak L^\psi$.
\begin{prop}
Assume that $\clM_\psi\in\mathcal X_1$. Then, a sufficient condition to have
$\mathfrak L_\psi= \clM_\psi$, is that
\begin{align}
\label{v}
\sup_{p,q>1}\exp\Big\{-\frac{q-1}{q(p-1)}\Big\}\,\frac{\|\psi'\|_p}{\|\psi'\|_{q}}<\infty.
\end{align}
\end{prop}
\noindent
\begin{proof}
Under the assumption that $\clM_\psi\in\mathcal X_1$, we have by Proposition \ref{eq-conds} (iii):
$$
\frac{1}{\psi(t)}\leq C\, \sup_{p>1}\frac{t^{1/p-1}}{\|\psi'\|_p}.
$$
Hence
\begin{align*}
\|\psi'\|_{\mathfrak L_\psi}=\sup_{p>1}\frac{\|\psi'\|_p}{\psi\big(e^{(p-1)^{-1}}\big)}
\leq C\,\sup_{p>1}\sup_{q>1}\frac{\|\psi'\|_p}{\|\psi'\|_q}\Big(e^{(p-1)^{-1}}\Big)^{1/q-1}
=C\, \sup_{p,q>1}\,\frac{\|\psi'\|_p}{\|\psi'\|_{q}}\,e^{-\frac{q-1}{q(p-1)}}.
\end{align*}
Thus if \eqref{v} holds true, the inequality above  shows that $\psi'\in \mathfrak L_\psi(\mathbb R_+)$,
 so
 Proposition \ref{conds-2} gives  $\mathfrak L_\psi= \clM_\psi$ as needed.
\hfill\end{proof}

\noindent
In order to produce examples of element $\psi\in\Omega_b$ such that $\mathfrak L_\psi=\clM_\psi=
\mathfrak L^\psi$,
we make the following easy observation:
\begin{lemma}
\label{lem:trivial}
Let $\psi\in\Omega_b$ be such that  $\clM_\psi$  is close to $\clL^1$. Then $\mathfrak L_\psi=\clM_\psi=
\mathfrak L^\psi$ if 
and only if there exists $C>0$
such that for all $\eps>0$, we have
\begin{align}
\label{trivial}
\int_{e^{\eps^{-1}}}^\infty \psi'(t)^{1+\eps} \,dt\leq C\,\psi\big(e^{\eps^{-1}}\big).
\end{align}
\end{lemma}
\noindent
\begin{proof}
For $\eps>0$ small enough, $\|\psi'\|_{1+\eps}\geq 1$ and thus since $\psi'$ is decreasing, 
\begin{align*}
\|\psi'\|_{1+\eps}\leq \|\psi'\|_{1+\eps}^{1+\eps}\leq\psi'(0)^\eps\int_0^{e^{\eps^{-1}}}\psi'(t)\,dt+
\int_{e^{\eps^{-1}}}^\infty \psi'(t)^{1+\eps} \,dt=\psi'(0)^\eps\psi\big(e^{\eps^{-1}}\big)
+\int_{e^{\eps^{-1}}}^\infty \psi'(t)^{1+\eps} \,dt.
\end{align*}
Hence, Proposition \ref{conds-2} (iii) shows that  \eqref{trivial}  implies that 
$\mathfrak L_\psi=\clM_\psi$. The converse is obvious by Proposition \ref{conds-2} (iii) 
again:
$$
\int_{e^{\eps^{-1}}}^\infty \psi'(t)^{1+\eps} \,dt\leq \|\psi'\|_{1+\eps}^{1+\eps}\leq C
\,\psi\big(e^{\eps^{-1}}\big)^{1+\eps}\leq C'\,\psi\big(e^{\eps^{-1}}\big),
$$
since $\psi(t)\leq \psi'(0) t$ and $\psi'(0)=\|\psi'\|_\infty<\infty$ as $\psi\in\Omega_b$.
\hfill\end{proof}

\noindent
We are now able to give an handy  criterion  for the  equality 
$\mathfrak L_\psi=\clM_\psi=\mathfrak L^\psi$:
\begin{prop}
\label{prop-eq}
Let $\psi\in\Omega_b$. Assume that  for all $\delta>0$,  the map $[t\mapsto t^{-\delta}\psi(t)]$ is decreasing 
 and  that there exists $\rho>0$ such that the map $[t\mapsto\psi\big(\exp(t^\rho)\big)]$
is still concave. Then $\mathfrak L_\psi=\clM_\psi=\mathfrak L^\psi$.
\end{prop}
\noindent
\begin{proof}
Let $\varphi(t):=\psi\big(e^{t^\rho}\big)$. By concavity of $\varphi$ we get 
$$
\rho\, t^\rho\, e^{t^\rho}\psi'\big(e^{t^\rho}\big)=t\,\varphi'(t)\leq \varphi(t)=\psi\big(e^{t^\rho}\big).
$$
Performing the change of variable $t\mapsto \log(t)^{1/\rho}$
  gives for
$t>1$:
$$
\psi'(t)\leq\frac{\psi(t)}{\rho\, t\log(t)}.
$$
This clearly implies that $\clM_\psi$ is close to $\clL^1$ 
as $\psi(t)\leq Ct^\delta$ for all $\delta>0$ by assumption, 
and moreover
\begin{align*}
\int_{e^{\eps^{-1}}}^\infty \psi'(t)^{1+\eps}\, dt&\leq \rho^{-1-\eps}
\int_{e^{\eps^{-1}}}^\infty \frac{\psi(t)^{1+\eps}}{t^{1+\eps}\log(t)^{1+\eps} }\,dt\\
&\leq \rho^{-1-\eps}\sup_{s\geq e^{\eps^{-1}}}\Big\{\frac{\psi(s)^{1+\eps}}{s^{\eps/2}\log(s)^{1+\eps} }\Big\}\;
\int_{e^{\eps^{-1}}}^\infty t^{-1-\eps/2 }\,dt
= \frac2{e\,\rho^{1+\eps}}\,\psi\big(e^{\eps^{-1}}\big)^{1+\eps}\,\eps^{\eps}.
\end{align*}
 In the last equality, we use the assumption that 
the map $[t\mapsto t^{-\delta}\psi(t)]$ is decreasing  for all $\delta>0$, which implies that
$$
\sup_{s\geq e^{\eps^{-1}}}\Big\{\frac{\psi(s)^{1+\eps}}{s^{\eps/2}\log(s)^{1+\eps} }\Big\}=
e^{-1/2}\eps^{1+\eps}\psi(e^{\eps^{-1}}),
$$
and the the final expression follows by performing the integral.
The conclusion is then a consequence of  Lemma \ref{lem:trivial}.
\hfill\end{proof}

\begin{rem}
The concavity of $[t\mapsto\psi\big(\exp(t^\rho)\big)]$ implies that 
$\lim_{t\to0} t\psi'(t)/\psi(t)=0$, 
which by Lemma \ref{eq-conds-1}  implies that condition \eqref{ex-sing-trace-easy} is satisfied.
\end{rem}

\begin{cor}
\label{obs}
Let $\psi \in\Omega_b$ be either 
the elements $\psi_{n,\beta}$ given in \eqref{ex1}, for $n\geq1$ and $\beta>0$ or
 the elements $\tilde\psi_{n,\beta}$ given in \eqref{ex2},  for $n>1$ and $\beta\in(0,1)$. 
 Then we have $\mathfrak L_\psi=\clM_\psi=\mathfrak L^\psi$.
 \end{cor}
 However, for $n=1$ and $\beta\in(0,1)$, 
there is no $\rho>0$ such that $[t\mapsto \tilde\psi_{1,\beta}\big(\exp(t^\rho)\big)]$ is concave. 
Moreover,  
$\mathfrak L_{\tilde\psi_{1,\beta}}\subsetneqq \clM_{\tilde\psi_{1,\beta}}$.
For example, take $\beta=1/2$. We then have
$\tilde\psi_{1,1/2}'(t)=(2t\sqrt{\log t})^{-1}e^{\sqrt{\log t}}$, so
 performing the change of variable $t\mapsto \log(t)^{1/2}$, we get:
 \begin{align*}
\int_1^\infty \tilde\psi_{1,1/2}'(t)^{1+\eps}\,dt&=\int_0^\infty e^{-\eps t^2} t^{-\eps}e^{(1+\eps)t}\,dt\\
&\geq \int_0^\infty e^{-\eps t^2} e^{t}\,dt
=\frac{\sqrt\pi}{2\sqrt\eps}\,e^{1/4\eps}\big(1+{\rm Erf}\big(
\frac{1}{2\sqrt\eps}\big)\big)\sim_{\eps\to0}\sqrt\frac{\pi}{\eps}\,e^{1/4\eps},
 \end{align*}
 while $\psi_{1,1/2}'(e^{1/\eps})=e^{1/\sqrt\eps}$, and we conclude using  Lemma \ref{lem:trivial}.
 However, we don't know in these cases if $ \clM_{\tilde\psi_{1,\beta}}\subsetneqq
  \mathfrak L^{\tilde\psi_{1,\beta}}$ or $ \clM_{\tilde\psi_{1,\beta}}= \mathfrak L^{\tilde\psi_{1,\beta}}$.

\noindent
Last, we  give stability properties for the class
 of elements $\psi$ in $\Omega_b$ such that $\mathfrak L_\psi=\clM_\psi=\mathfrak L^\psi$.
\begin{prop}
\label{prop-eq2}
{\rm (i)} Let $\psi\in\Omega_b$ be such that $\mathfrak L_\psi=\clM_\psi$ and  such that for all 
$\delta>0$,  the map $[t\mapsto t^{-\delta}\psi(t)]$ is decreasing. Then for any  $\alpha>0$, we have
 $\mathfrak L_{\psi^\alpha}=\clM_{\psi^\alpha}$.\\
{\rm (ii)} Let  $\psi_1,\psi_2\in\Omega_b$ be such that $\mathfrak L_{\psi_i}=\clM_{\psi_i}$
and such that for all 
$\delta>0$,  the map $[t\mapsto t^{-\delta}\psi_i(t)]$ is decreasing. Then we have
$\mathfrak L_{\psi_1\psi_2}=\clM_{\psi_1\psi_2}$.\\
{\rm (iii)} Let  $\psi_1,\psi_2\in\Omega_b$  be such  that for all 
$\delta>0$,  the map $[t\mapsto t^{-\delta}\psi_2(t)]$ is decreasing and such that 
there exists $\rho>0$ such that the map
$[t\mapsto\psi_2\big(\exp(t^\rho)\big)]$ is still concave, then
$\mathfrak L_{\psi_1\circ\psi_2}=\clM_{\psi_1\circ\psi_2}$.
\end{prop}
\noindent
\begin{proof}
To prove (i), fix $\delta>0$, and observe that by concavity of $\psi$, we have
$$
 \psi'(t)^{\eps\delta}\,\psi(t)^{(\alpha-1)(1+\eps)}\leq \psi(t)^{\eps\delta+(\alpha-1)(1+\eps)} 
 t^{-\eps\delta}\,,\qquad\forall\eps>0.
$$
Thus, we obtain:
\begin{align*}
\int_{e^{\eps^{-1}}}^\infty \psi'(t)^{1+\eps}\,\psi(t)^{(\alpha-1)(1+\eps)}\,dt&\leq
\sup_{s\geq e^{\eps^{-1}}}\big\{ s^{-\eps\delta}\,\psi(s)^{\eps\delta+(\alpha-1)(1+\eps)}\big\}
\,\|\psi'\|_{1+\eps(1-\delta)}^{1+\eps(1-\delta)}\\
&\quad=e^{-\eps\delta}\,
\psi\big(e^{\eps^{-1}}\big)^{\eps\delta+(\alpha-1)(1+\eps)}\,
\|\psi'\|_{1+\eps(1-\delta)}^{1+\eps(1-\delta)},
\end{align*}
since the map $[t\mapsto \psi(t)^{\eps\delta+(\alpha-1)(1+\eps)} 
 t^{-\eps\delta}]$ is decreasing. Letting $\delta$ going to zero, we finally deduce:
\begin{align*}
 \int_{e^{\eps^{-1}}}^\infty \psi'(t)^{1+\eps}\,\psi(t)^{(\alpha-1)(1+\eps)}\,dt&\leq
 \psi\big(e^{\eps^{-1}}\big)^{(\alpha-1)(1+\eps)}\,
\|\psi'\|_{1+\eps}^{1+\eps}
\leq  C\,\psi\big(e^{\eps^{-1}}\big)^{(\alpha-1)}\,
\|\psi'\|_{1+\eps}^{1+\eps},
\end{align*}
where we used $\psi(t)\leq \psi'(0) t$ and $\psi'(0)<\infty$ as $\psi\in\Omega_b$, when 
$\alpha\geq1$.
One concludes using Proposition \ref{conds-2} (iii) and Lemma \ref{lem:trivial}.
The proof of  (ii) is similar. On starts from:
$$
\|(\psi_1\psi_2)'\|_{1+\eps}\leq \|\psi_1\psi_2'\|_{1+\eps}+\|\psi_1'\psi_2\|_{1+\eps},
$$
and one uses the fact that $\psi_1$ is increasing to get:
\begin{align*}
\int_0^{e^{\eps^{-1}}} \psi_1(t)^{1+\eps}\psi_2'(t)^{1+\eps}\,dt&\leq 
\psi_1\big(e^{\eps^{-1}}\big)^{1+\eps}
\int_0^{e^{\eps^{-1}}} \psi_2'(t)^{1+\eps}\,dt
\leq 
 C\psi_1\big(e^{\eps^{-1}}\big)
\psi_2\big(e^{\eps^{-1}}\big).
\end{align*}
Then one uses the concavity of $\psi_2$, to estimate with $\delta>0$ arbitrary:
$$
\psi_1(t)^{1+\eps}\psi_2'(t)^{\eps\delta}\leq  \psi_1(t)^{1+\eps}\psi_2(t)^{\eps\delta}t^{-\eps\delta},
$$
which from the same reasoning as in the proof of (i) entails that:
$$
\int_{e^{\eps^{-1}}}^\infty \psi_1(t)^{1+\eps}\psi_2'(t)^{1+\eps}\,dt\leq
\psi_1\big(e^{\eps^{-1}}\big)^{1+\eps}\,
\|\psi_2'\|_{1+\eps}^{1+\eps}\leq  C
 \psi_1\big(e^{\eps^{-1}}\big)\,\psi_2\big(e^{\eps^{-1}}\big),
$$
and one concludes with Proposition \ref{conds-2} (iii). The case of $\psi_1'\psi_2$
is identical, proving (ii). Last, from $\psi_1(t)\leq \psi_1'(0)t$ we observe that the conditions
of Proposition \ref{prop-eq} are satisfied for $\psi_1\circ\psi_2$ and (iii) follows.
This concludes the proof.
\hfill\end{proof}

\section{Dixmier traces and $\zeta$-functions}
\label{3}

In this section, we will formulate the first main result of this  article, which is Theorem \ref{zeta}.
It answers positively the first part of the question raised in the introduction, for all the Lorentz spaces 
$\clM_\psi$ which are associated with elements
$\psi\in\Omega_b$ satisfying   condition \eqref{our-cond} and which are such that 
$\mathfrak L_{\psi}=\clM_\psi$.
Note that  this result  connects
section \ref{1} and  section \ref{2}. 
The proof of Theorem \ref{zeta} relies on Proposition \ref{calc}  and Proposition \ref{lem2}.
The latter can be thought as 
    a variant of  the weak$^*$-Karamata Theorem, generalizing  \cite[Theorem 2.2]{CPS2}.
 Its proof essentially relies on the following Lemma:

 \begin{lemma}
 \label{lem1}
Let  $\beta:[0,\infty)\to[0,\infty)$ be a non-negative, increasing  and right continuous function 
and let $\vf:[0,\infty)\to[0,\infty)$ be a measurable function satisfying
 \begin{equation}
 \label{xw}
 \lim_{r\to\infty}\frac{\vf(r/n)}{\vf(r)}=n^{-k},\quad \mbox{for all}\; n\in\mathbb N\;\mbox{and
 for  some } \;k\geq0.
 \end{equation}
 Suppose also that
 $$
 \Big[r\mapsto\frac1{\vf(r)}\int_0^\infty e^{-t/r}\,d\beta(t)\Big]\in  L^\infty(\mathbb R^*_+),
 $$
and let $\omega$ be a  dilation invariant state of $ L^\infty(\mathbb R_+^*)$. Then,
for every $f\in C[0,1]$, we have
\begin{align}
\label{quation}
&\omega\Big(\Big[r\mapsto 
\frac1{\vf(r)}\int_0^\infty f\big(e^{-t/r}\big)\,e^{-t/r}\,d\beta(t)\Big]\Big)
=C(f)\;\omega\Big(\Big[r\mapsto
\frac1{\vf(r)}\int_0^\infty e^{-t/r}\,d\beta(t)\Big]\Big),
 \end{align}
 where
$$
 C(f):=
 \begin{cases}
 \frac 1{\Gamma(k)}\,
\int_0^1 f(t)\,|\log(t)|^{k-1}\,dt,\quad &\mbox{for} \quad k>0,\\
f(1),\quad &\mbox{for} \quad k=0.
\end{cases}
$$
  \end{lemma}
 \noindent
 \begin{proof}
 First of all, observe that since $f$ is bounded,
  both sides of the equality \eqref{quation} are well defined.
  Indeed, for the left hand side we have
  $$
  \sup_{r>0} \frac1{\vf(r)}\Big|\int_0^\infty f\big(e^{-t/r}\big)\,e^{-t/r}\,d\beta(t)\Big|\leq\|f\|_\infty
  \sup_{r>0}   \frac1{\vf(r)}\int_0^\infty e^{-t/r}\,d\beta(t),
  $$
  which is finite by assumption, and for the right hand side we clearly have $|C(f)|\leq\|f\|_\infty$
  in both cases $k>0$ and $k=0$.
  In particular,  both sides of the equality \eqref{quation} are continuous in $f\in C[0,1]$
  for the uniform topology.
  By density of polynomials in $C[0,1]$ and by linearity, it is sufficient to verify  \eqref{quation}
  for $f_n(t):=t^n$, $n\in\mathbb N_0$. For such a function, we   have $C(f_n)=1$ when $k=0$ and when
  $k>0$:
  \begin{align*}
   \frac 1{\Gamma(k)}\,
\int_0^1 f_n(t)\,|\log(t)|^{k-1}\,dt=
 \frac 1{\Gamma(k)} \int_0^1 t^n\,|\log(t)|^{k-1}\,dt
= \frac 1{\Gamma(k)}\int_0^\infty e^{-u(n+1)}\,
 u^{k-1}\,du= (n+1)^{-k}.
 \end{align*}
 Now, for the left hand side of \eqref{quation}, we find for $k\geq 0$:
 \begin{align*}
 \omega\Big(\Big[r\mapsto 
\frac1{\vf(r)}\int_0^\infty e^{-(n+1)t/r}\,d\beta(t)\Big]\Big)&= \omega\Big(\Big[r\mapsto
\frac{\vf(r/(n+1))}{\vf(r)}\frac1{\vf(r/(n+1))}\int_0^\infty e^{-(n+1)t/r}\,d\beta(t)\Big]\Big)\\
&=(n+1)^{-k} \;\omega\Big(\Big[r\mapsto
\frac1{\vf(r/(n+1))}\int_0^\infty e^{-(n+1)t/r}\,d\beta(t)\Big]\Big)\\
&=(n+1)^{-k}\; \omega\Big(\Big[r\mapsto
\frac1{\vf(r)}\int_0^\infty e^{-t/r}\,d\beta(t)\Big]\Big),
 \end{align*}
 where the second equality follows from assumption \eqref{xw} and the third one
 follows from the dilation-invariance of $\omega$. 
 \hfill\end{proof}

  \begin{prop}
 \label{lem2}
 Under the  assumptions on the functions $\vf$, $\beta$ and on the extended limit
 $\omega$ of Lemma \ref{lem1}, we have
 \begin{align}
 \label{uation}
\omega\Big(\Big[r\mapsto
\frac{\beta(r)}{\vf(r)}\Big]\Big)=\frac 1{\Gamma(1+k)}\,
\omega\Big(\Big[r\mapsto
\frac1{\vf(r)}\int_0^\infty e^{-t/r}\,d\beta(t)\Big]\Big).
 \end{align}
 \end{prop}
 \noindent
 \begin{proof}
 Consider $L$, the continuous linear functional given by the left hand side of \eqref{quation},
 namely:
 $$
 L:C[0,1]\to\mathbb C,\quad f\mapsto 
\omega\Big(\Big[r\mapsto
\frac1{\vf(r)}\int_0^\infty f\big(e^{-t/r}\big)\,e^{-t/r}\,d\beta(t)\Big]\Big).
$$
$L$ defines a Radon measure on $[0,1]$ thus it extends to piecewise continuous functions. 
Now, consider 
$$
g(x):=\frac{\chi_{[e^{-1},1]}(x)}x.
$$
We then obtain
\begin{align*}
L(g)=\omega\Big(\Big[r\mapsto
\frac1{\vf(r)}\int_0^\infty \chi_{[e^{-1},1]}\big(e^{-t/r}\big)\,d\beta(t)\Big]\Big)
=\omega\Big(\Big[r\mapsto
\frac1{\vf(r)}\int_0^r \,d\beta(t)\Big]\Big)=\omega\Big(\Big[r\mapsto
\frac{\beta(r)}{\vf(r)}\Big]\Big).
\end{align*}
Thus, all what remains to do is to show that the right hand side of \eqref{uation},
coincides with $L(g)$ too. To see this, we 
choose  two sequences $\{g^\pm_n\}_{n\in\mathbb N}$ in $C[0,1]$, 
which converge point-wise   to $g$ and such that there exist 
$a<b\in\mathbb R$ with
$a\leq g^-_n\leq g\leq g^+_n\leq b$,  $\forall n\in\mathbb N$.
By positivity, we have $L(g_n^-)\leq L(g)\leq L(g_n^+)$ and thus it suffices 
to prove that the numerical sequences $\{L(g_n^\pm)\}_{n\in\N}$ converge to 
the right hand side of \eqref{uation}. We first  deduce from Lemma \ref{lem1} that
$$
L(g^\pm_n)=C(g^\pm_n)\;
\omega\Big(\Big[r\mapsto
\frac1{\vf(r)}\int_0^\infty e^{-t/r}\,d\beta(t)\Big]\Big),
$$
with 
$$
 C(g^\pm_n)=
 \begin{cases}
 \frac 1{\Gamma(k)}\,
\int_0^1 g^\pm_n(t)\,|\log(t)|^{k-1}\,dt,\quad &\mbox{for} \quad k>0,\\
g^\pm_n(1),\quad &\mbox{for} \quad k=0.
\end{cases}
$$
Accordingly,  it suffices 
to prove that the numerical sequences $\{C(g_n^\pm)\}_{n\in\N}$ converge to
$\Gamma(1+k)^{-1}$. For $k=0$ this is obvious since  $C(g^\pm_n)=g^\pm_n(1)$ converges to 
$g(1)=1$. For $k>0$, 
we use  dominated convergence to get that $C(g^\pm_n)$ converges to $C(g)$.
Last, we compute
\begin{align*}
C(g)=\frac 1{\Gamma(k)}\int_0^1 g(t)|\log(t)|^{k-1}\,dt&=\frac 1{\Gamma(k)}
\int_{e^{-1}}^1 |\log(t)|^{k-1}\,
\frac{dt}t=\frac 1{\Gamma(k)}\int_0^1u^{k-1}du=\frac1{k\Gamma(k)}=\frac 1{\Gamma(1+k)},
\end{align*}
which  completes the proof.
 \hfill\end{proof}

 \begin{thm}
\label{psi}
Let $\psi\in\Omega$ satisfying  condition \eqref{our-cond}.
Let also $k_\psi\geq 0$ be the constant  given  in Definition \ref{kpsi} and
let finally   $\omega$ be  an exponentiation invariant state of $ L^\infty(\mathbb R_+^*)$.
Then, for all 
 $ T\in \mathfrak L_\psi^+$ (the subcone of $\clM_\psi^+$  given in Definition \ref{UD}),  we have
\begin{align}
\label{main}
\tau_{\psi,\omega}(T)=
\frac1{\Gamma(1+k_\psi)}\,\omega\Big(\Big[r\mapsto \frac{\tau(T^{1+1/\log(r)})}{\psi(r)}\Big]\Big).
\end{align} 
 \end{thm}
\noindent
\begin{proof}
That $\tau_{\psi,\omega}$ defines a Dixmier trace for  an exponentiation invariant  
state $\omega$
and $\psi$ satisfying  condition \eqref{our-cond} follows from Proposition \ref{lin}
and the fact that condition \eqref{our-cond} implies condition \eqref{cond-mod}.
Indeed, we have 
$$
\omega\Big(\Big[t\mapsto\frac{\psi(t^\alpha)}{\psi(t)}\Big]\Big)=
\lim_{t\to\infty}\frac{\psi(t^\alpha)}{\psi(t)}=A_\psi(\alpha), \quad \alpha>1,
$$
 thus from Lemma \ref{F}, we deduce 
$$
\lim_{\alpha\downarrow 1}\,\omega\Big(\Big[t\mapsto\frac{\psi(t^\alpha)}{\psi(t)}\Big]\Big)=
\lim_{\alpha\downarrow 1}\,
A_\psi(\alpha)=1.
$$
Fix $T\in  \mathfrak L_\psi^+$,  set $E_T$ for the spectral family of $T$ and define
$$
\beta(t):=\tau\big(T  \,E_T(e^{-t})\big).
$$
Note that $\beta$ is  finite (since $T$ is $\tau$-compact), real valued, increasing  and 
right continuous. Define also
$\vf:=\psi\circ\exp$.
Now, observe that for $n\in\mathbb N$, we get by Lemma \ref{F}:
$$
\lim_{r\to\infty} \frac{\vf(r/n)}{\vf(r)}=\lim_{r\to\infty} \frac{\psi(e^{r/n})}{\psi(e^r)}=
\lim_{t\to\infty} \frac{\psi(t^{1/n})}{\psi(t)}=A_\psi(1/n)=n^{-k_\psi}.
$$
Next, we observe that by spectral theory
\begin{align*}
\int_0^\infty e^{-t/r}\,d\beta(t)=\int_0^\infty e^{-t/r}\,d\tau\big(T  \,E_T(e^{-t})\big)=
\tau\Big(T\int_0^\infty e^{-t/r}\,d  E_T(e^{-t})\Big)=\tau\big(T^{1+1/r}\big).
\end{align*}
By assumption we have $T\in\mathfrak L_\psi^+$ and so 
$$
 \Big[r\mapsto\frac1{\vf(r)}\int_0^\infty e^{-t/r}\,d\beta(t)\Big]\in  L^\infty(\mathbb R^*_+),
 $$
 where $\omega\circ\log$ is the dilation
invariant state  of $L^\infty(\mathbb R_+^*)$ defined by
$$
\omega\circ\log(f)=\omega\big(f\circ\log).
$$
Next, by  Proposition \ref{calc} and Remark \ref{ouf}, we have:
$$
\tau_{\psi,\omega}(T)=\omega\Big(\Big[t\mapsto\frac{1}{\psi(t)}\,\tau\big( T\,\chi_{(1/t,1)}(T)\big)\Big]
\Big)
=\omega\Big(\Big[t\mapsto\frac{\beta\big(\log(t)\big)}{\psi(t)}\Big]\Big)
=\omega\circ\log\Big(\Big[t\mapsto\frac{\beta(t)}{\vf(t)}\Big]\Big).
$$
Hence, the previous computations and Proposition
 \ref{lem2} applied to $\beta$ and $\vf$ as defined above, and to the dilation invariant state
$\omega\circ\log$, give:
\begin{align*}
\tau_{\psi,\omega}(T)&=\omega\circ\log\Big(\Big[t\mapsto\frac{\beta(t)}{\vf(t)}\Big]\Big)=
\frac 1{\Gamma(1+k_\psi)}\,
\omega\circ\log\Big(\Big[r\mapsto
\frac1{\vf(r)}\int_0^\infty e^{-t/r}\,d\beta(t)\Big]\Big)\\&=
\frac 1{\Gamma(1+k_\psi)}\,
\omega\circ\log\Big(\Big[r\mapsto
\frac{\tau\big(T^{1+1/r}\big)}{\psi(e^r)}\Big]\Big)=
\frac1{\Gamma(1+k_\psi)}\,\omega\Big(\Big[r\mapsto \frac{\tau(T^{1+1/\log(r)})}{\psi(r)}\Big]\Big).
\end{align*}
This completes the proof.\hfill
\end{proof}
 
 \noindent
The next result is  a generalization of \cite[Proposition 3.6]{CPS2}:
 \begin{thm}
 \label{zeta}
 Let $\psi\in\Omega$ satisfying  condition \eqref{our-cond}.
Let also   $\omega$ be an exponentiation invariant state of $ L^\infty(\mathbb R_+^*)$, 
 $ T\in \mathfrak L_\psi^+$ and $B\in\clN$. Then,  with the notations
 of Theorem \ref{psi} we have:
\begin{align}
\label{ty}
\tau_{\psi,\omega}(BT)=\frac1{\Gamma(1+k_\psi)}\,\omega\Big(\Big[r\mapsto\frac{\tau(BT^{1+1/\log(r)})}{\psi(r)}\Big]\Big).
\end{align}
 \end{thm}
 \noindent
 \begin{proof}
 First observe that both sides  of the equality \eqref{ty} are linear in $B$. We can therefore assume
 that $B\geq 0$. Next, assume for a moment that the result holds when $B\geq {\rm Id}_\clN$.
To simplify the notations,  set $R_T(B)$ for the right hand-side of \eqref{ty}. So for arbitrary
$B\geq0$, we deduce by Theorem \ref{psi} and by our (momentary) extra assumption:
\begin{align*}
\tau_{\psi,\omega}(BT)=\tau_{\psi,\omega}\big((B+{\rm Id}_\clN)T\big)-\tau_{\psi,\omega}(T)
&=R_T(B+{\rm Id}_\clN)-\tau_{\psi,\omega}(T)
\\&=R_T(B)+R_T({\rm Id}_\clN)-\tau_{\psi,\omega}(T)=
R_T(B).
\end{align*}
 Hence, we may assume without loss of generality that $B\geq  {\rm Id}_\clN$.
 Then, by traciality and Theorem \ref{psi},
  it suffices to prove that 
 \begin{equation}
 \label{eqB}
\omega\circ\log\Big(\Big[r\mapsto \frac{\tau\big((B^{1/2}TB^{1/2})^{1+1/r}
 \big)}{\psi(e^r)}\Big]\Big)=
 \omega\circ\log\Big(\Big[r\mapsto \frac{\tau\big(B^{1/2}T^{1+1/r}B^{1/2}
 \big)}{\psi(e^r)}\Big]\Big).
 \end{equation}
 But when $B\geq  {\rm Id}_\clN$, since $\mathfrak L_\psi\subset \clN$ (see remark \ref{rem2}), 
 we have by \cite[Lemma 3.3]{CPS2} the operator inequality:
 $$
 0\leq (B^{1/2}TB^{1/2})^{1+1/r}-B^{1/2}T^{1+1/r}B^{1/2}\leq\big(\|B\|^{1/r}-1\big) B^{1/2}T^{1+1/r}B^{1/2},\qquad\forall r>1,
 $$
 which gives the result as $\|B\|^{1/r}-1\to0$, $r\to\infty$.
 \hfill
 \end{proof}

\noindent
We conclude this section by applying our results to the elements of  $\Omega_b$ given in Example
\ref{ex}. First, for $\psi(t)=\log(1+t^{1/\beta})^\beta$, $\beta>0$,  condition \eqref{our-cond}
is satisfied with $k_\psi=\beta$ and   we have $\mathfrak L_\psi=\clM_\psi$ by Corollary \ref{obs}.
Therefore Theorem \ref{zeta} shows that for any $T\in\clM_\psi^+$, $B\in\clN$ and 
 exponentiation invariant state $\omega$, we have:
$$
\tau_{\psi,\omega}(BT)=\frac{\beta^\beta}{\Gamma(1+\beta)} \;\omega\circ\log\Big(\Big[r\mapsto
\frac{\tau(BT^{1+\frac1r})}{r^\beta}\Big]\Big).
$$

In particular, we observe that for the Lorentz space associated to the function
$\psi(t)=\log(1+t^{1/n})^n$ with $n\in\mathbb N$, 
when meromorphic, the $\zeta$-function $\zeta(T,z)=\tau(T^z)$ of $T\in\clM_\psi^+$  has a  pole of order $n$ 
at $z=1$. This generalizes the well known fact that, when meromorphic, the $\zeta$-function of a
positive operator in $\clM^{1,\infty}$ has a simple pole at $z=1$.

Next, take 
  $\psi(t)=\log(1+\log(1+t))$. Recall that we have then 
   $k_\psi=0$ and $\mathfrak L_\psi=\clM_\psi$ by Corollary \ref{obs}.
Hence,  for any $T\in\clM_\psi^+$, $B\in\clN$ and 
 an exponentiation invariant state $\omega$, we have
$$
\tau_{\psi,\omega}(BT)=\omega\circ\log\Big(\Big[r\mapsto
\frac{\tau(BT^{1+\frac1r})}{\log(r)}\Big]\Big).
$$

\section{Dixmier traces and heat kernels}
\label{4}

In this section, we answer positively the second part of the question raised in the introduction.
This is done in Theorem \ref{psi2}. As for Theorem \ref{zeta}, the answer is limited to
the Lorentz spaces 
$\clM_\psi$ which are associated with  elements
$\psi\in\Omega_b$ satisfying   condition \eqref{our-cond} and which are such that 
$\mathfrak L_{\psi}=\clM_\psi$.  We start by introducing the subset $\clC_\psi$ of $\tau$-measurable
operators in $\clN$ which are such that the right hand side of the equality \eqref{triple} is well
defined. We then show in Proposition \ref{C-M} and Proposition \ref{cont} that we 
have $\mathfrak L_\psi\subset \clC_\psi\subset \clM_\psi$ for all $\psi\in\Omega$ and
 with continuous inclusions. Moreover, we show in Corollary \ref{equal}, that
 condition \eqref{SSZ-cond} is sufficient to ensure that $\clC_\psi=\clM_\psi$.
This is why  we strongly  believe that
Theorem \ref{psi2} should hold with  condition \eqref{SSZ-cond} instead
of the stronger condition   \eqref{our-cond} that we use. 

\quad

 In what follows, for $T\in\clN^+$, we define  the operator $e^{-T^{-1}}$ to be zero on
${\rm ker}( T)$  and on ${\rm ker}( T)^\perp$ in the usual way by 
the functional calculus. 
\begin{defn}
\label{C}
For  $\psi\in\Omega$, we let $\clC_\psi(\clN,\tau)\equiv\clC_\psi$ to be the subset of 
 $\clL_0$ such that 
\begin{equation}
\label{CHK}
\sup_{\lambda>0}
\frac1{\psi(\lambda)}\int_0^\lambda \tau\big(e^{-u^{-1}|T|^{-1}}\big) \,\frac{du}{u^2}<\infty.
\end{equation}
\end{defn}

\noindent
Of course, we can write 
\begin{equation*}
\int_0^\lambda \tau\big(e^{-u^{-1}|T|^{-1}}\big) \,\frac{du}{u^2}=
 \tau\big(|T|e^{-\lambda^{-1}|T|^{-1}}\big),
\end{equation*}
so that condition \eqref{CHK} is equivalent to 
$$
\Big[\lambda\mapsto\frac{ \tau\big(|T|e^{-\lambda^{-1}|T|^{-1}}\big)}{\psi(\lambda)}\Big]\in 
L^\infty(\mathbb R_+^*).
$$
But  this observation is of little help for one of our motivations, that  is to derive formulas for Dixmier traces from
short-time expansion of the  heat-kernel type functions $\xi(T,t)=
\tau(e^{-t|T|^{-1}})$. 
With this regard, the correct way to understand  condition \eqref{CHK} goes as follows.
First write
$$
\frac1{\psi(\lambda)}\int_0^\lambda \tau\big(e^{-u^{-1}|T|^{-1}}\big) \,\frac{du}{u^2}
=\frac1{\psi(\lambda)}\int_0^\lambda \frac{\tau\big(e^{-u^{-1}|T|^{-1}}\big)}{u^2\psi'(u)} \,\psi'(u)\,du.
$$
Then, set $M$ for the Ces\`aro mean of the additive group 
$(\mathbb R,+)$:
$$
M:L^\infty(\R_+^*)\to L^\infty(\R_+^*),\quad f\mapsto \Big[\lambda\mapsto
\frac1{\lambda}\int_0^\lambda f(u) \,du\Big].
$$
and  observe that the map
$$
M_\psi:L^\infty(\R_+^*)\to L^\infty(\R_+^*),\quad f\mapsto \Big[\lambda\mapsto
\frac1{\psi(\lambda)}\int_0^\lambda f(u) \,\psi'(u)\,du\Big],
$$
coincides with the Ces\`aro mean twisted by $\psi$:
$$
M_\psi (f)=\big( M_+(f\circ \psi^{-1})\big)\circ\psi,\quad \forall f\in L^\infty(\R_+^*).
$$
Hence, we see that a sufficient condition to have $T\in\clC_\psi$ is that
$$
\Big[\lambda\mapsto\frac{\tau\big(e^{-\lambda^{-1}|T|^{-1}}\big)}{\lambda^2\psi'(\lambda)} \Big]
\in L^\infty(\mathbb R_+^*).
$$
For instance,  the bound 
$$
{\tau\big(e^{-t|T|^{-1}}\big)} \leq C t^{-1}\,|\log t|^{\beta},\quad \beta\geq -1,
$$
implies that $T\in\mathcal C_\psi$  for $\psi(t)=\log(1+t^{1/(1+\beta)})^{1+\beta}$, when $\beta>-1$.
For $\beta=-1$ it implies  that $T\in\mathcal C _\psi$  for
  $\psi(t)=\log(1+\log(1+t))$.

\quad

\noindent
For latter use, we first observe the following evident property:
\begin{lemma}
\label{lautre2}
Let $\psi\in\Omega$, $ T\in\clC_\psi^+$ and $S\in \clL^q$ for some $q\in(1,\infty)$. Then, we have
$$
\lim_{\lambda\to\infty}
\frac{1}{\psi(\lambda)}\int_1^\lambda\tau\big(S\,e^{-t^{-1}T^{-1}}
\big)\frac{dt}{t^2}=0.
$$
\end{lemma}
\noindent
\begin{proof}
Let $p=(1-q^{-1})^{-1}$.
The H\"older estimate on the von Neumann algebra $L^\infty([1,\lambda])\otimes \clN$ for the 
trace
$ \int_1^\lambda .\frac{dt}{t^{2}}\otimes\tau$ gives:
\begin{align*}
& \Big|\int_1^\lambda\tau\big(S\,e^{-t^{-1}T^{-1}}
\big)\frac{dt}{t^2}\Big|\leq \|S\|_q\Big(\int_1^\lambda t^{-2}\,dt\Big)^{1/q}
\Big(\int_1^\lambda\tau\big(e^{-pt^{-1} T^{-1}}\big) t^{-2}\,dt\Big)^{1/p}.
\end{align*}
Since  $T\in\clC_\psi^+$, we have
$$
\sup_{\lambda\geq 1}\frac{1}{\psi(\lambda)^{1/p}}\Big(\int_1^\lambda\tau\big(e^{-pt^{-1} T^{-1}}\big) t^{-2}\,dt\Big)^{1/p}\leq
\sup_{\lambda\geq 1}\frac{1}{\psi(\lambda)^{1/p}}\Big(\int_0^\lambda\tau\big(e^{-pt^{-1} T^{-1}}\big) t^{-2}\,dt\Big)^{1/p}<\infty.
$$
Last, we observe that
$$
\lim_{\lambda\to\infty}\frac{1}{\psi(\lambda)^{1/q}}\Big(\int_1^\lambda t^{-2}\,dt\Big)^{1/q}=0,
$$
which concludes the proof.
\hfill\end{proof}

\noindent
The two following  results  show that we  
have $\mathfrak L_\psi\subset\clC_\psi\subset\clM_\psi$, with almost no 
condition on $\psi\in\Omega$.
The first one is a direct generalization of   \cite[Proposition 3.8]{CGRS1}.
\begin{prop}
\label{C-M}
Let $T\in\clN$. Then we have
\begin{align}
\label{GGGG}
\int_0^t\mu(T,s) \,ds\leq \|T\|\,t\,\chi_{(0,1)}(t)+
\Big({\int_0^{t}\tau\big(e^{-\lambda^{-1}|T|^{-1}}\big)\,\frac{d\lambda}{\lambda^2}}+1\Big)
\chi_{[1,\infty)}(t).
\end{align}
Consequently, for
 $\psi\in\Omega$, we have $\clC_\psi\cap\clN\subset\clM_\psi$ with:
$$
\|T\|_{\clM_\psi}\leq\|T\|\sup_{t\in(0,1)}\frac t{\psi(t)}
+\sup_{t\geq 1}\frac1{\psi(t)}\Big({\int_0^{t}\tau\big(e^{-\lambda^{-1}|T|^{-1}}\big)\,\frac{d\lambda}
{\lambda^2}}+1\Big)\,,
\quad\forall T\in \clC_\psi\cap\clN.
$$
\end{prop}
\noindent
\begin{proof}
We may assume that $T\geq 0$.  For $t\in(0,1)$, we can use the estimates
$\int_0^t\mu(T,s) \,ds\leq \|T\|\,t$.
 For $t\geq 1$, 
 let
$$
A_t:=\int_0^{1/t}e^{-sT^{-1}}\,ds\quad \mbox{and}\quad B_t:=\int_{1/t}^\infty e^{-sT^{-1}}\,ds,
$$
so that by Laplace transform, we have $T=A_t+B_t$ for all $t\geq 1$.
Note first that
\begin{align*}
\|B_t\|_1=\tau(B_t)=\int_{1/t}^\infty\tau\big(e^{-sT^{-1}}\big)\,ds=
\int_0^{t}\tau\big(e^{-\lambda^{-1}T^{-1}}\big)\,\frac{d\lambda}{\lambda^2}.
\end{align*}
On the other hand, we have 
$$
\|A_t\|\leq \int_0^{1/t}\big\|e^{-sT^{-1}}\big\|\,ds\leq 1/t.
$$
Thus for $t>0$, we have
\begin{align*}
\int_0^t\mu(T,s)\,ds&=\inf\big\{\|B\|_1+t\|A\| \,:\,T=B+A\in\clL^1+\clN\big\}\\
&\leq \|B_t\|_1+t\|A_t\|\leq \int_0^{t}\tau\big(e^{-\lambda^{-1}T^{-1}}\big)\,\frac{d\lambda}{\lambda^2}+1.
\end{align*}
This concludes the proof. 
\hfill\end{proof}

\begin{prop}
\label{cont}
Let $\psi\in\Omega$ be such that $\clM_\psi$ is close to $\clL^1$, then
$\mathfrak L_\psi\subset \clC_\psi$.
\end{prop}
\noindent
\begin{proof}
We may assume that $T\geq0$. Recall that $\mathfrak L_\psi\subset\clN$ by Remark \ref{rem2}. We have for $\eps>0$:
\begin{align*}
&\frac1{\psi(e^{\eps^{-1}})}\int_0^{e^{\eps^{-1}}}
\tau\big(e^{-\lambda^{-1} T^{-1}}\big)\,\frac{d\lambda}{\lambda^2}
=\frac1{\psi(e^{\eps^{-1}})}\int_{e^{-\eps^{-1}}}^\infty\tau\big(e^{-t T^{-1}}\big)\,dt
=\frac1{\psi(e^{\eps^{-1}})}\,\tau\Big(\int_{e^{-\eps^{-1}}}^\infty e^{-t T^{-1}}\,dt\Big)\\
&\qquad\qquad=\frac{\tau\big(Te^{-e^{-\eps^{-1}}T^{-1}}\big)}{\psi(e^{\eps^{-1}})}
\leq\frac{\|T\|_{1+\eps}^{1+\eps}}{\psi(e^{\eps^{-1}})}\,\big\|T^{-\eps}e^{-e^{-\eps^{-1}}T^{-1}}\big\|
\leq \|T\|_{\mathfrak L_\psi}^{1+\eps}\,\psi(e^{\eps^{-1}})^\eps\,\big\|T^{-\eps}e^{-e^{-\eps^{-1}}T^{-1}}\big\|.
\end{align*}
By functional calculus, we easily deduce
$$
\big\|T^{-\eps}e^{-e^{-\eps^{-1}}T^{-1}}\big\|= e\,\eps^\eps\,e^{-\eps}\leq e\,,\quad \eps\in(0,1].
$$
Since $\psi(t)\leq Ct$ for $t$ large enough, we get
 $$
\sup_{\eps\in(0,1]} \frac1{\psi(e^{\eps^{-1}})}\int_0^{e^{\eps^{-1}}}
\tau\big(e^{-\lambda^{-1} T^{-1}}\big)\,\frac{d\lambda}{\lambda^2}\leq C\,
\|T\|_{\mathfrak L_\psi}^2,
 $$
 which concludes the proof.
\hfill\end{proof}

\noindent
In conclusion, we observe that for a Lorentz space close to $\clL^1$, one always has
$$
\mathfrak L_\psi\subset \clC_\psi\subset \clM_\psi\subset\mathfrak L^\psi\;,
$$
with continuous inclusions. Our next goal is to show that under condition \eqref{SSZ-cond}
(hence under the stronger condition \eqref{our-cond} too),
we have the equality $ \clC_\psi= \clM_\psi$. We will prove a more general
result, similar to \cite[Theorem 40]{SZ}:
\begin{prop}
\label{equal-f}
Let $\psi\in\Omega$ satisfying  condition \eqref{SSZ-cond}. Then, for all 
 $f\in L^\infty(\mathbb R_+^*)\cap C^2(\mathbb R_+^*)$ such that $f'(0)=f(0)=0$, for all 
 $T\in\clM_\psi^+$ and all $B\in\clN$, we have
$$
\Big[t\mapsto \frac{1}{\psi(t)}\int_0^t\tau\big(Bf(uT)\big) \,\frac{du}{u^2}\Big]\in L^\infty(\mathbb R_+^*).
$$
\end{prop}
\noindent
\begin{proof}
By assumptions on the function $f$, we have 
$|f(t)|\leq c\min\{1,t^2\}=ct^2\chi_{[0,1)}(t)+c\chi_{[1,\infty)}(t)$. Hence,
\begin{align}
\label{first}
\big|\tau\big(Bf(uT)\big)\big|\leq c\|B\|\Big(u^2\tau\big(T^2\chi_{[0,1/u)}(T)\big)+
\tau\big(\chi_{[1/u,\infty)}(T)\big)\Big).
\end{align}
For the second term in the inequality above, we get with $E_T(\lambda)$ the spectral projection of 
$T$:
\begin{align*}
\int_0^t\tau\big(\chi_{[1/u,\infty)}(T)\big) \,\frac{du}{u^2}=
\tau\Big(\int_0^t\Big(\int_0^\infty\chi_{[1/u,\infty)}(\lambda)\,dE_T(\lambda)\Big) \,\frac{du}{u^2}\Big).
\end{align*}
But
$$
\int_0^t\chi_{[1/u,\infty)}(\lambda)\,\frac{du}{u^2}=\int_0^t\chi_{[1/\lambda,\infty)}(u)\,\frac{du}{u^2}=
\int_{1/\lambda}^t\frac{du}{u^2} \;\chi_{[1/t,\infty)}(\lambda)\leq
\lambda\,\chi_{[1/t,\infty)}(\lambda),
$$
and thus by Fubini's Theorem, we get
\begin{align*}
\int_0^t\tau\big(\chi_{[1/u,\infty)}(T)\big) \,\frac{du}{u^2}\leq
\tau\Big(\int_0^\infty\lambda\,\chi_{[1/t,\infty)}(\lambda)\,dE_T(\lambda)\Big)=
\tau\big(T\,\chi_{[1/t,\infty)}(T)\big)
=\int_0^{n(T,{1/t})}\mu(T,s)\,ds,
\end{align*}
where $n(T)$ is the distribution function of $T$. Now, by \cite[Lemma 7]{SSZ}
(this is where condition \eqref{SSZ-cond} is used),
for all $\eps>0$, there exists $t_0>0$ such that for all $t\geq t_0$, we have
$n(T,{1/t})\leq(1+\eps)\|T\|_{\clM_\psi}\,t\,\psi(t)$. Hence, with $\eps>0$ fixed, 
we get for $t$ large enough
\begin{align*}
&\frac1{\psi(t)}\int_0^t\tau\big(\chi_{[1/u,\infty)}(T)\big) \,\frac{du}{u^2}\leq
\frac1{\psi(t)}\int_0^{(1+\eps)\|T\|_{\clM_\psi}\,t\,\psi(t)}\mu(T,s)\,ds
\leq
\frac{\psi\big((1+\eps)\|T\|_{\clM_\psi}\,t\,\psi(t)\big)}{\psi(t)}\,\|T\|_{\clM_\psi}.
\end{align*}
Observe then that the  right most  side of the inequality above is bounded. Indeed,   condition 
\eqref{SSZ-cond}
implies condition \eqref{ex-sing-trace-easy} and thus
$$
\frac{\psi\big((1+\eps)\|T\|_{\clM_\psi}\,t\,\psi(t)\big)}{\psi(t)}=
\frac{\psi\big((1+\eps)\|T\|_{\clM_\psi}\,t\,\psi(t)\big)}{\psi\big(t\psi(t)\big)}
\frac{\psi\big(t\psi(t)\big)}{\psi(t)}\to1\,,\quad t\to\infty.
$$
For the contribution corresponding to the first term in the right hand side of \eqref{first}, we get
\begin{align*}
\int_0^t\tau\big(T^2\chi_{[1/u,\infty)}(T)\big) \,du=\tau\Big(\int_0^\infty\lambda^2\,
\Big(\int_0^t\chi_{[0,1/u)}(\lambda)\,du\Big)\,dE_T(\lambda)\Big).
\end{align*}
Since
$$
\int_0^t\chi_{[0,1/u)}(\lambda)\,du=\min\{t,1/\lambda\}=t\chi_{[0,1/t)}(\lambda)+\lambda^{-1}
\chi_{[1/t,\infty)}(\lambda),
$$
we deduce
$$
\int_0^t\tau\big(T^2\chi_{[1/u,\infty)}(T)\big) \,du=t\,\tau\big(T^2\,\chi_{[0,1/t)}(T)\big)+
\tau\big(T\,\chi_{[1/t,\infty)}(T)\big).
$$
The conclusion for the contribution coming from the second term above is reached using
the same argument as above. Thus, all that remains to show  is that
$$
\tau\big(T^2\,\chi_{[0,1/t)}(T)\big)=O\big(\psi(t)/t\big).
$$
To this end, we write:
\begin{align*}
\tau\big(T^2\,\chi_{[0,1/t)}(T)\big)&=\int_0^\infty \mu(T,s)^2\,\chi_{\{\mu(T)\leq1/t\}}(s)\,ds\\
&=
\int_0^{\|T\|_{\clM_\psi} t\psi(t)} \mu(T,s)^2\,\chi_{\{\mu(T)\leq1/t\}}(s)\,ds+\int_{\|T\|_{\clM_\psi} t\psi(t)}^\infty
 \mu(T,s)^2\,\chi_{\{\mu(T)\leq1/t\}}(s)\,ds\\
 &\leq \|T\|_{\clM_\psi}\frac{\psi(t)}t+\int_{\|T\|_{\clM_\psi} t\psi(t)}^\infty \mu(T,s)^2\,ds.
\end{align*}
Next, we use
$$
\|T\|_{\clM_\psi}\geq\frac{1}{\psi(t)}\int_0^t\mu(T,s)\,ds\geq\frac{t\mu(T,t)}{\psi(t)},
$$
to get for all $\delta\in(0,1)$:
\begin{align*}
&\int_{\|T\|_{\clM_\psi} t\psi(t)}^\infty \mu(T,s)^2\,ds\leq\|T\|_{\clM_\psi}^2\int_{\|T\|_{\clM_\psi} t\psi(t)}^\infty s^{-2}\,
\psi(s)^2\,ds
\\&
\leq  \|T\|_{\clM_\psi}^2\, \psi\big(\|T\|_{\clM_\psi} \,t\,\psi(t)\big)^2\,(\|T\|_{\clM_\psi} \,t\,\psi(t))^{-\delta}\,
\int_{\|T\|_{\clM_\psi} t\psi(t)}^\infty s^{-2+\delta}\,ds = \frac{\|T\|_{\clM_\psi}}{1-\delta}\, 
\frac{\psi\big(\|T\|_{\clM_\psi} \,t\,\psi(t)\big)^2}{t\psi(t)}.
\end{align*}
Letting $\delta\to 0$, we eventually get:
$$
\tau\big(T^2\,\chi_{[0,1/t)}(T)\big)\leq \|T\|_{\clM_\psi}\frac{\psi(t)}t\,\Big(1+ 
\frac{\psi\big(\|T\|_{\clM_\psi}\, t\,\psi(t)\big)^2}{\psi(t)^2}\Big),
$$
which concludes the proof since the right most fraction above is bounded as already shown.
\hfill\end{proof}

\begin{cor}
\label{equal}
Let $\psi\in\Omega$ satisfying  condition \eqref{SSZ-cond}. Then we have $ \clC_\psi= \clM_\psi$.
\end{cor}
\noindent
\begin{proof}
By Proposition \ref{C-M}, we only need to prove that if $T\in\clM_\psi^+$ then 
$$
\sup_{\lambda>0}
\frac1{\psi(\lambda)}\int_0^\lambda \tau\big(e^{-t^{-1}T^{-1}}\big) \,\frac{dt}{t^2}<\infty.
$$
But the result follows from Proposition \ref{equal-f} applied to the function  $f(x)=e^{-x^{-1}}$
and to $B=1$.
\hfill\end{proof}

\noindent
We can now state the main result of this section, which relates Dixmier traces and heat kernels,
when the strongest condition  \eqref{our-cond} holds.

 \begin{thm}
\label{psi2}
Let $\psi\in\Omega$ satisfying  condition \eqref{our-cond}.
Then, for any exponentiation invariant state   $\omega$ of $ L^\infty(\mathbb R_+^*)$,  any 
 $ T\in \mathfrak L_\psi^+$ ($\subset\clC_\psi^+=\clM_\psi^+$ in this 
 case) and any $B\in\clN$,  we have
$$
\tau_{\psi,\omega}(BT)=
\omega\Big(\Big[\lambda\mapsto \frac1{\psi(\lambda)}\int_1^\lambda \tau\big(Be^{-t^{-1}T^{-1}}\big)
 \,\frac{dt}{t^2}\Big]\Big).
 $$
\end{thm}
\noindent
\begin{proof}
By linearity, we may assume $B\in\clN^+$ and by traciality and Theorem \ref{zeta},
 we are left to prove that:
$$
\omega\Big(\Big[\lambda\mapsto \frac1{\psi(\lambda)}\int_0^\lambda \tau\big(B^{1/2}
e^{-t^{-1}T^{-1}}B^{1/2}\big) \,\frac{dt}{t^2}\Big]\Big)
=\frac1{\Gamma(1+k_\psi)}\,\,\omega\circ\log\Big(\Big[r\mapsto\frac1{\psi(e^r)}
{\tau(B^{1/2}T^{1+1/r}B^{1/2})}\Big]\Big),
$$
with $k_\psi$ the constant associated to $\psi\in\Omega$,  given as in Definition \ref{kpsi}.
Using the integral representation $T^{1+1/r}=\Gamma(1+1/r)^{-1}\int_0^\infty t^{1/r}e^{-tT^{-1}}dt$
and since $\Gamma(1+1/r)\to1$, $r\to\infty$,
we deduce
\begin{align*}
\omega\circ\log\Big(\Big[r\mapsto\frac1{\psi(e^r)}{\tau(B^{1/2}T^{1+1/r}B^{1/2})}\Big]\Big)=
\omega\circ\log\Big(\Big[r\mapsto \frac1{\psi(e^r)}
{\int_0^\infty t^{1/r}\tau(B^{1/2}e^{-tT^{-1}}B^{1/2})}\,dt\Big]\Big).
\end{align*}
 Next, observe that for $t>1$, we have since $\mathfrak L_\psi\subset\clN$:
$$
B^{1/2}e^{-tT^{-1}}B^{1/2}\leq e^{-(t-1)/\|T\|}\, B^{1/2}e^{-T^{-1}}B^{1/2},
$$
which implies that
$$
\int_1^\infty t^{1/r}\tau(B^{1/2}e^{-tT^{-1}}B^{1/2})\,dt\leq \tau(B^{1/2}e^{-T^{-1}}B^{1/2})\,
\int_1^\infty t\, e^{-(t-1)/\|T\|}\,dt,
$$
for $r\geq 1$. Hence, 
$$
\lim_{r\to \infty} \frac1{\psi(e^r)}{\int_1^\infty t^{1/r}\tau(B^{1/2}e^{-tT^{-1}}B^{1/2})}\,dt
=0,
$$
and thus
\begin{align*}
\omega\circ\log\Big(\Big[r\mapsto \frac1{\psi(e^r)}{\tau(B^{1/2}T^{1+1/r}B^{1/2})}\Big]\Big)=
\omega\circ\log\Big(\Big[r\mapsto \frac1{\psi(e^r)}{\int_0^1 t^{1/r}\tau(B^{1/2}e^{-tT^{-1}}B^{1/2})}\,dt
\Big]\Big).
\end{align*}
Set 
$$
\beta(u):=\int_0^ue^{-\nu}\tau(B^{1/2}e^{-e^{-\nu}T^{-1}}B^{1/2})\,d\nu.
$$
Then, $\beta$ is real valued, increasing  and right continuous  on $\mathbb R_+^*$ 
 with $\beta(0)=0$ and a little computation shows that
 $$
 \int_0^\infty e^{-\mu/r}\,d\beta(\mu)=\int_0^1 t^{1/r}\tau(B^{1/2}e^{-tT^{-1}}B^{1/2})\,dt.
 $$
 Therefore, we get
 \begin{align*}
\omega\circ\log\Big(\Big[r\mapsto\frac1{\psi(e^r)}{\tau(B^{1/2}T^{1+1/r}B^{1/2})}\Big]\Big)=
\omega\circ\log\Big(\Big[r\mapsto\frac1{\psi(e^r)}{ \int_0^\infty e^{-\mu/r}\,d\beta(\mu)}\Big]\Big).
\end{align*}
Since $T\in\mathfrak L_\psi^+$, we also get from what precedes: 
$$
\Big[r\mapsto \frac1{\psi(e^r)}{ \int_0^\infty e^{-\mu/r}\,d\beta(\mu)}\Big]\in L^\infty(\mathbb R_+^*).
$$
Thus, we can apply Proposition \ref{lem2}, to deduce:
 \begin{align*}
 \omega\circ\log\Big(\Big[r\mapsto \frac1{\psi(e^r)}{ \int_0^\infty e^{-\mu/r}\,d\beta(\mu)}\Big]\Big)&=
\Gamma(1+k_\psi)\,\, \omega\circ\log\Big(\Big[r\mapsto \frac{ \beta(r)}{\psi(e^r)}\Big]\Big)
\\& =\Gamma(1+k_\psi)\,\,
 \omega\Big(\Big[r\mapsto\frac{ \beta(\log(r))}{\psi(r)}\Big]\Big).
 \end{align*}
 But
 \begin{align*}
\frac{ \beta(\log(r))}{\psi(r)}=\frac1{\psi(r)}
\int_0^{\log(r)}e^{-\nu}\tau(B^{1/2}e^{-e^{-\nu}T^{-1}}B^{1/2})\,d\nu=
\frac1{\psi(r)}\int_1^r\tau(B^{1/2}e^{-\lambda^{-1}T^{-1}}B^{1/2})\,\frac{d\lambda}{\lambda^2},
 \end{align*}
 and finally,
  \begin{align*}
\frac1{\Gamma(1+k_\psi)}\,\omega\circ\log\Big(\Big[r\mapsto
 \frac{\tau(B^{1/2}T^{1+1/r}B^{1/2})}{\psi(e^r)}\Big]\Big)&=
 \omega\Big(\Big[r\mapsto \frac1{\psi(r)}\int_1^r\tau(B^{1/2}
 e^{-\lambda^{-1}T^{-1}}B^{1/2})\,\frac{d\lambda}{\lambda^2}\Big]\Big),
 \end{align*}
which completes the proof.
\hfill\end{proof}

\noindent
Last, we apply the previous result in the case of a  short-time asymptotic expansion
of the heat trace function with logarithms in the leading term.

\begin{cor}
\label{OUOU}
Let $T\in\clN^+$ and $\beta\in[-1,\infty)$ be such that there exists $a(T)>0$ with 
\begin{align}
\label{ASS}
\xi(T,t):={\tau\big(e^{-tT^{-1}}\big)} \sim a (T)\,t^{-1}\,|\log t|^{\beta},\quad t\downarrow 0.
\end{align}
Then  $T\in\mathcal \clM_\psi$  for $\psi(t)=\log(1+t^{1/(1+\beta)})^{1+\beta}$ when $\beta>-1$ and 
for  $\psi(t)=\log(1+\log(1+t))$ when $\beta=-1$. Moreover, for any  exponentiation invariant state 
$\omega$  on $ L^\infty(\mathbb R_+^*)$, we have
$$
\tau_{\psi,\omega}(T)=C(\beta)\,a(T)\quad\mbox{with}\quad
C(\beta)=\begin{cases}
(1+\beta)^\beta, &\beta>-1,\\
1,&\beta=-1.
\end{cases}
$$
\end{cor}
\noindent
\begin{proof}
Note first that since $T\in\clN^+$, the map  
$\big[t\in\R_+^*\mapsto{\tau\big(e^{-tT^{-1}}\big)}\big]$ is decreasing and thus the behavior 
\eqref{ASS} gives rise to a  bound ${\tau\big(e^{-tT^{-1}}\big)} \leq C\,t^{-1}\,|\log t|^{\beta}$
for all $t>0$ and some constant $C>0$. Then, the discussion right after Definition \ref{C} combined
with Proposition \ref{equal-f} shows that $T\in\clM_\psi$, for the elements $\psi\in\Omega_b$ 
given as above. Then Theorem \ref{psi2} gives the result for $\mathfrak L_\psi^+$,  the  subcone
 of $\clN^+$. But Proposition \ref{prop-eq} implies that in these cases, we have 
$\mathfrak L_\psi=\mathcal \clM_\psi$, concluding the proof.
\hfill\end{proof}

\noindent
One may wonder about the interest of Corollary \ref{OUOU} knowing that even in the most exotic 
(commutative) geometric situations (e.g$.$ manifolds with boundaries or with singularities),
the $\log$-terms in the asymptotic expansion (near zero) of the trace of the  heat semigroup 
generated by  an elliptic pseudo-differential operator, never show up in the leading divergence
(see for instance \cite{Lesch} and the references therein).  
However, this is no longer the case for noncommutative spaces, where
logarithms may already appear in the leading divergence of the trace of the heat semigroup. 
Indeed, this happens  for the Podl\`es spheres, a quantum homogeneous
space for the quantum group $SU_q(2)$.
 Corollary \ref{OUOU} thus
 produces  examples of (noncommutative) geometric 
  operators living in a Lorentz space
different from $\clM^{1,\infty}$.
\begin{example}
The authors of \cite{EIS,KM},  consider 
 a positive operator $P$ on the Podl\`es sphere (with deformation parameter $0<q<1$)
 such that 
\begin{align}
\label{QP}
{\rm Tr}\big(e^{-t P}\big)\sim_{t\to 0} A(q)\,|\log(t)|^{2},
\end{align}
where $\rm Tr$ is the standard operator trace.
In \cite{EIS}, $P$ is the absolute value of a Dirac-type operator and $A(q)=2/(\log q)^2$, while 
in \cite{KM}, $P$ is a $q$-analogue of the 
Casimir-type operator and $A(q)=1/4(\log q)^2$. The natural question
we address here,  is to know whether the constant
$A(q)$ can be recovered as the value of a Dixmier trace of an operator an operator
related to $P$.  
Of course, this might seem hopeless since the equivalence
 \eqref{QP}  implies that $(1+P)^{-1}$ is of trace class. This 
(well known) fact can be easily deduced from  the general inequality \eqref{GGGG}.
However, if one considers the operator $P\otimes 1+1\otimes \Delta_{\mathbb T^2}$ instead,
where $\Delta_{\mathbb T^2}$ denotes the ordinary Laplacian on  the $2$-torus\footnote{One can 
take any 2-dimensional compact Riemannian manifold instead of $\mathbb T^2$.}
$\mathbb T^2=\R^2/\mathbb Z^2$, then from \cite[p. 307]{Chavel}, 
$$
{\rm Tr}\big(e^{-t \Delta_{\mathbb T^2}}\big)= \frac{{\rm Vol}(\mathbb T^2)}{4\pi t} \sum_{k\in\mathbb Z^2}e^{-| k|^2/4t}
\sim_{t\to 0}\frac{{\rm Vol}(\mathbb T^2)}{4\pi t} ,
$$
 we get
\begin{align*}
{\rm Tr}\big(e^{-t (P\otimes 1+1\otimes \Delta_{\mathbb T^2})}\big)={\rm Tr}\big(e^{-t P}\big)
{\rm Tr}\big(e^{-t \Delta_{\mathbb T^2})}\big)\sim_{t\to 0} \frac{A(q){\rm Vol}(\mathbb T^2)}{4\pi} \frac{|\log(t)|^{2}}{t}.
\end{align*}
Hence, Corollary \ref{OUOU} implies that the operator 
$ (P\otimes 1+1\otimes \Delta_{\mathbb T^2})^{-1}$ belongs to the Lorentz space 
(of ordinary compact operators on a separable infinite dimensional Hilbert space) associated
with $\psi(t)=\log(1+t^{1/3})^3$
and that for any  exponentiation invariant state 
$\omega$  on $ L^\infty(\mathbb R_+^*)$, we have
$$
A(q)=\frac{4\pi}{9{\rm Vol}(\mathbb T^2)}\, {\rm Tr}_{\psi,\omega}\big( (P\otimes 1+1\otimes \Delta_{\mathbb T^2})^{-1}\big)
\,,\qquad \psi(t)=\log(1+t^{1/3})^3.
$$
\end{example}

\section{Applications to pseudo-differential operators on $\mathbb R^n$}
\label{psido}
The interconnections between the theory of Dixmier traces and that of pseudo-differential operators are well known.
Without pretending to be exhaustive,  these interconnections include Connes Trace Theorem \cite{Co0, KLPS, LPS, LSZ}, 
Toeplitz and Hankel operators \cite{ER}, 
noncompact isospectral deformations \cite{G1,GIV}, pseudo-differential operators with 
$\log$-polyhomogeneous symbols \cite{Lesch}, pseudo-differential operators on manifolds
with boundaries \cite{FGLS,NS}, holomorphic families of pseudo-differential operators \cite{N}, 
H\"ormander-Weyl pseudo-differential
operators \cite{NR}, pseudo-differential operators on manifolds
with conical singularities \cite{S}, pseudo-differential
operators on $CR$- and contact manifolds \cite{Po}. However, and with the exception of \cite{N},
in all these works only the ideal $\clM^{1,\infty}$ is concerned. 

\quad

In this last section, our main goal in to relate Dixmier traces and pseudo-differential operators
for more general Lorentz ideal.
More specifically, we apply here our results to the setting of H\"ormander-Weyl pseudo-differential
operators on $\R^n$. 
The framework that we will introduce in a moment,  is strongly related to the work of
Nicola and Rodino in \cite{NR}. In fact,
the first main result of this section, Theorem \ref{Nique}, generalizes their main result (which is stated 
as Theorem 1.1 
there) in two directions. First,  it works for more general  Lorentz spaces $\clM_\psi$  than the 
dual to the Macaev ideal $\clM^{1,\infty}$. Second, 
 it is not restricted to the non-closed subspace characterized by $\mu(T)=O(\psi')$. But we 
  go even further: we are able to express the Dixmier trace of an 
 H\"ormander-Weyl pseudo-differential operator in term of its symbol only. This is done
 in Theorem \ref{ils} and Corollary \ref{mam} and we believe that they are important results.
 
 \quad

\noindent
Recall that the Weyl quantization map, ${\rm OP}_W$, is a
 continuous linear map from the space of tempered distributions on $\R^{2n}$ (with the strong
dual topology) to the space of continuous linear maps acting from the space of Schwartz functions on
 $\R^n$ 
(with its standard Fr\'echet topology) to the space of tempered distributions on $\R^{n}$:
$$
{\rm OP}_W\in\clL\left(\clS'(\mathbb R^{2n}),\clL\big( \clS(\mathbb R^{n}), \clS'(\mathbb R^{n})\big)
\right),
$$
given, with a little abuse of notations, by
$$
\big({\rm OP}_W(T) \phi\big)(x)=(2\pi)^{-n}\int_{\mathbb R^{2n}}
e^{i\xi.(x-y)}\,T\big(\frac{x+y}2,\xi\big)\,\phi(y)\,d^ny\,d^n\xi,
\quad T\in \clS'(\mathbb R^{2n}),\;\phi\in \clS(\mathbb R^{n}).
$$
For $T\in\mathcal S'(\mathbb R^{2n})$, the linear operator 
${\rm OP}_W(T)$ from $\clS(\mathbb R^{n})$ to $\clS'(\mathbb R^{n})$ is called the Weyl
pseudo-differential operator with symbol $T$. The most important properties of this
quantization scheme, not shared by other variants of pseudo-differential calculus
on $\R^n$, is that it maps real symbols to self-adjoint operators: 
$$
{\rm OP}_W(\overline T)={\rm OP}_W(T)^*,
$$
and defines a unitary operator from the Hilbert 
  space of $L^2$-function on $\R^{2n}$, to the Hilbert space of Hilbert-Schmidt operators acting on
  $L^2$-functions on $\R^n$:
  $$
  {\rm OP}_W\in\mathcal U\left(L^2(\R^{2n}), \clL^2\big(L^2(\R^n)\big)\right),
  $$
that is to say
  \begin{align}
  \label{unit}
    {\rm Tr}\big(
 {\rm OP}_W(f_1)^*{\rm OP}_W(f_2)\big)=
 \int_{\R^{2n}} \overline{ f_1(x,\xi)}\,f_2(x,\xi)\,d^nx\,d^n\xi,
 \qquad \forall f_1,f_2\in L^2(\R^{2n}).
 \end{align}
It is then easy to see that  the relation \eqref{unit}  still holds for
a pair of symbols $(f_1,f_2)\in L^\infty(\R^{2n})\times L^1(\R^{2n})$ 
 such that ${\rm OP}_W(f_1)$ is bounded and
${\rm OP}_W(f_2)$ is trace class (and vice verca).

\quad

\noindent
We consider here symbols in
 $S(m,g)$, the  H\"ormander space. This class of symbols is
associated to $g$, a slowly varying and $\sigma$-temperate metric in $\mathbb R^{2n}$,
satisfying the uncertainty principle 
and to $m:\mathbb R^{2n}\to \mathbb R_+^*$, a $g$-continuous and $(\sigma,g)$-temperate
weight function, see \cite[Definition 18.5.1]{Hor}. 
The most important feature of $S(m,g)$ is that it is
a Fr\'echet space for the   topology  associated
to the seminorms:
\begin{align}
\label{Smg}
\|f\|_{k;m,g}:=
\sup_{(x,\xi)\in\mathbb R^{2n}}\;\sup_{X_1,\dots,X_k\in \,T_{x,\xi}\mathbb R^{2n}}\;
\frac{\big|f^{(k)}\big((x,\xi);X_1,\dots,X_k\big)\big|}{m(x,\xi)\,g_{x,\xi}(X_1)^{1/2}\dots g_{x,\xi}(X_k)^{1/2}}
<\infty, \quad  k\in\mathbb N_0.
\end{align}
Here,  $f^{(k)}\big((x,\xi);.,\dots,.\big)$ denotes
 the $k$-multilinear form on  
$T_{(x,\xi)}\mathbb R^{2n}$, the tangent space of $\R^{2n}$ at $(x,\xi)$,
 given by the differential of order $k\in\mathbb N_0$ of the function
$f$ at the point $(x,\xi)\in\mathbb R^{2n}$.
Accordingly,
we let ${\rm OP}S(m,g)$ to be the class of Weyl pseudo-differential operators with
symbols  in the H\"ormander class $S(m,g)$. 
We call the elements of ${\rm OP}S(m,g)$  the H\"ormander-Weyl pseudo-differential operators.
We also define  $g^\sigma$, the symplectic dual metric of $g$, by 
$$
g^\sigma_{x,\xi}(t,\tau):=\sup\big\{\sigma(t,\tau;y,\eta)^2\;:\;g_{x,\xi}(y,\eta)=1\big\},\qquad
(x,\xi),(t,\tau)\in\mathbb R^{2n},
$$
where $\sigma$ is the standard symplectic form of $\mathbb R^{2n}\simeq T^*\mathbb R^{n}$.
We also let $h_g$ be the so called Planck function. It is defined by: 
\begin{align}
\label{cond-NRh}
h^2_g(x,\xi):=\sup_{(t,\tau)\in\mathbb R^{2n}}\frac{g_{x,\xi}(t,\tau)}{g_{x,\xi}^\sigma(t,\tau)}.
\end{align}
The uncertainty principle mentioned above, corresponds to the condition that
 $h_g\leq1$. Here, we make the further assumption that there exists $q\in(1,\infty)$ 
 such that $h_g\in L^q(\R^{2n})$. Observe however that this condition is weaker than
 the one used in \cite{NR}, which is $h_g(x,\xi)\leq C(1+|x|+|\xi|)^{-\delta}$, for some 
 $C,\delta>0$.
 
 \quad

We refer to  \cite[Chapter XVIII]{Hor} for precise definitions and details as here 
we just need to know the following facts. First,
${\rm OP}_W$ sends $S(m,g)$ to $\mathcal B\big(L^2(\R^n)\big)$  continuously, if and
only if $m$ is bounded \cite[Theorem 18.6.3]{Hor}. Second,
${\rm OP}_W$ sends $S(m,g)$ to $\mathcal K\big(L^2(\R^n)\big)$  continuously, if and
only if $m\to 0$ at $\infty$  \cite[Theorem 18.6.6]{Hor}. If moreover $m\in S(m,g)$ then
 ${\rm OP}_W(m)$ is self-adjoint and its spectrum
 is bounded from below and if  $m\to\infty$ at $\infty$,  then the spectrum
  is discrete \cite[Theorem 3.4]{Hor2}.   Following the discussion of \cite[page 143]{Hor}, given
 a $g$-continuous and $(\sigma,g)$-temperate
weight function $m$, we can always find another  $g$-continuous and $(\sigma,g)$-temperate
weight function $m'$ such that $m'\in S(m',g)$ and $S(m',g)=S(m,g)$. Hence, the assumption
$m\in S(m,g)$ (that we will make in what follows)  is, in fact,  irrelevant. 
We will also use the following trace norm estimate, proved in 
   \cite[Theorem 3.9]{Hor2}: for all $k\in\N$, there exists $C_k>0$ such that for all $f\in S(m,g)$,
   we have
   \begin{align}
   \label{trest}
   \|{\rm OP}_W(f)\|_1\leq C_k(\|f\|_1+\|h_g^km\|_1\|f\|_{k;m,g}).
   \end{align}
   Set $\star$ for the composition product associated to the Weyl calculus, i.e$.$ defined by
   the relation
   $$
   {\rm OP}_W(f_1){\rm OP}_W(f_2)={\rm OP}_W(f_1\star f_2).
   $$
 The last property  we need to know, is that the operation  $\star$ defines a continuous bilinear
  map:
  $$
  \star:S(m_1,g)\times S(m_2,g)\to S(m_1m_2,g),
  $$
   and if $f_j\in S(m_j,g)$, $j=1,2$, we have
  \begin{align}
  \label{comp}
  f_1\star f_2-f_1f_2\in S(m_1m_2h_g,g).
  \end{align}
  This is proven in \cite[Theorem 18.5.4]{Hor}.
  To simplify the discussion, it is useful to introduce the following terminology:
  \begin{defn}
  \label{HP}
  An  H\"ormander pair $(g,m)$ consists of a slowly varying and $\sigma$-temperate metric $g$ in 
  $\mathbb R^{2n}$ such that the Planck function $h_g$, given in \eqref{cond-NRh}, belongs to $L^q(\R^{2n})\cap L^\infty(\R^{2n})$
  for some $q\in(1,\infty)$ and 
a weight function  $m:\mathbb R^{2n}\to \mathbb R_+^*$ which is $g$-continuous, $(\sigma,g)$-temperate
and which satisfies $m\in S(m,g)$.
  \end{defn}
  
Finally, we need to record some  results of \cite{NR}.
\begin{prop}
\label{NRprop}
Let $(g,m)$  be an H\"ormander pair such that $m\to 0$ at $\infty$. \\
(i) There exists $c>0$, such that  setting $\tilde m:=(c+m^{-1})^{-1}$, the 
operator ${\rm OP}_W(\tilde m^{-1})$ is positive, boundedly invertible and with compact inverse 
in ${\rm OP}S(m,g)$.\\
(ii) For all $t\geq 0$, we have
$$
e^{-t\,{\rm OP}_W(\tilde m^{-1})}={\rm OP}_W(b_t)+S_t,
$$
where $\{b_t\}_{t>0}$ is a  bounded family of symbols in $S(1,g)$ 
  and $\{S_t\}_{t>0}$ is a family of trace class operators such that
  we have $\|S_t\|_1\leq C\,t$, for some $C>0$ and all $t>0$. \\
(iii) For fixed
$t>0$,  $b_t\in\cap_{l\in \N}S(m^l,g)$\\
(iv) The symbol $b_t$ can be written as
 \begin{equation*}
  b_t=e^{-t\tilde m^{-1}}+\sum_{j=1}^Nb_{t,j}\quad\mbox{with}\quad
  |b_{t,j}|\leq C_j\, e^{-t\tilde m/2}\,h_g^j,
  \end{equation*} 
  for  some constants $C_j>0$ independent of $t>0$. In particular, we have
  for some $C>0$:
  $$
  | b_t|\leq C e^{-t\tilde m/2}.
  $$
\end{prop}
\noindent
\begin{proof}
(i) is a restatement of  \cite[Lemma 3.2]{NR}.  (ii)  is \cite[Theorem 3.3]{NR}. 
(iii) and the first part of (iv) follow from the intermediate estimates used in the proof of 
\cite[Theorem 3.3]{NR} and the last part of (iv) follows from $h_g\in L^\infty(\R^{2n})$, 
by assumption.
\end{proof}

\noindent
The following result is mostly a consequence of Proposition
 \ref{C-M}.
The Lorentz spaces we consider there are associated with the type $I_\infty$
  factor
 of all  bounded operators on $L^2(\mathbb R^n)$ and  with the standard trace.
 To make the notations explicit, we denote this Lorentz space by $\clM_\psi(L^2(\R^n))$.  
We stress that in this type $I_\infty$ factor setting, there is no need to distinguish $\Omega$ and 
$\Omega_b$.
\begin{prop}
\label{elle}
Let $\psi\in\Omega_b$
and let $(g,m)$  be an H\"ormander pair.  If furthermore
\begin{align}
\label{cond-NR}
\sup_{\lambda>0}\frac{1}{\psi(\lambda)}\int_0^\lambda \int_{\R^{2n}} e^{-t^{-1}
 m^{-1}(x,\xi)}
\,d^nx\,d^n\xi\frac{dt}{t^2}<\infty,
\end{align}
then ${\rm OP}S(m,g)\subset\clM_\psi(L^2(\R^n))$ continuously.
\end{prop}
\noindent
\begin{proof}
We start by observing that the condition \eqref{cond-NR} implies that 
$m\to0$ at $\infty$.
Indeed assume it does not. Then, in the proof of \cite[Lemma 3.1]{NR} and under the 
assumption that $m$ is $g$-admissible, the authors have constructed a  Borel set
$B$ of $\mathbb R^{2n}$ and of infinite Lebesgue measure ($B=\cup_{j=1}^\infty B_j$
within the notations of \cite[Lemma 3.1]{NR})
and such that $B\subset\{m> C\}$, for a constant $C>0$. This entails that 
$$
\int_{\R^{2n}} e^{-t m^{-1}(x,\xi)}
\,d^nx\,d^n\xi\geq \int_B e^{-t m^{-1}(x,\xi)}\,d^nx\,d^n\xi\geq  \int_B e^{-t C^{-1}}
\,d^nx\,d^n\xi=e^{-t C^{-1}}{\rm meas}(B)=\infty,
$$
a contradiction with the assumption \eqref{cond-NR}. 
Thus, we can apply Proposition \ref{NRprop}. We adopt the following notations:
\begin{align}
\label{bob}
\tilde m_\delta:=(c+\delta+m^{-1})^{-1}, \quad \mbox{for} \quad
\delta>0,\qquad \mbox{and}\qquad\tilde m:=\tilde m_0.
\end{align}
Next, we show that \eqref{cond-NR} 
implies that
 ${\rm OP}_W(\tilde m^{-1}_\delta)^{-1}$ 
belongs to $\clM_\psi(L^2(\R^n))$ for all $\delta>0$.
Indeed, since  ${\rm OP}_W(\tilde m^{-1}_\delta)^{-1}$ is bounded and positive by Proposition 
\ref{NRprop} (i),
we get from the first estimate of Proposition \ref{C-M} and
 for  $t>0$:
\begin{align}
\label{fgh}
&\int_0^t\mu\big({\rm OP}_W(\tilde m^{-1}_\delta)^{-1},s\big)\,ds \nonumber\\
&\qquad\qquad\qquad\leq
 \|{\rm OP}_W(\tilde m^{-1}_\delta)^{-1}\|\,t\,\chi_{(0,1)}(t) +\Big(
\int_0^t {\rm Tr}\big(e^{-u^{-1}\,{\rm OP}_W(\tilde m^{-1}_\delta)}\big)
\,\frac{du}{u^2}+1
\Big)\chi_{[1,\infty)}(t).
\end{align} 
Writing,
${\rm OP}_W(\tilde m^{-1}_\delta)={\rm OP}_W(\tilde m^{-1})+\delta$, we deduce by
 Proposition \ref{NRprop} (ii):
\begin{align*}
\int_0^\lambda {\rm Tr}\big(e^{-u^{-1}\,{\rm OP}_W(\tilde m^{-1}_\delta)}\big)
\,\frac{du}{u^2}&\leq \int_0^\lambda \Big(\|{\rm OP}_W(b_{u^{-1}})\|_1
+\|S_{u^{-1}}\|_1\Big)\,e^{-u^{-1}\delta}\,\frac{du}{u^2}\\
&\leq \int_0^\lambda \Big(\|{\rm OP}_W(b_{u^{-1}})\|_1
+\frac Cu\Big)\,e^{-u^{-1}\delta}
\,\frac{du}{u^2}.
\end{align*}
The last term is  bounded uniformly  in $\lambda>0$ 
 and since $\psi(t)=O(t)$, $t\to 0$ (by the 
assumption that $\psi\in\Omega_b$), the corresponding contribution to \eqref{fgh} is 
bounded.
For the
first term, we use  the inequality given in \eqref{trest} to get that for every $N\in\mathbb N$, there exists 
$k\in\mathbb N$ such that
we have
$$
 \|{\rm OP}_W(b_{u^{-1}})\|_1\leq C_N\big(\|b_{u^{-1}}\|_1+\|h_g^N\|_1\|b_{u^{-1}}\|_{k;1,g}\big).
$$
By assumption, we have $\|h^N\|_1<\infty$ for $N\geq q$ 
 and we  know by Proposition \ref{NRprop} (ii) that $[u\mapsto \|b_{u^{-1}}\|_{k;1,g}]$ belongs to
$ L^\infty(\R^*_+)$ for each $k\in\mathbb N$. Thus, the corresponding contribution in \eqref{fgh}
is finite.
 For the contribution coming from $ \|b_{u^{-1}}\|_1$,  Proposition \ref{NRprop} (iv)
 shows that
 $$
 \|b_{u^{-1}}\|_1\leq C \int_{\R^{2n}} e^{-u^{-1}
 m^{-1}(x,\xi)/2}
\,d^nx\,d^n\xi,
 $$
 and we conclude using
 our condition \eqref{cond-NR}.
 Now, take any symbol $f\in S(m,g)$ then we write
 $$
 {\rm OP}_W(f)= {\rm OP}_W(f)\,{\rm OP}_W(\tilde m^{-1}_\delta)\, {\rm OP}_W(\tilde m^{-1}_
 \delta)^{-1}.
 $$
As $\tilde m^{-1}_\delta\in S(m^{-1},g)$ we have by \cite[Theorem 18.5.4]{Hor} that 
${\rm OP}_W(f)\,{\rm OP}_W(m^{-1}_\delta)\in {\rm OP}S(1,g)$ and we conclude using the fact
that ${\rm OP}S(1,g)$ in contained in the set of bounded operators on $L^2(\mathbb R^n)$ by
\cite[Theorem 18.6.3]{Hor}.
This is enough to conclude as $\clM_\psi(L^2(\R^n))$ is an ideal of $\clB\big(L^2(\mathbb R^n)\big)$.
\hfill\end{proof}

\noindent
Observe (see Definition \ref{C}) that the finiteness of the quantity \eqref{cond-NR}
just means that the weight function $m$ belongs to the space $\clC_\psi$ for the 
commutative von Neumann algebra $L^\infty(\mathbb R^{2n})$ with Lebesgue integral
for semi-finite normal trace. Also,  when $\psi\in\Omega_b$ satisfies condition \eqref{SSZ-cond},
Corollary \ref{equal} shows that the latter space coincides with the commutative Lorentz
space $\clM_\psi(\R^{2n})$.
Hence,   Proposition \ref{elle} yields the following `existence result':
\begin{thm}
\label{Nique}
Let $\psi\in\Omega_b$ satisfying condition \eqref{SSZ-cond} and 
$(g,m)$   an H\"ormander pair.
If $m\in  \clM_\psi(\mathbb R^{2n})$, then
 ${\rm OP}S(m,g)\subset\clM_\psi\big(L^2(\mathbb R^n)\big)$ continuously.
\end{thm}

Our final aim is to prove that  the Dixmier
trace of an H\"ormander-Weyl pseudo-differential operator coincides with the Dixmier
trace of its symbol, for the commutative von Neumann algebra of essentially bounded functions
on $\mathbb R^{2n}$ with trace given by the Lebesgue integral. This result nicely complement
the classical fact that, when trace class, the trace  of a Weyl pseudo-differential operator 
coincide with the Lebesgue integral of its symbol\footnote{This follows from the fact that, for a 
kernel operator on $L^2(\R^n)$ of trace class, its trace coincides with the integral of its kernel on the diagonal.}.
 This result relies on the
properties of the H\"ormander-Weyl calculus and on the
following simple lemma.

\begin{lemma}
\label{lautre}
Let $\psi\in\Omega_b$ satisfying condition \eqref{SSZ-cond} and 
$(g,m)$   an H\"ormander pair.
Let also $0<f\in\clM_\psi(\R^{2n})$. Then, we have
$$
\lim_{\lambda\to\infty}
\frac{1}{\psi(\lambda)}\int_1^\lambda \Big(\int_{\R^{2n}} 
h_g(x,\xi)\,e^{-t^{-1} f^{-1}(x,\xi)}
\,d^nx\,d^n\xi\Big)\frac{dt}{t^2}=0.
$$
\end{lemma}
\noindent
\begin{proof}
Since $\psi$ satisfies condition \eqref{SSZ-cond}, we deduce by
Corollary \ref{equal} that $\clC_\psi(\R^{2n})= \clM_\psi(\R^{2n})$. By assumption 
that $(g,m)$ is  an H\"ormander pair, we have $h_g\in L^q(\R^{2n})$ for some $q\in(1,\infty)$. 
Thus, the proof follows directly from Lemma \ref{lautre2}.
\hfill\end{proof}

We can now formulate the second main result of this section.  Here,  $\int_{\psi,\omega}$ denotes a 
Dixmier trace for the commutative Lorentz space $\clM_\psi(\R^{2n})$ associated with 
$L^\infty(\R^{2n})$ and with the Lebesgue integral. For this result to hold, condition \eqref{SSZ-cond} 
is not enough and we need to consider the stronger condition \eqref{our-cond}.

\begin{thm}
\label{ils}
Let $\psi\in\Omega_b$ satisfying  condition \eqref{our-cond} and assume further that 
$\mathfrak L_\psi=\clM_\psi$.
Let also $(g,m)$  be an  H\"ormander pair such that 
$m\in\clM_\psi(\mathbb R^{2n})$.
Then, for any symbol $f\in S(m,g)$ and any exponentiation invariant state $\omega$
of $L^\infty(\R_+^*)$, we have
$$
{\rm Tr}_{\psi,\omega}\big({\rm OP}_W(f)\big)=\int_{\psi,\omega} f.
$$
\end{thm}
\noindent
\begin{proof}
First, fix $\delta>0$ and write as in Proposition \ref{elle} and with the notations given in \eqref{bob}:
\begin{align}
\label{eux}
 {\rm OP}_W(f)= {\rm OP}_W(f)\,{\rm OP}_W(\tilde m_\delta^{-1})\, {\rm OP}_W(\tilde m_\delta^{-1})^{-1}
 = {\rm OP}_W(f\star \tilde m_\delta^{-1})\, {\rm OP}_W(\tilde m_\delta^{-1})^{-1}.
\end{align}
We have shown in the proof of Proposition \ref{elle},  that
$$
{\rm OP}_W(\tilde m_\delta^{-1})^{-1}\in\clM_\psi\big(L^2(\R^n)\big)^+
=\mathfrak L_\psi\big(L^2(\R^n)\big)^+.
$$
Since $f\in S(m,g)$ and $\tilde m_\delta^{-1}\in  S(m^{-1},g)$, we deduce by
\eqref{comp} that $f\star\tilde m_\delta^{-1}\in  S(1,g)$
and by \cite[Theorem 18.6.3]{Hor}  that
$$
 {\rm OP}_W(f\star\tilde m_\delta^{-1})\in 
 \mathcal B\big(L^2(\R^n)\big),
 $$
hence, using \eqref{eux},  we get  ${\rm OP}_W(f)\in\clM_\psi\big(L^2(\R^n)\big)$.
Thus, we can apply Theorem \ref{psi2}, to get:
$$
{\rm Tr}_{\psi,\omega}\big( {\rm OP}_W(f)\big)=
\omega\Big(\Big[\lambda\mapsto \frac1{\psi(\lambda)}\int_1^\lambda {\rm Tr}\big(
 {\rm OP}_W(f\star\tilde m_\delta^{-1})e^{-t^{-1}{\rm OP}_W(\tilde m_\delta^{-1})}\big)
 \,\frac{dt}{t^2}\Big]\Big).
 $$
Next,  by Proposition \ref{NRprop} (ii) we have 
$$
e^{-t^{-1}\,{\rm OP}_W(\tilde m_\delta^{-1})}=
e^{-t^{-1}\,{\rm OP}_W(\tilde m^{-1})}e^{-\delta/t}
=\big({\rm OP}_W(b_{t^{-1}})+S_{t^{-1}}\big)e^{-\delta/t}\quad\mbox{with}\quad
 \|S_{t^{-1}}\|_1\leq Ct^{-1}.
$$
Hence,  as
$$
\big| {\rm Tr}\big(
 {\rm OP}_W(f\star\tilde m_\delta^{-1})\,S_{t^{-1}}\big)\big|\leq
  \big\|{\rm OP}_W(f\star\tilde m_\delta^{-1})\big\|
\, \|S_{t^{-1}}\|_1,
$$
we get
$$
\Big|\frac1{\psi(\lambda)}\int_1^\lambda {\rm Tr}\big(
 {\rm OP}_W(f\star\tilde m_\delta^{-1})\,S_{t^{-1}}\big)
 \,\frac{dt}{t^2}\Big|
 \leq C\big\|{\rm OP}_W(f\star\tilde m_\delta^{-1})\big\|
 \frac1{\psi(\lambda)}\int_1^\lambda 
 \,\frac{dt}{t^3}\to0,\quad \lambda\to\infty.
 $$
 Accordingly, we get
\begin{align*}
{\rm Tr}_{\psi,\omega}\big( {\rm OP}_W(f)\big)&=
\omega\Big(\Big[\lambda\mapsto \frac1{\psi(\lambda)}\int_1^\lambda {\rm Tr}\big(
 {\rm OP}_W(f\star\tilde m_\delta^{-1}){\rm OP}_W(b_{t^{-1}})\big)e^{-\delta/t}
 \,\frac{dt}{t^2}\Big]\Big).
 \end{align*}
 Since $f\star\tilde m_\delta^{-1}\in L^\infty(\R^{2n})$, $b_{t^{-1}}\in L^1(\R^{2n})$, 
   ${\rm OP}_W(f\star\tilde m_\delta^{-1})$ is bounded and ${\rm OP}_W(b_{t^{-1}})$ is trace
 class, we can employ  the  relation \eqref{unit}, to get   
\begin{align*}
  {\rm Tr}\big(
 {\rm OP}_W(f\star\tilde m_\delta^{-1}){\rm OP}_W( b_{t^{-1}})\big)&=
\int_{\R^{2n}}
f\star\tilde m_\delta^{-1}(x,\xi)\, b_{t^{-1}}(x,\xi)\,d^nx d^n\xi.
 \end{align*}
 Consequently
  \begin{align*}
{\rm Tr}_{\psi,\omega}\big( {\rm OP}_W(f)\big)&=
\omega\Big(\Big[\lambda\mapsto \frac1{\psi(\lambda)}\int_1^\lambda\int_{\R^{2n}}
f\star\tilde m_\delta^{-1}(x,\xi)\, b_{t^{-1}}(x,\xi)\,d^nx d^n\xi\,e^{-\delta/t}
 \,\frac{dt}{t^2}\Big]\Big).
 \end{align*}
But by \cite[Theorem 18.5.4]{Hor}, we have 
 $f\star \tilde m_\delta^{-1}-f\tilde m_\delta^{-1}\in S(h_g,g)$, which, in particular, gives 
 $$
 |f\star \tilde m_\delta^{-1}-f\tilde m_\delta^{-1}|\leq C h_g.
 $$
 Hence, since
$|b_{t^{-1}}|\leq C e^{-t^{-1}\tilde m/2}$ by Proposition \ref{NRprop} (iv),
 since $\tilde m\in \clM_\psi(\R^{2n})$, we can use Lemma \ref{lautre} to deduce that
  $$
  \omega\Big(\Big[\lambda\mapsto \frac1{\psi(\lambda)}\int_1^\lambda\int_{\R^{2n}}
\big(f\star\tilde m_\delta^{-1}(x,\xi)-f(x,\xi)\tilde m_\delta^{-1}(x,\xi)\big)\, b_{t^{-1}}(x,\xi)\,d^nx d^n\xi
\,e^{-\delta/t} \,\frac{dt}{t^2}\Big]\Big)=0
  $$
  Accordingly, we get
  \begin{align*}
{\rm Tr}_{\psi,\omega}\big( {\rm OP}_W(f)\big)=
\omega\Big(\Big[\lambda\mapsto \frac1{\psi(\lambda)}\int_1^\lambda \int_{\R^{2n}}
f(x,\xi)\;\tilde m_\delta^{-1}(x,\xi)\;b_{t^{-1}}(x,\xi)\,d^nx d^n\xi\,e^{-\delta/t}
 \,\frac{dt}{t^2}\Big]\Big).
 \end{align*}
Similarly, we have by Proposition \ref{NRprop} (iv) that
$$
|b_{t^{-1}}-e^{-{t^{-1}}\tilde m^{-1}}|\leq C\, h_g\,e^{-{t^{-1}}\tilde m^{-1}/2}.
$$
Hence, from Lemma \ref{lautre} again,  we deduce that 
\begin{align*}
{\rm Tr}_{\psi,\omega}\big( {\rm OP}_W(f)\big)&=
\omega\Big(\Big[\lambda\mapsto \frac1{\psi(\lambda)}\int_1^\lambda \int_{\R^{2n}}
f(x,\xi)\;\tilde m_\delta^{-1}(x,\xi)\;e^{-{t^{-1}}\tilde m_\delta^{-1}(x,\xi)}\,d^nx d^n\xi
 \,\frac{dt}{t^2}\Big]\Big).
  \end{align*}
However, since $f\,\tilde m_\delta^{-1}\in S(1,g)\subset L^\infty(\R^{2n})^+$  and since 
 $\tilde m_\delta\in \clM_\psi(\R^{2n})^+=\mathfrak L_\psi(\R^{2n})^+$, we can use 
 Theorem \ref{psi2} for the von Neumann algebra $L^\infty(\R^{2n})$ with Lebesgue integral, to 
 deduce 
 $$
 \omega\Big(\Big[\lambda\mapsto \frac1{\psi(\lambda)}\int_1^\lambda\int_{\R^{2n}}
f(x,\xi)\;\tilde m_\delta^{-1}(x,\xi)\;e^{-t\tilde m_\delta^{-1}(x,\xi)}\,d^nx d^n\xi
 \,\frac{dt}{t^2}\Big]\Big)
 =\int_{\psi,\omega}(f\,\tilde m_\delta^{-1}\,\tilde m_\delta)=
 \int_{\psi,\omega}f.
 $$
This completes the proof.
\hfill\end{proof}

Our last result combines Theorem \ref{zeta} and  Theorem \ref{ils} 
(with $B={\rm Sign}(f)$ and $T=|f|$). It complements 
 Theorem \ref{ils}  as
it gives a very simple way to compute the Dixmier trace of an H\"ormander-Weyl
pseudo-differential operator. 
\begin{corollary}
\label{mam}
Let $\psi\in\Omega_b$ satisfying  condition \eqref{our-cond} and assume further that 
$\mathfrak L_\psi=\clM_\psi$.
Let also $(g,m)$ be an H\"ormander pair such that $m\in  \clM_\psi(\mathbb R^{2n})$.
Then, for any  symbol $ f\in S(m,g)$ and any exponentiation invariant state $\omega$
of $L^\infty(\R_+^*)$, we have:
$$
{\rm Tr}_{\psi,\omega}\big({\rm OP}_W(f)\big)=\frac1{\Gamma(1+k_\psi)}\,
\omega\Big(\Big[r\mapsto\frac1{\psi(r)}\int_{\R^{2n}}f(x,\xi) \;|f(x,\xi)|^{1/\log(r)}\;d^nxd^n\xi\Big]\Big),
$$
where $k_\psi$ is the constant associated to $\psi$ given in Definition \ref{kpsi}.
\end{corollary}

\begin{rem}
The results of this section certainly extend to  more general contexts. In particular,
to pseudo-differential operators on Riemannian manifolds of bounded geometry with symbols
in  Shubin's classes  \cite{Shu4}.
\end{rem}

\end{document}